\documentclass[10pt]{amsart} 
\usepackage[all]{xy}
\usepackage{float}
\usepackage{tikz-cd}
\usepackage{verbatim}
\usepackage{mathtools}
\usepackage{ stmaryrd }
\usepackage{amsmath}
\usepackage{amssymb}
\usetikzlibrary{intersections, matrix,decorations.pathreplacing}
\usetikzlibrary{calc}
\usepackage{longtable}

\tikzset{
    ncbar angle/.initial=90,
    ncbar/.style={
        to path=(\tikztostart)
        -- ($(\tikztostart)!#1!\pgfkeysvalueof{/tikz/ncbar angle}:(\tikztotarget)$)
        -- ($(\tikztotarget)!($(\tikztostart)!#1!\pgfkeysvalueof{/tikz/ncbar angle}:(\tikztotarget)$)!\pgfkeysvalueof{/tikz/ncbar angle}:(\tikztostart)$)
        -- (\tikztotarget)
    },
    ncbar/.default=0.5cm,
}

\tikzset{square left brace/.style={ncbar=0.5cm}}
\tikzset{square right brace/.style={ncbar=-0.5cm}}

\tikzset{round left paren/.style={ncbar=0.5cm,out=120,in=-120}}
\tikzset{round right paren/.style={ncbar=0.5cm,out=60,in=-60}}

\newtheorem{theorem}{Theorem}[section]
\newtheorem{lemma}[theorem]{Lemma}
\newtheorem{cor}[theorem]{Corollary}
\newtheorem*{lemma*}{Lemma}

\newtheorem{prop}[theorem]{Proposition}
\newtheorem*{prop*}{Proposition}

\theoremstyle{definition}

\newtheorem{definition}[theorem]{Definition}

\theoremstyle{remark}
\newtheorem{remark}[theorem]{Remark}

\numberwithin{equation}{section}

\newcommand{\A}{\mathcal{A}}
\newcommand{\B}{\mathcal{B}}
\newcommand{\C}{\mathcal{C}}
\newcommand{\OmegaMat}{\Omega_{\boxplus}}
\DeclareFontFamily{U}{wncy}{}
    \DeclareFontShape{U}{wncy}{m}{n}{<->wncyr10}{}
    \DeclareSymbolFont{mcy}{U}{wncy}{m}{n}
    \DeclareMathSymbol{\Sh}{\mathord}{mcy}{"58} 

\newcommand\restr[2]{{% we make the whole thing an ordinary symbol
  \left.\kern-\nulldelimiterspace % automatically resize the bar with \right
  #1 % the function
  \vphantom{\big|} % pretend it's a little taller at normal size
  \right|_{#2} % this is the delimiter
  }}
\allowdisplaybreaks

\begin{document}

\title{The Zigzag Hochschild Complex \\  \today }

\author[C.~Miller]{Cheyne ~Miller}
  \address{Cheyne Miller, Department of Mathematics, The Graduate Center, City University of New York, 365 Fifth Ave, NY 10016}
  \email{cmiller@gradcenter.cuny.edu}
\begin{abstract}
In this paper, a new higher Hochschild Complex is defined with an Iterated Integral map to locally model differential forms on the space of bigons on $M$. In particular, given the local data for a gerbe with structure 2-group given by a crossed module of matrix-groups, there is an element in our curved zigzag Hochschild complex associated to the local 2-holonomy given by such a gerbe. This paper introduces an initial construction central to the author's PhD Thesis. 
\end{abstract}

\maketitle

\tableofcontents

\section{The zigzag Hochschild Model: An Illustrated Introduction} 
The motivation for this paper is born out of modeling iterated integrals \cite{C} in a non-commutative setting.  While it is not an issue to simply compute an iterated integral in a non-abelian setting, finding the correct algebraic structure which respects the product and differential in that setting is a problem of interest to not only the study of 2-holonomy for non-abelian gerbes, but also other fields such as, for example, in the study of Quantum Control Theory and the study of Multiple Dedekind Zeta Values.  While the latter two fields may find some tools in this paper useful, we will focus our attention on differential forms on smooth spaces with values in some not-necessarily-commutative-algebra -- so that we generally restrict ourselves to Matrix Lie Algebras -- which we will denote $\OmegaMat(M)$ for a manifold, $M$.

In a commutative setting, suppose we would like to consider the wedge product of two iterated integrals (each thought of as a differential form on the path space $PM$ of a manifold $M$) in  $\Omega(PM)$.  Then we are able to write 
\begin{align*}
&\int\limits_{\Delta^n} a_1(\tau_1) \dots a_n(\tau_n) d\tau_1 \dots d\tau_n \wedge \int\limits_{\Delta^m} a_{n+1}(\tau_{n+1}) \dots a_{n+m}(\tau_{n+m}) d\tau_{n+1} \dots d\tau_{n+m} \\
= &\sum\limits_{\sigma} \pm \int\limits_{\Delta^{n+m}} a_{\sigma^{-1}(1)}(\tau_{\sigma^{-1}(1)}) a_{\sigma^{-1}(n+m)}(\tau_{\sigma^{-1}(n+m)}) d \tau_{\sigma^{-1}(1)}  \dots  d \tau_{\sigma^{-1}(n+m)}
\end{align*} where the sum is over order-preserving shuffles, $\sigma$.  The picture corresponding to this wedge product is the illustration for the shuffle product on the interval Hochschild complex:  
\begin{itemize}
\item[]
\resizebox{3cm}{!}{\begin{tikzpicture}
\draw[thick] (0,0) -- coordinate[pos=.13] (a1)  coordinate[pos=.32] (a2) coordinate[pos=.48] (a3) coordinate[pos=.68] (a4) coordinate[pos=.88] (a5)(2,0);
\foreach \x in {1,2,...,5}
{  \draw ($(a\x)+(0,5pt)$) -- ($(a\x)-(0,5pt)$);
\node[above= 0.1cm] at (a\x) {\small$a_{\x}$};}
\end{tikzpicture}} $\odot$  \resizebox{3cm}{!}{\begin{tikzpicture}
\draw[thick] (0,0) -- coordinate[pos=.13] (b1)  coordinate[pos=.32] (b2) coordinate[pos=.48] (b3) coordinate[pos=.68] (b4) (2,0);
\foreach \x in {1,2,...,4}
{  \draw ($(b\x)+(0,5pt)$) -- ($(b\x)-(0,5pt)$);
\node[above= 0.1cm] at (b\x) {\small$b_{\x}$};}
\end{tikzpicture}} 
\item[] $= \sum\limits_{\text{shuffles}}$ \resizebox{6cm}{!}{\begin{tikzpicture}
\draw[thick] (0,0) -- coordinate[pos=.13] (x1) coordinate[pos=.20] (x2) coordinate[pos=.32] (x3) coordinate[pos=.40] (x4)coordinate[pos=.48] (x5) coordinate[pos=.56] (x6) coordinate[pos=.68] (x7) coordinate[pos=.75] (x8)  coordinate[pos=.90] (x9) (5,0);
\foreach \x in {1,3, 5, 6, 9}
{  \draw ($(x\x)+(0,5pt)$) -- ($(x\x)-(0,5pt)$);}
\node[above= 0.1cm] at (x1) {\small$a_{1}$};
\node[above= 0.1cm] at (x3) {\small$a_{2}$};
\node[above= 0.1cm] at (x5) {\small$a_{3}$};
\node[above= 0.1cm] at (x6) {\small$a_{4}$};
\node[above= 0.1cm] at (x9) {\small$a_{5}$};

\foreach \x in {2,4, 7, 8}
{  \draw ($(x\x)+(0,5pt)$) -- ($(x\x)-(0,5pt)$);}
\node[above= 0.1cm] at (x2) {\small$b_{1}$};
\node[above= 0.1cm] at (x4) {\small$b_{2}$};
\node[above= 0.1cm] at (x7) {\small$b_{3}$};
\node[above= 0.1cm] at (x8) {\small$b_{4}$};
\end{tikzpicture}} 
\end{itemize}
However, when we are in the non-commutative setting, the best we can do when writing the wedge product out is:
\begin{align*}
&\int\limits_{\Delta^n} a_1(\tau_1) \dots a_n(\tau_n) d\tau_1 \dots d\tau_n \wedge \int\limits_{\Delta^m} b_1(\tau_{n+1}) \dots b_m(\tau_{n+m}) d\tau_{n+1} \dots d\tau_{n+m} \\
= &\sum\limits_{\sigma} \pm \int\limits_{\Delta^{n+m}} a_1(\tau_{\sigma^{-1}(1)})\dots a_n(\tau_{\sigma^{-1}(n)})b_1(\tau_{\sigma^{-1}(n+1)}) \dots b_m(\tau_{\sigma^{-1}(n+m)}) d\tau_{1} \dots d\tau_{n+m}
\end{align*} 
We would like an illustration for this wedge product in some Hochschild complex.  Note that the order of the differential forms is preserved but the coordinates are shuffled around.  Moreover, by integrating all of the $a$'s we have passed over the interval (from left to right) once and then we pass over the interval a second time while integrating all of the $b$'s.  The corresponding picture below is our shuffle product in the zigzag Hochschild complex which preserves the order of the differential forms but shuffles the coordinates:
\begin{itemize}
\setlength{\itemsep}{10pt}
\item[]
\resizebox{3cm}{!}{\begin{tikzpicture}
\draw[thick] (0,0) -- coordinate[pos=.13] (a1)  coordinate[pos=.32] (a2) coordinate[pos=.48] (a3) coordinate[pos=.68] (a4) coordinate[pos=.88] (a5)(2,0);
\foreach \x in {1,2,...,5}
{  \draw ($(a\x)+(0,5pt)$) -- ($(a\x)-(0,5pt)$);
\node[above= 0.1cm] at (a\x) {\small$a_{\x}$};}
\end{tikzpicture}} $\odot$  \resizebox{3cm}{!}{\begin{tikzpicture}
\draw[thick] (0,0) -- coordinate[pos=.13] (b1)  coordinate[pos=.32] (b2) coordinate[pos=.48] (b3) coordinate[pos=.68] (b4) (2,0);
\foreach \x in {1,2,...,4}
{  \draw ($(b\x)+(0,5pt)$) -- ($(b\x)-(0,5pt)$);
\node[above= 0.1cm] at (b\x) {\small$b_{\x}$};}
\end{tikzpicture}} 
\item[]  $= \sum\limits_{\text{shuffles}}$ \resizebox{!}{3cm}{\begin{tikzpicture}[baseline={([yshift=-.5ex]current bounding box.center)},vertex/.style={anchor=base,
    circle,fill=black!25,minimum size=18pt,inner sep=2pt}]
\draw[thick] (0,10)  --coordinate[pos=.13] (x1) coordinate[pos=.20] (x2) coordinate[pos=.32] (x3) coordinate[pos=.40] (x4)coordinate[pos=.48] (x5) coordinate[pos=.56] (x6) coordinate[pos=.68] (x7) coordinate[pos=.75] (x8)  coordinate[pos=.90] (x9) 
(12,8) -- (0,6) -- coordinate[pos=.13] (y1) coordinate[pos=.20] (y2) coordinate[pos=.32] (y3)
 coordinate[pos=.40] (y4)coordinate[pos=.48] (y5) coordinate[pos=.56] (y6) coordinate[pos=.68] (y7) coordinate[pos=.75] (y8)  coordinate[pos=.90] (y9) (12,4) -- (0,2);
\foreach \x in {1,...,9}
{  \draw ($(x\x)+(0,5pt)$) -- ($(x\x)-(0,5pt)$);
    \draw ($(y\x)+(0,5pt)$) -- ($(y\x)-(0,5pt)$);}

\node[above= 0.1cm] at (x1) {\large$a_{1}$};
\node[above= 0.1cm] at (x3) {\large$a_{2}$};
\node[above= 0.1cm] at (x5) {\large$a_{3}$};
\node[above= 0.1cm] at (x6) {\large$a_{4}$};
\node[above= 0.1cm] at (x9) {\large$a_{5}$};
\node[above= 0.1cm] at (x2) {\large$1$};
\node[above= 0.1cm] at (x4) {\large$1$};
\node[above= 0.1cm] at (x7) {\large$1$};
\node[above= 0.1cm] at (x8) {\large$1$};

\node[below= 0.1cm] at (y2) {\large$b_{1}$};
\node[below= 0.1cm] at (y4) {\large$b_{2}$};
\node[below= 0.1cm] at (y7) {\large$b_{3}$};
\node[below= 0.1cm] at (y8) {\large$b_{4}$};
\node[below= 0.1cm] at (y1) {\large$1$};
\node[below= 0.1cm] at (y3) {\large$1$};
\node[below= 0.1cm] at (y5) {\large$1$};
\node[below= 0.1cm] at (y6) {\large$1$};
\node[above= 0.1cm] at (y9) {\large$1$};
\end{tikzpicture}} 
\end{itemize}
In the above diagram our zigzags, although they require some vertical space to draw, are actually passing back and forth over the same interval.  With this motivation for our zigzags in hand, Section 2.1 provides the formal details for the zigzag Hochschild complex,  $CH^{ZZ}(\A)$, of an associative (unital) algebra, $\A$.  Section 2.2 deals with certain special cases and the compatibility between the interval Hochschild complex and our zigzag Hochschild complex for these cases.  Chapter 3 deals with an iterated integral map $CH^{ZZ}(\OmegaMat(M)) \xrightarrow{It} \OmegaMat(PM)$.  

In Chapter 4, we explore the higher (2-dimensional) Hochschild complexes.   In a commutative setting it is sufficient to use a bi-simplicial model, $CH^{I\times I}(\A)$, consisting of monomials of the form

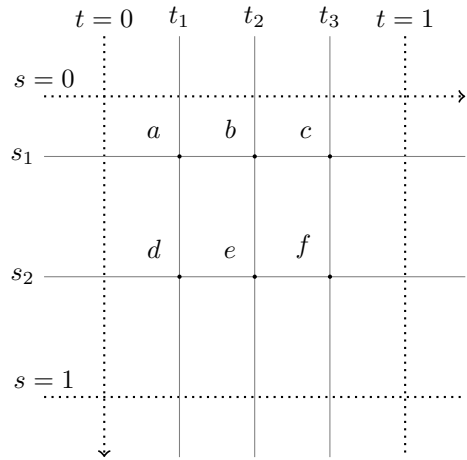
\begin{equation*}
\begin{tikzpicture}[scale=0.4]
%[thick,scale=0.6, every node/.style={transform shape}]
%help lines grid
%\draw[help lines] (0.1,0.1) grid (9.9,9.9);

% s axis
\draw[thick][dotted][->] (0,12) -- (0,-2);
\node [left] at (0,-2) {};
% t axis
\draw[thick][dotted][->] (-2,10) -- (12,10);
\node [above] at (12,10) {};
%label axes
%\node [above left] at (0,10) {0};
\node[above] at (0,12) {$t=0$};
\draw[thick][dotted] (10,12) -- (10,-2);
\node[above] at(10,12) {$t=1$};
\draw[thick][dotted] (-2,0) -- (12,0);
\node[above] at(-2,0) {$s=1$};
\node[above] at (-2,10) {$s=0$};

%s-lines
%s_1

%s_2
%label si's
\draw[name path=s1, thin, gray] (-2,8) -- (12,8);
\node[left] at (-2,8) {$s_1$};
\draw[name path=s2, thin, gray] (-2,4) -- (12,4);
\node[left] at (-2,4) {$s_2$};

%t-lines

%t1
\draw[name path=t1, thin, gray] (2.5,-2) -- (2.5,12);
\node[above] at (2.5,12) {$t_1$};

%t2
\draw[name path=t2, thin, gray] (5,-2) -- (5,12);
\node[above] at (5,12) {$t_2$};

%t3
\draw[name path=t3, thin, gray] (7.5,-2) -- (7.5,12);
\node[above] at (7.5,12) {$t_3$};

\fill [name intersections={of=t1 and s1, by={a}}]
       (a) circle (2pt);
\node[label=above left:$a$] at (a) {};

\fill [name intersections={of=t1 and s2, by={a}}]
       (a) circle (2pt);
\node[label=above left:$d$] at (a) {};

\fill [name intersections={of=t2 and s1, by={a}}]
       (a) circle (2pt);
\node[label=above left:$b$] at (a) {};

\fill [name intersections={of=t2 and s2, by={a}}]
       (a) circle (2pt);
\node[label=above left:$e$] at (a) {};

\fill [name intersections={of=t3 and s1, by={a}}]
       (a) circle (2pt);
\node[label=above left:$c$] at (a) {};

\fill [name intersections={of=t3 and s2, by={a}}]
       (a) circle (2pt);
\node[label=above left:$f$] at (a) {};

\end{tikzpicture}
\end{equation*}
Here, the variables $t_1, t_2, s_1$, and $s_2$ describe the coordinates by which a 2-dimensional iterated integral map $CH^{I\times I}(\OmegaMat(M)) \to \OmegaMat(M^{I \times I})$ has to be evaluated.  The differential for $CH^{I \times I}(A)$ consists of two components: one which  ``collapses'' vertically and one which ``collapses'' horizontally.   We see below that $D^2=0$ is not possible in the non-abelian case, where we have highlighted the significant areas of the pictures:

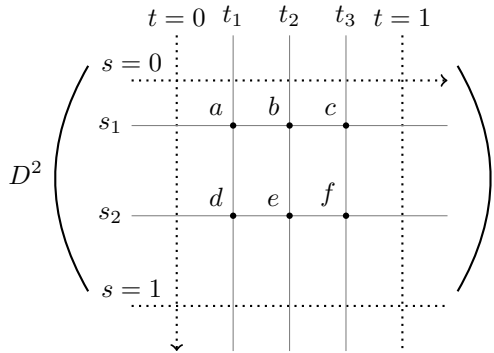
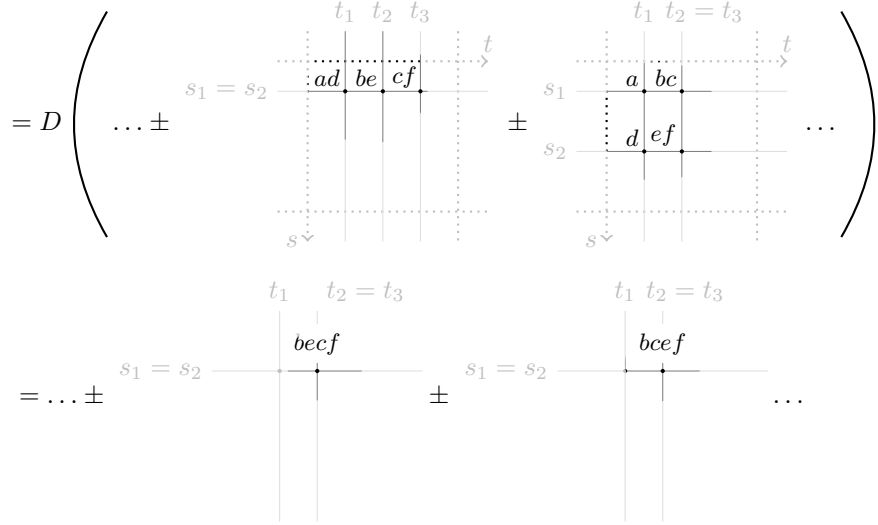
\begin{equation*} \hspace*{-4cm} 
  D^2
  \begin{tikzpicture}[scale= 0.2, baseline={([yshift=-.5ex]current bounding box.center)},vertex/.style={anchor=base,circle,fill=black!25,minimum size=18pt,inner sep=2pt}]
  \draw [thick] (0,0) to [round left paren ] (0,15);
  \end{tikzpicture}
   \begin{tikzpicture}[scale= 0.3, baseline={([yshift=-.5ex]current bounding box.center)},vertex/.style={anchor=base,circle,fill=black!25,minimum size=18pt,inner sep=2pt}]
  
%[thick,scale=0.6, every node/.style={transform shape}]
%help lines grid
%\draw[help lines] (0.1,0.1) grid (9.9,9.9);

% s axis
\draw[thick][dotted][->] (0,12) -- (0,-2);
\node [left] at (0,-2) {};
% t axis
\draw[thick][dotted][->] (-2,10) -- (12,10);
\node [above] at (12,10) {};
%label axes
%\node [above left] at (0,10) {0};
\node[above] at (0,12) {$t=0$};
\draw[thick][dotted] (10,12) -- (10,-2);
\node[above] at(10,12) {$t=1$};
\draw[thick][dotted] (-2,0) -- (12,0);
\node[above] at(-2,0) {$s=1$};
\node[above] at (-2,10) {$s=0$};

%s-lines
%s_1

%s_2
%label si's
\draw[name path=s1, thin, gray] (-2,8) -- (12,8);
\node[left] at (-2,8) {$s_1$};
\draw[name path=s2, thin, gray] (-2,4) -- (12,4);
\node[left] at (-2,4) {$s_2$};

%t-lines

%t1
\draw[name path=t1, thin, gray] (2.5,-2) -- (2.5,12);
\node[above] at (2.5,12) {$t_1$};

%t2
\draw[name path=t2, thin, gray] (5,-2) -- (5,12);
\node[above] at (5,12) {$t_2$};

%t3
\draw[name path=t3, thin, gray] (7.5,-2) -- (7.5,12);
\node[above] at (7.5,12) {$t_3$};

\fill [name intersections={of=t1 and s1, by={a}}]
       (a) circle (4pt);
\node[above left = 0.01cm] at (a) {$a$};

\fill [name intersections={of=t1 and s2, by={a}}]
       (a) circle (4pt);
\node[above left = 0.01cm] at (a) {$d$};

\fill [name intersections={of=t2 and s1, by={a}}]
       (a) circle (4pt);
\node[above left = 0.01cm] at (a) {$b$};

\fill [name intersections={of=t2 and s2, by={a}}]
       (a) circle (4pt);
\node[above left = 0.01cm] at (a) {$e$};

\fill [name intersections={of=t3 and s1, by={a}}]
       (a) circle (4pt);
\node[above left = 0.01cm] at (a) {$c$};

\fill [name intersections={of=t3 and s2, by={a}}]
       (a) circle (4pt);
\node[above left = 0.01cm] at (a) {$f$};

\end{tikzpicture}
 \begin{tikzpicture}[scale= 0.2, baseline={([yshift=-.5ex]current bounding box.center)},vertex/.style={anchor=base,circle,fill=black!25,minimum size=18pt,inner sep=2pt}]
  \draw [thick] (0,0) to [round right paren ] (0,15);
  \end{tikzpicture}
  \end{equation*}
  \begin{equation*}
 = D 
  \begin{tikzpicture}[scale= 0.2, baseline={([yshift=-.5ex]current bounding box.center)},vertex/.style={anchor=base,circle,fill=black!25,minimum size=18pt,inner sep=2pt}]
  \draw [thick] (0,0) to [round left paren ] (0,15);
  \end{tikzpicture}
  \ldots \pm 
  \begin{tikzpicture}[scale = 0.2, baseline={([yshift=-.5ex]current bounding box.center)},vertex/.style={anchor=base,circle,fill=black!25,minimum size=18pt,inner sep=2pt}]
  
%[thick,scale=0.6, every node/.style={transform shape}]
%help lines grid
%\draw[help lines] (0.1,0.1) grid (9.9,9.9);

% s axis
\draw[thick][dotted][->] (0,12) -- (0,-2);
\node [left] at (0,-2) {$s$};
% t axis
\draw[thick][dotted][->] (-2,10) -- (12,10);
\node [above] at (12,10) {$t$};
%label axes
%\node [above left] at (0,10) {0};
%\node[above] at (0,12) {$t=0$};
\draw[thick][dotted] (10,12) -- (10,-2);
%\node[above] at(10,12) {$t=1$};
\draw[thick][dotted] (-2,0) -- (12,0);
%\node[above] at(-2,0) {$s=1$};
%\node[above] at (-2,10) {$s=0$};

%s-lines
%s_1

%s_2
%label si's
\draw[name path=s12, thin, gray] (-2,8) -- (12,8);
\node[left] at (-2,8) {$s_1= s_2$};
%\draw[name path=s2, thin, gray] (-2,4) -- (12,4);
%$\node[left] at (-2,4) {$s_2$};

%t-lines

%t1
\draw[name path=t1, thin, gray] (2.5,-2) -- (2.5,12);
\node[above] at (2.5,12) {$t_1$};

%t2
\draw[name path=t2, thin, gray] (5,-2) -- (5,12);
\node[above] at (5,12) {$t_2$};

%t3
\draw[name path=t3, thin, gray] (7.5,-2) -- (7.5,12);
\node[above] at (7.5,12) {$t_3$};

\fill [name intersections={of=t1 and s12, by={a}}]
       (a) circle (4pt);
\node[above left= -0.05cm] at (a) {\small$ad$};

%\fill [name intersections={of=t1 and s2, by={a}}]
%       (a) circle (4pt);
%\node[above left = 0.01cm] at (a) {$d$};

\fill [name intersections={of=t2 and s12, by={a}}]
       (a) circle (4pt);
\node[above  left= -0.05cm] at (a) {\small$be$};

%\fill [name intersections={of=t2 and s2, by={a}}]
%       (a) circle (4pt);
%\node[above left = 0.01cm] at (a) {$e$};

\fill [name intersections={of=t3 and s12, by={a}}]
       (a) circle (4pt);
\node[above  left = -0.08cm] at (a) {\small$cf$};

\fill[white,fill opacity=.75]
  (current bounding box.south west) rectangle (current bounding box.north east)
  (4, 8.5) circle[radius=4cm+.5\pgflinewidth];

%\fill [name intersections={of=t3 and s2, by={a}}]
   %    (a) circle (4pt);
%\node[above left = 0.01cm] at (a) {$f$};
\end{tikzpicture}  \pm
\begin{tikzpicture}[scale= 0.2, baseline={([yshift=-.5ex]current bounding box.center)},vertex/.style={anchor=base,circle,fill=black!25,minimum size=18pt,inner sep=2pt}]
  
%[thick,scale=0.6, every node/.style={transform shape}]
%help lines grid
%\draw[help lines] (0.1,0.1) grid (9.9,9.9);

% s axis
\draw[thick][dotted][->] (0,12) -- (0,-2);
\node [left] at (0,-2) {$s$};
% t axis
\draw[thick][dotted][->] (-2,10) -- (12,10);
\node [above] at (12,10) {$t$};
%label axes
%\node [above left] at (0,10) {0};
%\node[above] at (0,12) {$t=0$};
\draw[thick][dotted] (10,12) -- (10,-2);
%\node[above] at(10,12) {$t=1$};
\draw[thick][dotted] (-2,0) -- (12,0);
%\node[above] at(-2,0) {$s=1$};
%\node[above] at (-2,10) {$s=0$};

%s-lines
%s_1

%s_2
%label si's
\draw[name path=s1, thin, gray] (-2,8) -- (12,8);
\node[left] at (-2,8) {$s_1$};
\draw[name path=s2, thin, gray] (-2,4) -- (12,4);
\node[left] at (-2,4) {$s_2$};

%t-lines

%t1
\draw[name path=t1, thin, gray] (2.5,-2) -- (2.5,12);
\node[above] at (2.5,12) {$t_1$};

%t2
\draw[name path=t23, thin, gray] (5,-2) -- (5,12);
\node[above, xshift=0.3cm ] at (5,12) {$t_2= t_3$};

%t3
%\draw[name path=t3, thin, gray] (7.5,-2) -- (7.5,12);
%\node[above] at (7.5,12) {$t_3$};

\fill [name intersections={of=t1 and s1, by={a}}]
       (a) circle (4pt);
\node[above left = -0.05cm] at (a) {\small$a$};

\fill [name intersections={of=t1 and s2, by={a}}]
       (a) circle (4pt);
\node[above left = -0.05cm] at (a) {\small$d$};

\fill [name intersections={of=t23 and s1, by={a}}]
       (a) circle (4pt);
\node[above left = -0.05cm] at (a) {\small$bc$};

\fill [name intersections={of=t23 and s2, by={a}}]
       (a) circle (4pt);
\node[above left = -0.05cm] at (a) {\small$ef$};

%\fill [name intersections={of=t3 and s1, by={a}}]
%       (a) circle (4pt);
%\node[above left = 0.01cm] at (a) {$c$};

%\fill [name intersections={of=t3 and s2, by={a}}]
%       (a) circle (4pt);
%\node[above left = 0.01cm] at (a) {$f$};
\fill[white,fill opacity=.75]
  (current bounding box.south west) rectangle (current bounding box.north east)
  (3.5, 6) circle[radius=4cm+.5\pgflinewidth];
\end{tikzpicture}
\ldots 
 \begin{tikzpicture}[scale= 0.2, baseline={([yshift=-.5ex]current bounding box.center)},vertex/.style={anchor=base,circle,fill=black!25,minimum size=18pt,inner sep=2pt}]
  \draw [thick] (0,0) to [round right paren ] (0,15);
  \end{tikzpicture}
\end{equation*}

\begin{equation*}\hspace*{-1cm} 
 =  \ldots \pm 
  \begin{tikzpicture}[scale = 0.2, baseline={([yshift=-.5ex]current bounding box.center)},vertex/.style={anchor=base,circle,fill=black!25,minimum size=18pt,inner sep=2pt}]
  
%[thick,scale=0.6, every node/.style={transform shape}]
%help lines grid
%\draw[help lines] (0.1,0.1) grid (9.9,9.9);

% s axis
%\draw[thick][dotted][->] (0,12) -- (0,-2);
%\node [left] at (0,-2) {$s$};
% t axis
%\draw[thick][dotted][->] (-2,10) -- (12,10);
%\node [above] at (12,10) {$t$};
%label axes
%\node [above left] at (0,10) {0};
%\node[above] at (0,12) {$t=0$};
%\draw[thick][dotted] (10,12) -- (10,-2);
%\node[above] at(10,12) {$t=1$};
%\draw[thick][dotted] (-2,0) -- (12,0);
%\node[above] at(-2,0) {$s=1$};
%\node[above] at (-2,10) {$s=0$};

%s-lines
%s_1

%s_2
%label si's
\draw[name path=s12, thin, gray] (-2,8) -- (12,8);
\node[left] at (-2,8) {$s_1= s_2$};
%\draw[name path=s2, thin, gray] (-2,4) -- (12,4);
%$\node[left] at (-2,4) {$s_2$};

%t-lines

%t1
\draw[name path=t1, thin, gray] (2.5,-2) -- (2.5,12);
\node[above] at (2.5,12) {$t_1$};

%t2
\draw[name path=t23, thin, gray] (5,-2) -- (5,12);
\node[above right] at (5,12) {$t_2= t_3$};

%t3
%\draw[name path=t3, thin, gray] (7.5,-2) -- (7.5,12);
%\node[above] at (7.5,12) {$t_3$};

\fill [name intersections={of=t1 and s12, by={a}}]
       (a) circle (4pt);
\node[above left= -0.05cm] at (a) {};

%\fill [name intersections={of=t1 and s2, by={a}}]
%       (a) circle (4pt);
%\node[above left = 0.01cm] at (a) {$d$};

\fill [name intersections={of=t23 and s12, by={a}}]
       (a) circle (4pt);
\node[above= 0.1cm, fill=white] at (a) {\small$becf$};

%\fill [name intersections={of=t2 and s2, by={a}}]
%       (a) circle (4pt);
%\node[above left = 0.01cm] at (a) {$e$};

%\fill [name intersections={of=t3 and s12, by={a}}]
%       (a) circle (4pt);
%\node[above  left = -0.08cm] at (a) {\small$cf$};

\fill[white,fill opacity=0.75]
  (current bounding box.south west) rectangle (current bounding box.north east)
  (5.5, 8.5) circle[radius=2.5cm+.5\pgflinewidth];

%\fill [name intersections={of=t3 and s2, by={a}}]
   %    (a) circle (4pt);
%\node[above left = 0.01cm] at (a) {$f$};
\end{tikzpicture}  \pm
\begin{tikzpicture}[scale= 0.2, baseline={([yshift=-.5ex]current bounding box.center)},vertex/.style={anchor=base,circle,fill=black!25,minimum size=18pt,inner sep=2pt}]
  
%[thick,scale=0.6, every node/.style={transform shape}]
%help lines grid
%\draw[help lines] (0.1,0.1) grid (9.9,9.9);

% s axis
%\draw[thick][dotted][->] (0,12) -- (0,-2);
%\node [left] at (0,-2) {$s$};
% t axis
%\draw[thick][dotted][->] (-2,10) -- (12,10);
%\node [above] at (12,10) {$t$};
%label axes
%\node [above left] at (0,10) {0};
%\node[above] at (0,12) {$t=0$};
%\draw[thick][dotted] (10,12) -- (10,-2);
%\node[above] at(10,12) {$t=1$};
%\draw[thick][dotted] (-2,0) -- (12,0);
%\node[above] at(-2,0) {$s=1$};
%\node[above] at (-2,10) {$s=0$};

%s-lines
%s_1

%s_2
%label si's
\draw[name path=s12, thin, gray] (-2,8) -- (12,8);
\node[left] at (-2,8) {$s_1= s_2$};
%\draw[name path=s2, thin, gray] (-2,4) -- (12,4);
%\node[left] at (-2,4) {$s_2$};

%t-lines

%t1
\draw[name path=t1, thin, gray] (2.5,-2) -- (2.5,12);
\node[above] at (2.5,12) {$t_1$};

%t2
\draw[name path=t23, thin, gray] (5,-2) -- (5,12);
\node[above, xshift=0.3cm ] at (5,12) {$t_2= t_3$};

%t3
%\draw[name path=t3, thin, gray] (7.5,-2) -- (7.5,12);
%\node[above] at (7.5,12) {$t_3$};

\fill [name intersections={of=t1 and s12, by={a}}]
       (a) circle (4pt);
\node[above left = -0.05cm] at (a) {};

%\fill [name intersections={of=t1 and s2, by={a}}]
%       (a) circle (4pt);
%\node[above left = -0.05cm] at (a) {\small$d$};

\fill [name intersections={of=t23 and s12, by={a}}]
       (a) circle (4pt);
\node[above = 0.1cm, fill=white] at (a) {\small$bcef$};

%\fill [name intersections={of=t23 and s2, by={a}}]
%       (a) circle (4pt);
%\node[above left = -0.05cm] at (a) {\small$ef$};

%\fill [name intersections={of=t3 and s1, by={a}}]
%       (a) circle (4pt);
%\node[above left = 0.01cm] at (a) {$c$};

%\fill [name intersections={of=t3 and s2, by={a}}]
%       (a) circle (4pt);
%\node[above left = 0.01cm] at (a) {$f$};
\fill[white,fill opacity=.75]
  (current bounding box.south west) rectangle (current bounding box.north east)
  (5, 8.5) circle[radius=2.5cm+.5\pgflinewidth];
\end{tikzpicture}
\ldots 
\end{equation*}

Instead of using just one row per $s_i$-coordinate, we adopt a model\footnote{Of course, one could use a "square model" or even $CH^{ZZ}(CH^{ZZ}(\OmegaMat(M)))$, but we adopt the rectangular model out of personal preference and comment on the relationships later.}  $CH^{ZZ}_{Rec}(\A)$  where we start with $k$-many zigzags at each $s_i$ and then allow our differential forms to be placed at the intersection-points of $t_i$'s and these zigzags.  The corresponding picture looks as follows:
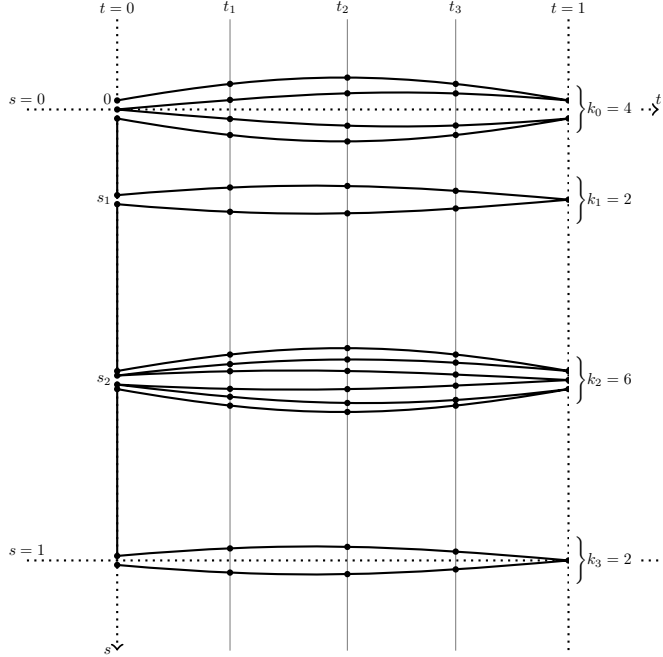
\begin{equation*}
\begin{tikzpicture}[thick,scale=0.6, every node/.style={transform shape}]
%help lines grid
%\draw[help lines] (0.1,0.1) grid (9.9,9.9);

% s axis
\draw[thick][dotted][->] (0,12) -- (0,-2);
\node [left] at (0,-2) {$s$};
% t axis
\draw[thick][dotted][->] (-2,10) -- (12,10);
\node [above] at (12,10) {$t$};
%label axes
\node [above left] at (0,10) {0};
\node[above] at (0,12) {$t=0$};
\draw[thick][dotted] (10,12) -- (10,-2);
\node[above] at(10,12) {$t=1$};
\draw[thick][dotted] (-2,0) -- (12,0);
\node[above] at(-2,0) {$s=1$};
\node[above] at (-2,10) {$s=0$};

% s=0 zigzags
\draw[name path=zigzag] (0,10.2) to [out=10, in=170] (10,10.2) to [out=175, in=5](0,10)
to [out=-5, in=185] (10,9.8) to [out=190, in=-10] (0,9.8)
% s=s_1 zigzazags
-- (0,8.1) to [out=5, in=175] (10,8) to [out=185, in=-5](0,7.9)

%s=s_2 zigzags
--(0, 4.2) to  [out=10, in=170] (10,4.2) to [out=175, in=7](0,4.1)
to [out=3, in=177] (10,4) to [out=183, in=-3] (0,3.9)
to [out=-7, in=187] (10,3.8) to [out=190, in=-10](0,3.8)

%s=1 zigzags
-- (0,0.1) to [out=5, in=175] (10,0) to [out=185, in=-5](0,-0.1);

\fill (0,10.2) circle (2pt); 
\fill (10,10.2) circle (2pt); 
\fill (0,10) circle (2pt); 
\fill (10,9.8) circle (2pt); 
\fill (0,9.8) circle (2pt); 
\fill (0,8.1) circle (2pt); 
\fill (10,8) circle (2pt); 
\fill (0,7.9) circle (2pt); 
\fill (0, 4.2) circle (2pt); 
\fill (10,4.2) circle (2pt); 
\fill (0,4.1) circle (2pt); 
\fill (10,4) circle (2pt); 
\fill (0,3.9) circle (2pt); 
\fill (10,3.8) circle (2pt); 
\fill (0,3.8) circle (2pt); 
\fill (0,0.1) circle (2pt); 
\fill (10,0) circle (2pt); 
\fill (0,-0.1) circle (2pt); 

%label si's
\node[left] at (0,8) {$s_1$};
\node[left] at (0,4) {$s_2$};

%label ki's
\node[right, fill=white] at (10,10) { $\Bigg\} k_0 = 4$};
\node[right, fill=white] at (10,8) { $\Bigg\} k_1 = 2$};
\node[right, fill=white] at (10,4) { $\Bigg\} k_2 = 6$};
\node[right, fill=white] at (10,0) { $\Bigg\} k_3 = 2$};

%t-lines

%t1
\draw[name path=t1, thin, gray] (2.5,-2) -- (2.5,12);
\node[above] at (2.5,12) {$t_1$};

%t2
\draw[name path=t2, thin, gray] (5.1,-2) -- (5.1,12);
\node[above] at (5,12) {$t_2$};

%t3
\draw[name path=t3, thin, gray] (7.5,-2) -- (7.5,12);
\node[above] at (7.5,12) {$t_3$};

\fill [name intersections={of=t1 and zigzag, by={a1, a2, a3, a4, a5, a6, a7, a8, a9, a10, a11, a12, a13, a14}}]
        \foreach \y in {1, 2, ..., 14}
        {(a\y) circle (2pt)};

 \fill [name intersections={of=t3 and zigzag, by={a1, a2, a3, a4, a5, a6, a7, a8, a9, a10, a11, a12, a13, a14}}]
        \foreach \y in {1, 2, ..., 14}
        {(a\y) circle (2pt)};
        
\fill [name intersections={of=t2 and zigzag, by={a1, a2, a3, a4, a5, a6, a7, a8, a9, a10, a11, a12, a13, a14}}]
        \foreach \y in {1, 2, ..., 14}
        {(a\y) circle (2pt)};

%\node[label=above:$r$] at (a) {};
\end{tikzpicture}
\end{equation*}
This model resolves the issues involved when $\OmegaMat(M)$ is non-abelian, since a collapse of rows only stacks zigzags together without changing the order in which the forms appear.  After defining $CH^{ZZ}_{Rec}(A)$, we proceed in Chapter 4 to define an iterated integral $It: CH^{ZZ}_{Rec}(\OmegaMat(M)) \to \OmegaMat(BM)$, for the space $BM$ of squares\footnote{We proceed with the notation $BM$ having ``bigons'' in mind for our last chapter.} $[0,1]^2 \to M$.

\begin{figure}
\begin{tikzpicture}[scale= 0.8, baseline={([yshift=-.5ex]current bounding box.center)},vertex/.style={anchor=base,circle,fill=black!25,minimum size=18pt,inner sep=2pt}]
%help lines grid
%\draw[help lines] (0.1,0.1) grid (9.9,9.9);

% s axis
\draw[thick][dotted][->] (0,12) -- (0,-2);
\node [left, lightgray] at (0,-2) {$s$};
% t axis
\draw[thick, dotted, lightgray][->] (-2,10) -- (12,10);
\node [above, lightgray] at (12,10) {$t$};
%label axes
\node [above left] at (0,10) {0};
\node[above, lightgray] at (0,12) {$t=0$};
\draw[thick][dotted] (10,12) -- (10,-2);
\node[above, lightgray] at(10,12) {$t=1$};
\draw[thick][dotted] (-2,0) -- (12,0);
\node[above,lightgray] at(-2,0) {$s=1$};
\node[above, lightgray] at (-2,10) {$s=0$};

% s=0 zigzags
\draw[name path=zigzag] (0,10) %to [out=10, in=170] (10,10.2) to [out=175, in=5](0,10)
%to [out=-5, in=185] (10,9.8) to [out=190, in=-10] (0,9.8)
% s=s_1 zigzazags
-- coordinate[pos=.26] (A21)  coordinate[pos=.4] (A22) coordinate[pos=.56] (A23) 
(0,8.1)  to [out=5, in=175]
coordinate[pos=.02] (A1)  coordinate[pos=.05] (A2)  coordinate[pos=.08] (A3)  coordinate[pos=.11] (A4) coordinate[pos=.14] (A5)
 coordinate[pos=.32] (A6)coordinate[pos=.38] (A7) coordinate[pos=.42] (A8) coordinate[pos=.54] (A9)(10,8) to [out=185, in=-5]    coordinate[pos=.34] (A10) coordinate[pos=.38] (A11) 
coordinate[pos=.57] (A12)  coordinate[pos=.61] (A13)  coordinate[pos=.65] (A14) coordinate[pos=.69] (A15)  
coordinate[pos=.85] (A16)  coordinate[pos=.89] (A17)  coordinate[pos=.93] (A18)
   (0,7.9)

%s=s_2 zigzags
--  coordinate[pos=.43] (A19)  coordinate[pos=.60] (A20) 
 (0,6.1) to [out=5, in=175] coordinate[pos=.13]  (A21) coordinate[pos=.33] (A22) coordinate[pos=.57] (A23) coordinate[pos=.60] (A24) coordinate[pos=.63] (A25) (10,6) to [out=185, in=-5] coordinate[pos=.44] (A26) coordinate[pos=.82] (A27) coordinate[pos=.86] (A28) coordinate[pos=.90] (A29)  coordinate[pos=.96] (A30) (0,5.9)
%--(0, 4.2) to  [out=10, in=170] (10,4.2) to [out=175, in=7](0,4.1)
%to [out=3, in=177] (10,4) to [out=183, in=-3] (0,3.9)
%to [out=-7, in=187] (10,3.8) to [out=190, in=-10](0,3.8)

%s=1 zigzags
--  coordinate[pos=.10] (A31) coordinate[pos=.18] (A32) coordinate[pos=.23] (A33) coordinate[pos=.30] (A34) coordinate[pos=.35] (A35)  (0,0);% to [out=5, in=175] (10,0) to [out=185, in=-5](0,-0.1);
\foreach \x in {1,2,...,9}
{  \draw ($(A\x)+(0,5pt)$) -- ($(A\x)-(0,5pt)$);

\node[above= 0.1cm, gray] at (A\x) {\small$A$};}

\foreach \x in {10,...,18}
{  \draw ($(A\x)+(0,5pt)$) -- ($(A\x)-(0,5pt)$);

\node[below= 0.1cm, gray] at (A\x) {\small$A$};}

\foreach \x in {19, 20}
{  \draw ($(A\x)+(5pt,0)$) -- ($(A\x)-(5pt,0)$);

\node[left= 0.1cm, gray] at (A\x) {\small$A$};}

\foreach \x in {21, 22, ..., 25}
{  \draw ($(A\x)+(0,5pt)$) -- ($(A\x)-(0,5pt)$);

\node[above= 0.1cm, gray] at (A\x) {\small$A$};}

\foreach \x in {26, 27, ..., 30}
{  \draw ($(A\x)+(0,5pt)$) -- ($(A\x)-(0,5pt)$);

\node[below= 0.1cm, gray] at (A\x) {\small$A$};}

\foreach \x in {31, 32, ..., 35}
{  \draw ($(A\x)+(5pt, 0)$) -- ($(A\x)-(5pt, 0)$);

\node[left= 0.1cm, gray] at (A\x) {\small$A$};}

%label si's
\node[left, lightgray] at (0,8) {$s_1$};
\node[left, lightgray] at (0,6) {$s_2$};

%label ki's
%\node[right, fill=white] at (10,10) { $\Bigg\} k_0 = 4$};
%\node[right, fill=white] at (10,8) { $\Bigg\} k_{s_1} = 2$};
%\node[right, fill=white] at (10,4) { $\Bigg\} k_{s_2} = 6$};
%\node[right, fill=white] at (10,0) { $\Bigg\} k_1 = 2$};

%t-lines

%t1
\draw[thin, lightgray, name path=t1] (2.5,-2) -- (2.5,12);
\node[above, lightgray] at (2.5,12) {$t_1$};

%t2
\draw[thin, lightgray, name path=t2] (5,-2) -- (5,12);
\node[above, lightgray] at (5,12) {$t_2$};

%t3
\draw[thin, gray, name path=t3] (7.5,-2) -- (7.5,12);
\node[above] at (7.5,12) {$t_3$};

\fill [name intersections={of=t1 and zigzag, by={a,b, c, d}}]
        (a) circle (2pt)
        (b) circle (2pt)
        (c) circle (2pt)
        (d) circle (2pt);
\node[above left] at (a) {$B$};
\node[below left = -0.02cm] at (b) {$1$};
\node[above left = -0.02cm] at (c) {$1$};
\node[below left = -0.02cm] at (d) {$1$};

\fill [name intersections={of=t2 and zigzag, by={a,b, c, d}}]
        (a) circle (2pt)
        (b) circle (2pt)
        (c) circle (2pt)
        (d) circle (2pt);
        
\node[above left = -0.02cm] at (a) {$1$};
\node[below left = -0.02cm] at (b) {$1$};
\node[above left] at (c) {$B$};
\node[below left = -0.02cm] at (d) {$1$};

\fill[white,fill opacity=.75]
  (current bounding box.south west) rectangle (current bounding box.north east)
  (2, 8) circle[radius=5cm+.5\pgflinewidth];

\end{tikzpicture}
\caption{An illustration of our curved iterated integral, $It^A$, applied to one term of $exp(B)$.}
\label{fig:It^A(exp(B))}
\end{figure}
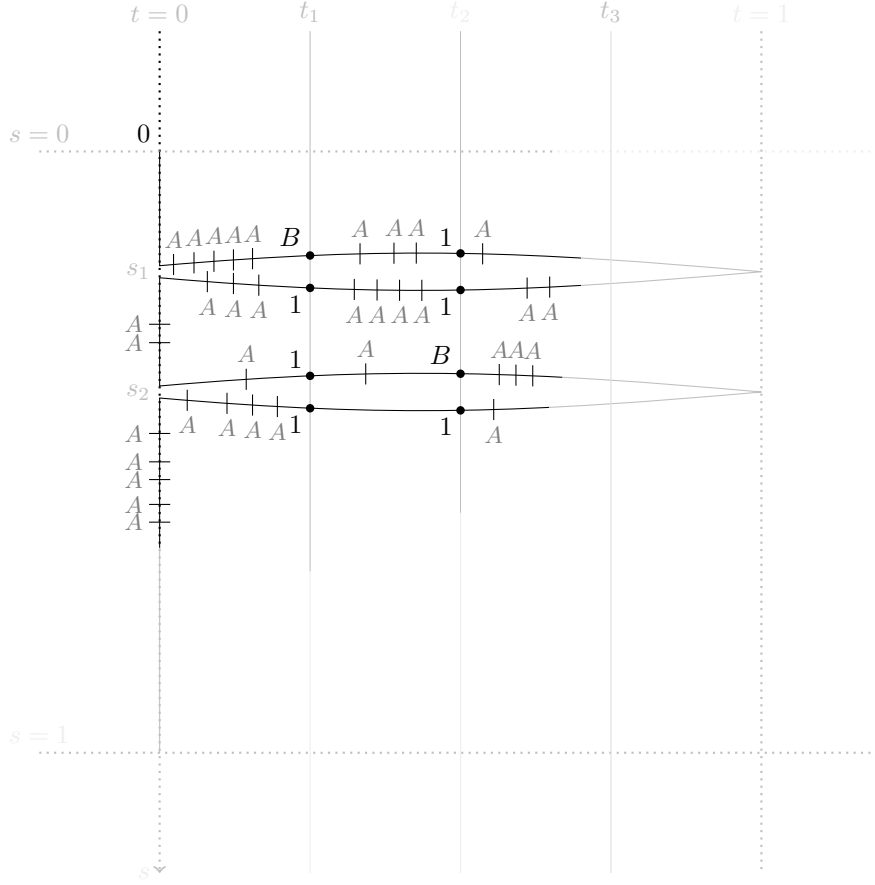

As mentioned at the beginning, our goal is to have an element in our Hochschild complex which maps to holonomy under our iterated integral maps.  In the commutative setting, such an element representing 2-holonomy in $CH^{I \times I}(\Omega(M))$ is given by an infinite sum (the exponential of some element) of monomials of the form, see \cite{TWZ},
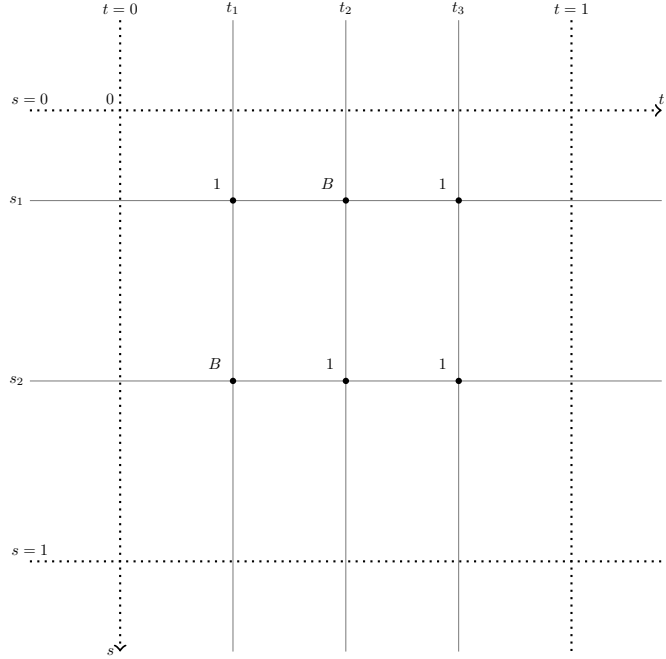
\begin{equation*}
\begin{tikzpicture}[thick,scale=0.6, every node/.style={transform shape}]
%help lines grid
%\draw[help lines] (0.1,0.1) grid (9.9,9.9);

% s axis
\draw[thick][dotted][->] (0,12) -- (0,-2);
\node [left] at (0,-2) {$s$};
% t axis
\draw[thick][dotted][->] (-2,10) -- (12,10);
\node [above] at (12,10) {$t$};
%label axes
\node [above left] at (0,10) {0};
\node[above] at (0,12) {$t=0$};
\draw[thick][dotted] (10,12) -- (10,-2);
\node[above] at(10,12) {$t=1$};
\draw[thick][dotted] (-2,0) -- (12,0);
\node[above] at(-2,0) {$s=1$};
\node[above] at (-2,10) {$s=0$};

%s-lines
%s_1

%s_2
%label si's
\draw[name path=s1, thin, gray] (-2,8) -- (12,8);
\node[left] at (-2,8) {$s_1$};
\draw[name path=s2, thin, gray] (-2,4) -- (12,4);
\node[left] at (-2,4) {$s_2$};

%t-lines

%t1
\draw[name path=t1, thin, gray] (2.5,-2) -- (2.5,12);
\node[above] at (2.5,12) {$t_1$};

%t2
\draw[name path=t2, thin, gray] (5,-2) -- (5,12);
\node[above] at (5,12) {$t_2$};

%t3
\draw[name path=t3, thin, gray] (7.5,-2) -- (7.5,12);
\node[above] at (7.5,12) {$t_3$};

\fill [name intersections={of=t1 and s1, by={a}}]
       (a) circle (2pt);
\node[label=above left:$1$] at (a) {};

\fill [name intersections={of=t1 and s2, by={a}}]
       (a) circle (2pt);
\node[label=above left:$B$] at (a) {};

\fill [name intersections={of=t2 and s1, by={a}}]
       (a) circle (2pt);
\node[label=above left:$B$] at (a) {};

\fill [name intersections={of=t2 and s2, by={a}}]
       (a) circle (2pt);
\node[label=above left:$1$] at (a) {};

\fill [name intersections={of=t3 and s1, by={a}}]
       (a) circle (2pt);
\node[label=above left:$1$] at (a) {};

\fill [name intersections={of=t3 and s2, by={a}}]
       (a) circle (2pt);
\node[label=above left:$1$] at (a) {};
\end{tikzpicture}
\end{equation*}
where $B \in \Omega^2(U)$ is a connection 2-form on an open set $U \subset M$.  Through the works of \cite{BaSc}, \cite{PM1}, \cite{PM2}, \cite{ScWa} et al we observe that the 2-holonomy we are interested in (for the non-commutative setting) can be expressed loosely as $$\iint \alpha_*(B)$$ where $\alpha_*$ is the action data coming from a crossed module $(\mathfrak{h} \xrightarrow{t} \mathfrak{g} \xrightarrow{\alpha_*} Der(\mathfrak{h}))$.  The formula is similar to our abelian situation but is more accurately written as $$P\exp \left( \int_0^t \mathcal{P}_{(t',s)} (B)(s,t') \mathcal{P}_{(t',s)}^{-1} dt' \right),$$ involving a path ordered exponential and parallel transport, $P_{(t',s)}$, from $(0,0)$ to $(t',s)$, via a given connection 1-form, $A \in \OmegaMat^1(U)$.  All of this means we are interested in integrating over a sum of monomials as illustrated in Figure \ref{fig:It^A(exp(B))} on page \pageref{fig:It^A(exp(B))}.
%%%%%%%%%%%%%%%%%%%%%%%%%%%%%%%%%%%%%%%%%%%%%%%%%%%%
In all of the above structures, we obtain our differentials on the various Hochschild complexes by requiring the iterated integral to be a chain map.  For example, when working in the 1-dimensional non-curved zigzag Hochschild Complex, $CH^{ZZ}(\A)$, we require our differential to be defined $D= (d+b)$ where $b$ and $d$ come from the terms in Stokes' formula,  $$d\int\limits_{F} \omega = (-1)^{dim{F}} \int\limits_{F} d\omega + (-1)^{dim(F) -1} \int\limits_{\partial F} \omega,$$ which differentiate the integrand and take the boundary of the fiber, respectively.  For the 2-dimensional non-curved zigzag Hochschild Complex, $CH^{ZZ}_{Rec}(\A)$, the boundary term in Stokes' formula also gives a vertical boundary, which we will call $\star$, so that our differential becomes $D = (d + b+ \star)$.  If we use the above figure as a guide to what happens in the curved 2-dimensional case, where we apply parallel transport between the forms originally placed along our zigzags, the boundary term in Stokes' also yields terms of the form $A \wedge A$, whereas the differential of the integrand adds terms $dA$.  Putting these terms together we see that $d_{DR}$ of our Iterated Integral will add terms $R:= dA + A \wedge A$ in between our differential forms on the zigzags.  For this reason, the curved zigzag (both 1-d and 2-d) Hochschild complex has a $c$ component in the differential, which precisely shuffles in these $R$'s.  In addition, one can see that the boundary term in Stokes' results in terms $A\wedge \omega$ and $\omega \wedge A$, where $\omega$ is some differential form originally placed on the zigzag.  These terms require us to replace the $d$ component of our differential with a component $\nabla = d + [A, -]$.  Since $\nabla^2 = [R, -] \ne 0$ in this case, we call these complexes {\it curved} and we use the differentials $D= (\nabla + b + c)$ for the 1-dimensional curved zigzag Hochschild complex, $CH^{ZZ}(\A)$, and $D= (\nabla + b + c + \star)$ for the 2-dimensional curved zigzag Hochschild complex, $CH^{ZZ}_{Rec}(\A)$.  Chapter 5 deals with our zigzag Hochschild complexes in this curved setting, complete with iterated integral maps.  

Finally, in Chapter 6 we verify that we have elements in $CH^{ZZ}(\OmegaMat(M))$ which map (in their limit) to holonomy, and elements in $CH^{ZZ}_{Rec}(\OmegaMat(M))$ which map to 2-holonomy under the curved iterated integral.  
\subsection{Acknowledgements}
I would like to thank Thomas Tradler for his insights, attention to detail, and constant encouragement.   I would also like to thank Mahmoud Zeinalian for his helpful conversations.
 \section{The one-dimensional zigzag Hochschild complex}
In this chapter, we define the complex $CH^{ZZ}(\A)$ along with a product making it into a DGA.
\subsection{For a (non-abelian) DGA}
Let $(\A, \cdot, d)$ be a (possibly non-commutative) associative, unital DGA over a commutative ring $S$.  We have the well-known interval Hochschild complex\footnote{Also known as the two-sided bar construction where the modules are chosen to be $\A$ in this case}, 
$CH^I(\A):= \bigoplus_{n \ge 0} \A \otimes (\A[1])^{\otimes n}\otimes \A$,
with its usual differential, $D$, given by
\begin{align*}
D(\omega_0 \otimes \omega_1\otimes \ldots \otimes \omega_n \otimes \omega_{n+1}):= &\sum\limits_{i=0}^{n+1} (-1)^{n+\beta_i}\omega_0 \otimes \ldots \otimes d(\omega_i) \otimes \ldots \otimes \omega_{n+1}\\
&\sum\limits_{i=0}^{n} (-1)^{i}\omega_0 \otimes \ldots \otimes( \omega_i \cdot \omega_{i+1})\otimes \ldots \otimes \omega_{n+1}
\end{align*}
However under the usual shuffle product, $D$ is not a derivation; that would require $\A$ to be a commutative algebra.  In this section, we define a new Hochschild complex, $(CH^{ZZ}(\A), D)$ with an associative shuffle product, $\odot$, for which our differential $D$ is a derivation.  \\
Before providing a formal definition, we will give the underlying ideas for its definition.  For our underlying vector space, which will consist of tensor products of elements in $\A$, monomials can easily be represented by diagrams 
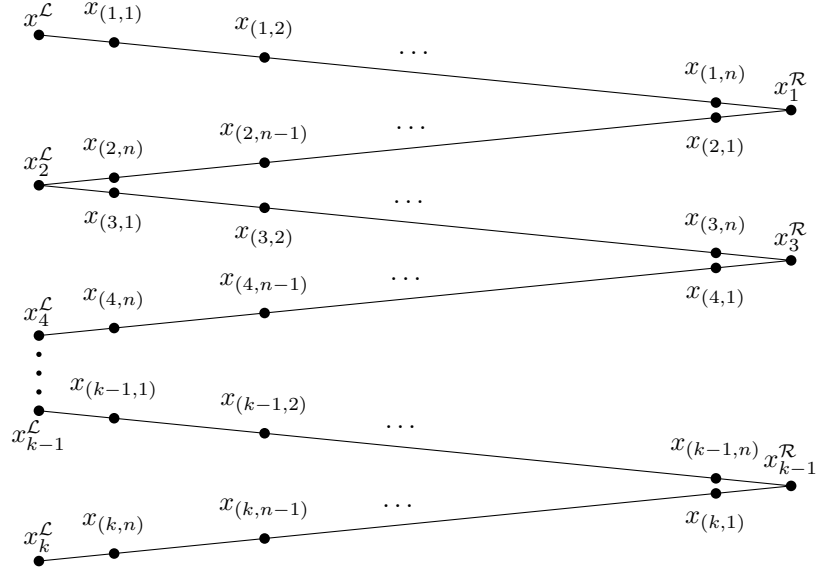
\begin{equation*}\label{general zigzag}
\begin{tikzpicture}
\draw[name path=zigzagtop] (0,10) to (10,9) to (0,8) to (10,7) to (0,6);
 \fill (0,5.75) circle (1pt);
 \fill (0,5.5) circle (1pt);
 \fill (0,5.25) circle (1pt);
\draw[name path=zigzagbot] (0,5) to (10,4) to (0,3); 
\fill (0,10) circle (2pt) node[above] {$x^{\mathcal{L}}$};
\fill (0,8) circle (2pt) node[above] {$x_2^{\mathcal{L}}$};
\fill (0,6) circle (2pt) node[above] {$x_4^{\mathcal{L}}$};
\fill (0,5) circle (2pt) node[below] {$x_{k-1}^{\mathcal{L}}$};
\fill (0,3) circle (2pt) node[above] {$x_k^{\mathcal{L}}$};
\fill (10,9) circle (2pt) node[above] {$x_1^{\mathcal{R}}$};
\fill (10,7) circle (2pt) node[above] {$x_3^{\mathcal{R}}$};
\fill (10,4) circle (2pt) node[above] {$x_{k-1}^{\mathcal{R}}$};
\draw[opacity=0,name path=t1] (1, 10) to (1,2.9);
\fill [name intersections={of=t1 and zigzagtop, by={a1,b1, c1, d1}}]
        (a1) circle (2pt)
        (b1) circle (2pt)
        (c1) circle (2pt)
        (d1) circle (2pt);
\node[label=above:{$x_{(1,1)}$}] at (a1) {};
\node[label=above:$x_{(2,n)}$] at (b1) {};
\node[label=below:$x_{(3,1)}$] at (c1) {};
\node[label=above:$x_{(4,n)}$] at (d1) {};

\fill [name intersections={of=t1 and zigzagbot, by={e1,f1}}]
        (e1) circle (2pt)
        (f1) circle (2pt);
\node[label=above:{$x_{(k-1,1)}$}] at (e1) {};
\node[label=above:$x_{(k,n)}$] at (f1) {};

\draw[opacity=0,name path=t2] (3, 10) to (3,2.9);
\fill [name intersections={of=t2 and zigzagtop, by={a2,b2, c2, d2}}]
        (a2) circle (2pt)
        (b2) circle (2pt)
        (c2) circle (2pt)
        (d2) circle (2pt);
\node[label=above:{$x_{(1,2)}$}] at (a2) {};
\node[label=above:$x_{(2,n-1)}$] at (b2) {};
\node[label=below:$x_{(3,2)}$] at (c2) {};
\node[label=above:$x_{(4,n-1)}$] at (d2) {};

\fill [name intersections={of=t2 and zigzagbot, by={e2,f2}}]
        (e2) circle (2pt)
        (f2) circle (2pt);
\node[label=above:{$x_{(k-1,2)}$}] at (e2) {};
\node[label=above:$x_{(k,n-1)}$] at (f2) {};

%\draw[opacity=0, name path=t2] (2, 10) to (2,0);

%\draw[opacity=0, name path=t3] (8, 10) to (8,0);

\draw[opacity=0, name path=t4] (9, 10) to (9,2.9);
\fill [name intersections={of=t4 and zigzagtop, by={a4,b4, c4, d4}}]
        (a4) circle (2pt)
        (b4) circle (2pt)
        (c4) circle (2pt)
        (d4) circle (2pt);
\node[label=above:{$x_{(1,n)}$}] at (a4) {};
\node[label=below:$x_{(2,1)}$] at (b4) {};
\node[label=above:$x_{(3,n)}$] at (c4) {};
\node[label=below:$x_{(4,1)}$] at (d4) {};

\fill [name intersections={of=t4 and zigzagbot, by={e4,f4}}]
        (e4) circle (2pt)
        (f4) circle (2pt);
\node[label=above:{$x_{(k-1,n)}$}] at (e4) {};
\node[label=below:$x_{(k,1)}$] at (f4) {};

%zigzapellipses
\draw[opacity=0, name path=dotL] (5, 10) to (4.75,2.9);
\path [name intersections={of=dotL and zigzagtop, by={dotL1, dotL2, dotL3, dotL4}}];

\foreach \n in {1,...,4} {
\node[label=above: $\ldots$] at (dotL\n) {};}

\path [name intersections={of=dotL and zigzagbot, by={dotL5, dotL6}}];
\node[label=above: $\ldots$] at (dotL5) {};
\node[label=above: $\ldots$] at (dotL6) {};
\end{tikzpicture}
\end{equation*}
where we've replaced the interval with $k$-many zigzags going back and forth over the same interval.  An odd-numbered trip over the interval (left-to-right) is referred to as a ``zig'' and an even-numbered trip over the interval (right-to-left) is called a ``zag''.  The elements of $\A$ are placed at the points where the zigzags cross the $n$-many columns\footnote{These columns represent the ``time-slots''  $0\le t_1 \le \ldots \le t_n \le 1$ which we will eventually integrate over when considering differential forms and an iterated integral.}.  Consider the following monomial with $n=2$ and $k=4$ for further clarification:
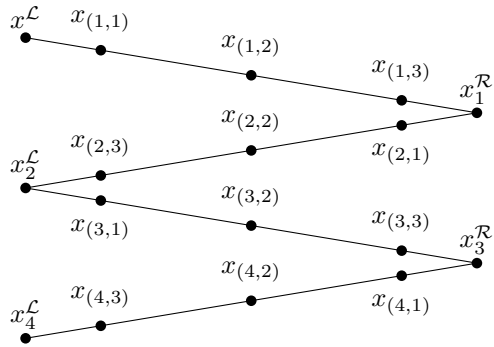
\abovedisplayskip=5pt
 \belowdisplayskip=5pt
\begin{equation*} \label{simple monomial pic}
\begin{tikzpicture}
\draw[name path=zigzagtop] (0,10) to (6,9) to (0,8) to (6,7) to (0,6); 
\fill (0,10) circle (2pt) node[above] {$x^{\mathcal{L}}$};
\fill (0,8) circle (2pt) node[above] {$x_2^{\mathcal{L}}$};
\fill (0,6) circle (2pt) node[above] {$x_4^{\mathcal{L}}$};
\fill (6,9) circle (2pt) node[above] {$x_1^{\mathcal{R}}$};
\fill (6,7) circle (2pt) node[above] {$x_3^{\mathcal{R}}$};
\draw[opacity=0,name path=t1] (1, 10) to (1,6);
\fill [name intersections={of=t1 and zigzagtop, by={a1,b1, c1, d1}}]
        (a1) circle (2pt)
        (b1) circle (2pt)
        (c1) circle (2pt)
        (d1) circle (2pt);
\node[label=above:{$x_{(1,1)}$}] at (a1) {};
\node[label=above:$x_{(2,3)}$] at (b1) {};
\node[label=below:$x_{(3,1)}$] at (c1) {};
\node[label=above:$x_{(4,3)}$] at (d1) {};

\draw[opacity=0,name path=t2] (3, 10) to (3,6);
\fill [name intersections={of=t2 and zigzagtop, by={a2,b2, c2, d2}}]
        (a2) circle (2pt)
        (b2) circle (2pt)
        (c2) circle (2pt)
        (d2) circle (2pt);
\node[label=above:{$x_{(1,2)}$}] at (a2) {};
\node[label=above:$x_{(2,2)}$] at (b2) {};
\node[label=above:$x_{(3,2)}$] at (c2) {};
\node[label=above:$x_{(4,2)}$] at (d2) {};

\draw[opacity=0, name path=t3] (5, 10) to (5,6);
\fill [name intersections={of=t3 and zigzagtop, by={a3,b3, c3, d3}}]
        (a3) circle (2pt)
        (b3) circle (2pt)
        (c3) circle (2pt)
        (d3) circle (2pt);
\node[label=above:{$x_{(1,3)}$}] at (a3) {};
\node[label=below:$x_{(2,1)}$] at (b3) {};
\node[label=above:$x_{(3,3)}$] at (c3) {};
\node[label=below:$x_{(4,1)}$] at (d3) {};
\end{tikzpicture}
\end{equation*}
The differential $D$ will consist of two components: $d$, which applies the differential coming from $\A$ to each element in the tensor product \label{dandbpics}
\begin{equation*}\resizebox{12cm}{!}{
 d{\begin{tikzpicture}[ baseline={([yshift=-.5ex]current bounding box.center)},vertex/.style={anchor=base,circle,fill=black!25,minimum size=18pt,inner sep=2pt}]
  \draw [thick] (0,0) to [round left paren ] (0,5);
  \end{tikzpicture}}
{ \begin{tikzpicture}[ baseline={([yshift=-.5ex]current bounding box.center)},vertex/.style={anchor=base,circle,fill=black!25,minimum size=18pt,inner sep=2pt}]
  \draw[name path=zigzag1] (0,10) to (3,9) to (0,8) to (3,7) to (0,6);
   %increase bounding box
   %\draw[opacity=0] (-0.2, -0.2) rectangle (18, 12.2);
   
   \draw[opacity=0,name path=t1] (0, 10) to (0,5);
   \fill [name intersections={of=t1 and zigzag1, by={a,blank, i,q}}]
        (a) circle (2pt)
        (i) circle (2pt)
        (q) circle (2pt);
\node[label=left:{$a$}] at (a) {};
\node[label=left:{$i$}] at (i) {};
\node[label=left:{$q$}] at (q) {};

\draw[opacity=0,name path=t2] (0.5, 10) to (0.5,6);
   \fill [name intersections={of=t2 and zigzag1, by={b,h, j, p}}]
        \foreach \x in {b, h, p}
        {(\x) circle (2pt)}
        (j) circle (2pt); 
    \foreach \x in {b, h, p}
        {\node[above=2pt] at (\x) {$\x$};}
        \node[below] at (j) {$j$};

\draw[opacity=0,name path=t3] (1.5, 10) to (1.5,6);
   \fill [name intersections={of=t3 and zigzag1, by={c,g, k, o}}]
        \foreach \x in {c,g, k, o}
        {(\x) circle (2pt)};
        \foreach \x in {c,g, k, o}
        {\node[label=above:{$\x$}] at (\x) {};}
        
\draw[opacity=0,name path=t4] (2, 10) to (2,6);
   \fill [name intersections={of=t4 and zigzag1, by={d,f, l, n}}]
        \foreach \x in {d,f, l, n}
        {(\x) circle (2pt)};
        \foreach \x in {d,f, l, n}
        {\node[label=above:{$\x$}] at (\x) {};}
        
\draw[opacity=0,name path=t5] (3, 10) to (3,6);
   \fill [name intersections={of=t5 and zigzag1, by={e,blank, m}}]
        \foreach \x in {e,m}
        {(\x) circle (2pt)};
        \foreach \x in {e,m}
        {\node[label=above:{$\x$}] at (\x) {};}
\end{tikzpicture}}
{\begin{tikzpicture}[ baseline={([yshift=-.5ex]current bounding box.center)},vertex/.style={anchor=base,circle,fill=black!25,minimum size=18pt,inner sep=2pt}]
  \draw [thick] (0,0) to [round right paren ] (0,5);
  \end{tikzpicture}} = 
   { \begin{tikzpicture}[ baseline={([yshift=-.5ex]current bounding box.center)},vertex/.style={anchor=base,circle,fill=black!25,minimum size=18pt,inner sep=2pt}]
  \draw[name path=zigzag1] (0,10) to (3,9) to (0,8) to (3,7) to (0,6);
   %increase bounding box
   %\draw[opacity=0] (-0.2, -0.2) rectangle (18, 12.2);
   
   \draw[opacity=0,name path=t1] (0, 10) to (0,5);
   \fill [name intersections={of=t1 and zigzag1, by={a,blank, i,q}}]
        (a) circle (2pt)
        (i) circle (2pt)
        (q) circle (2pt);
\node[label=left:{$d(a)$}] at (a) {};
\node[label=left:{$i$}] at (i) {};
\node[label=left:{$q$}] at (q) {};

\draw[opacity=0,name path=t2] (0.5, 10) to (0.5,6);
   \fill [name intersections={of=t2 and zigzag1, by={b,h, j, p}}]
        \foreach \x in {b, h, p}
        {(\x) circle (2pt)}
        (j) circle (2pt); 
    \foreach \x in {b, h, p}
        {\node[above=2pt] at (\x) {$\x$};}
        \node[below] at (j) {$j$};

\draw[opacity=0,name path=t3] (1.5, 10) to (1.5,6);
   \fill [name intersections={of=t3 and zigzag1, by={c,g, k, o}}]
        \foreach \x in {c,g, k, o}
        {(\x) circle (2pt)};
        \foreach \x in {c,g, k, o}
        {\node[label=above:{$\x$}] at (\x) {};}
        
\draw[opacity=0,name path=t4] (2, 10) to (2,6);
   \fill [name intersections={of=t4 and zigzag1, by={d,f, l, n}}]
        \foreach \x in {d,f, l, n}
        {(\x) circle (2pt)};
        \foreach \x in {d,f, l, n}
        {\node[label=above:{$\x$}] at (\x) {};}
        
\draw[opacity=0,name path=t5] (3, 10) to (3,6);
   \fill [name intersections={of=t5 and zigzag1, by={e,blank, m}}]
        \foreach \x in {e,m}
        {(\x) circle (2pt)};
        \foreach \x in {e,m}
        {\node[label=above:{$\x$}] at (\x) {};}
\end{tikzpicture}}
$\pm$  { \begin{tikzpicture}[ baseline={([yshift=-.5ex]current bounding box.center)},vertex/.style={anchor=base,circle,fill=black!25,minimum size=18pt,inner sep=2pt}]
  \draw[name path=zigzag1] (0,10) to (3,9) to (0,8) to (3,7) to (0,6);
   %increase bounding box
   %\draw[opacity=0] (-0.2, -0.2) rectangle (18, 12.2);
   
   \draw[opacity=0,name path=t1] (0, 10) to (0,5);
   \fill [name intersections={of=t1 and zigzag1, by={a,blank, i,q}}]
        (a) circle (2pt)
        (i) circle (2pt)
        (q) circle (2pt);
\node[label=left:{$a$}] at (a) {};
\node[label=left:{$i$}] at (i) {};
\node[label=left:{$q$}] at (q) {};

\draw[opacity=0,name path=t2] (0.5, 10) to (0.5,6);
   \fill [name intersections={of=t2 and zigzag1, by={b,h, j, p}}]
        \foreach \x in {b, h, p}
        {(\x) circle (2pt)}
        (j) circle (2pt); 
    \foreach \x in {h, p}
        {\node[above=2pt] at (\x) {$\x$};}
        \node[above=2pt] at (b) {$d(b)$};
        \node[below] at (j) {$j$};

\draw[opacity=0,name path=t3] (1.5, 10) to (1.5,6);
   \fill [name intersections={of=t3 and zigzag1, by={c,g, k, o}}]
        \foreach \x in {c,g, k, o}
        {(\x) circle (2pt)};
        \foreach \x in {c,g, k, o}
        {\node[label=above:{$\x$}] at (\x) {};}
        
\draw[opacity=0,name path=t4] (2, 10) to (2,6);
   \fill [name intersections={of=t4 and zigzag1, by={d,f, l, n}}]
        \foreach \x in {d,f, l, n}
        {(\x) circle (2pt)};
        \foreach \x in {d,f, l, n}
        {\node[label=above:{$\x$}] at (\x) {};}
        
\draw[opacity=0,name path=t5] (3, 10) to (3,6);
   \fill [name intersections={of=t5 and zigzag1, by={e,blank, m}}]
        \foreach \x in {e,m}
        {(\x) circle (2pt)};
        \foreach \x in {e,m}
        {\node[label=above:{$\x$}] at (\x) {};}
\end{tikzpicture}} $\pm \ldots \pm$ 
 { \begin{tikzpicture}[ baseline={([yshift=-.5ex]current bounding box.center)},vertex/.style={anchor=base,circle,fill=black!25,minimum size=18pt,inner sep=2pt}]
  \draw[name path=zigzag1] (0,10) to (3,9) to (0,8) to (3,7) to (0,6);
   %increase bounding box
   %\draw[opacity=0] (-0.2, -0.2) rectangle (18, 12.2);
   
   \draw[opacity=0,name path=t1] (0, 10) to (0,5);
   \fill [name intersections={of=t1 and zigzag1, by={a,blank, i,q}}]
        (a) circle (2pt)
        (i) circle (2pt)
        (q) circle (2pt);
\node[label=left:{$a$}] at (a) {};
\node[label=left:{$i$}] at (i) {};
\node[label=left:{$d(q)$}] at (q) {};

\draw[opacity=0,name path=t2] (0.5, 10) to (0.5,6);
   \fill [name intersections={of=t2 and zigzag1, by={b,h, j, p}}]
        \foreach \x in {b, h, p}
        {(\x) circle (2pt)}
        (j) circle (2pt); 
    \foreach \x in {b, h, p}
        {\node[above=2pt] at (\x) {$\x$};}
        \node[below] at (j) {$j$};

\draw[opacity=0,name path=t3] (1.5, 10) to (1.5,6);
   \fill [name intersections={of=t3 and zigzag1, by={c,g, k, o}}]
        \foreach \x in {c,g, k, o}
        {(\x) circle (2pt)};
        \foreach \x in {c,g, k, o}
        {\node[label=above:{$\x$}] at (\x) {};}
        
\draw[opacity=0,name path=t4] (2, 10) to (2,6);
   \fill [name intersections={of=t4 and zigzag1, by={d,f, l, n}}]
        \foreach \x in {d,f, l, n}
        {(\x) circle (2pt)};
        \foreach \x in {d,f, l, n}
        {\node[label=above:{$\x$}] at (\x) {};}
        
\draw[opacity=0,name path=t5] (3, 10) to (3,6);
   \fill [name intersections={of=t5 and zigzag1, by={e,blank, m}}]
        \foreach \x in {e,m}
        {(\x) circle (2pt)};
        \foreach \x in {e,m}
        {\node[label=above:{$\x$}] at (\x) {};}
\end{tikzpicture}}}
  \end{equation*}
  and $b$ which identifies two columns. 
   \begin{equation*}\resizebox{12cm}{!}{
 b{\begin{tikzpicture}[ baseline={([yshift=-.5ex]current bounding box.center)},vertex/.style={anchor=base,circle,fill=black!25,minimum size=18pt,inner sep=2pt}]
  \draw [thick] (0,0) to [round left paren ] (0,5);
  \end{tikzpicture}}
{ \begin{tikzpicture}[ baseline={([yshift=-.5ex]current bounding box.center)},vertex/.style={anchor=base,circle,fill=black!25,minimum size=18pt,inner sep=2pt}]
  \draw[name path=zigzag1] (0,10) to (3,9) to (0,8) to (3,7) to (0,6);
   %increase bounding box
   %\draw[opacity=0] (-0.2, -0.2) rectangle (18, 12.2);
   
   \draw[opacity=0,name path=t1] (0, 10) to (0,5);
   \fill [name intersections={of=t1 and zigzag1, by={a,blank, i,q}}]
        (a) circle (2pt)
        (i) circle (2pt)
        (q) circle (2pt);
\node[label=left:{$a$}] at (a) {};
\node[label=left:{$i$}] at (i) {};
\node[label=left:{$q$}] at (q) {};

\draw[opacity=0,name path=t2] (0.5, 10) to (0.5,6);
   \fill [name intersections={of=t2 and zigzag1, by={b,h, j, p}}]
        \foreach \x in {b, h, p}
        {(\x) circle (2pt)}
        (j) circle (2pt); 
    \foreach \x in {b, h, p}
        {\node[above=2pt] at (\x) {$\x$};}
        \node[below] at (j) {$j$};

\draw[opacity=0,name path=t3] (1.5, 10) to (1.5,6);
   \fill [name intersections={of=t3 and zigzag1, by={c,g, k, o}}]
        \foreach \x in {c,g, k, o}
        {(\x) circle (2pt)};
        \foreach \x in {c,g, k, o}
        {\node[label=above:{$\x$}] at (\x) {};}
        
\draw[opacity=0,name path=t4] (2, 10) to (2,6);
   \fill [name intersections={of=t4 and zigzag1, by={d,f, l, n}}]
        \foreach \x in {d,f, l, n}
        {(\x) circle (2pt)};
        \foreach \x in {d,f, l, n}
        {\node[label=above:{$\x$}] at (\x) {};}
        
\draw[opacity=0,name path=t5] (3, 10) to (3,6);
   \fill [name intersections={of=t5 and zigzag1, by={e,blank, m}}]
        \foreach \x in {e,m}
        {(\x) circle (2pt)};
        \foreach \x in {e,m}
        {\node[label=above:{$\x$}] at (\x) {};}
\end{tikzpicture}}
{\begin{tikzpicture}[ baseline={([yshift=-.5ex]current bounding box.center)},vertex/.style={anchor=base,circle,fill=black!25,minimum size=18pt,inner sep=2pt}]
  \draw [thick] (0,0) to [round right paren ] (0,5);
  \end{tikzpicture}} = 
   { \begin{tikzpicture}[ baseline={([yshift=-.5ex]current bounding box.center)},vertex/.style={anchor=base,circle,fill=black!25,minimum size=18pt,inner sep=2pt}]
  \draw[name path=zigzag1] (0,10) to (3,9) to (0,8) to (3,7) to (0,6);
   %increase bounding box
   %\draw[opacity=0] (-0.2, -0.2) rectangle (18, 12.2);
   
   \draw[opacity=0,name path=t1] (0, 10) to (0,5);
   \fill [name intersections={of=t1 and zigzag1, by={a,blank, i,q}}]
        (a) circle (2pt)
        (i) circle (2pt)
        (q) circle (2pt);
\node[label=left:{$ab$}] at (a) {};
\node[label=left:{$hij$}] at (i) {};
\node[label=left:{$pq$}] at (q) {};

%\draw[opacity=0,name path=t2] (0.5, 10) to (0.5,6);
%   \fill [name intersections={of=t2 and zigzag1, by={b,h, j, p}}]
 %       \foreach \x in {b, h, p}
  %      {(\x) circle (2pt)}
  %      (j) circle (2pt); 
  %  \foreach \x in {b, h, p}
  %      {\node[above=2pt] at (\x) {$\x$};}
  %      \node[below] at (j) {$j$};

\draw[opacity=0,name path=t3] (1.5, 10) to (1.5,6);
   \fill [name intersections={of=t3 and zigzag1, by={c,g, k, o}}]
        \foreach \x in {c,g, k, o}
        {(\x) circle (2pt)};
        \foreach \x in {c,g, k, o}
        {\node[label=above:{$\x$}] at (\x) {};}
        
\draw[opacity=0,name path=t4] (2, 10) to (2,6);
   \fill [name intersections={of=t4 and zigzag1, by={d,f, l, n}}]
        \foreach \x in {d,f, l, n}
        {(\x) circle (2pt)};
        \foreach \x in {d,f, l, n}
        {\node[label=above:{$\x$}] at (\x) {};}
        
\draw[opacity=0,name path=t5] (3, 10) to (3,6);
   \fill [name intersections={of=t5 and zigzag1, by={e,blank, m}}]
        \foreach \x in {e,m}
        {(\x) circle (2pt)};
        \foreach \x in {e,m}
        {\node[label=above:{$\x$}] at (\x) {};}
\end{tikzpicture}}
 $\pm$  { \begin{tikzpicture}[ baseline={([yshift=-.5ex]current bounding box.center)},vertex/.style={anchor=base,circle,fill=black!25,minimum size=18pt,inner sep=2pt}]
  \draw[name path=zigzag1] (0,10) to (3,9) to (0,8) to (3,7) to (0,6);
   %increase bounding box
   %\draw[opacity=0] (-0.2, -0.2) rectangle (18, 12.2);
   
   \draw[opacity=0,name path=t1] (0, 10) to (0,5);
   \fill [name intersections={of=t1 and zigzag1, by={a,blank, i,q}}]
        (a) circle (2pt)
        (i) circle (2pt)
        (q) circle (2pt);
\node[label=left:{$a$}] at (a) {};
\node[label=left:{$i$}] at (i) {};
\node[label=left:{$q$}] at (q) {};

\draw[opacity=0,name path=t2] (0.5, 10) to (0.5,6);
   \fill [name intersections={of=t2 and zigzag1, by={b,h, j, p}}]
        \foreach \x in {b, h, p}
        {(\x) circle (2pt)}
        (j) circle (2pt); 
      \node[above=2pt] at (b) {$bc$};
      \node[above=2pt] at (h) {$gh$};
      \node[above=2pt] at (p) {$op$};
        \node[below] at (j) {$jk$};

%\draw[opacity=0,name path=t3] (1.5, 10) to (1.5,6);
   %\fill [name intersections={of=t3 and zigzag1, by={c,g, k, o}}]
      %  \foreach \x in {c,g, k, o}
        %{(\x) circle (2pt)};
        %\foreach \x in {c,g, k, o}
        %{\node[label=above:{$\x$}] at (\x) {};}
        
\draw[opacity=0,name path=t4] (2, 10) to (2,6);
   \fill [name intersections={of=t4 and zigzag1, by={d,f, l, n}}]
        \foreach \x in {d,f, l, n}
        {(\x) circle (2pt)};
        \foreach \x in {d,f, l, n}
        {\node[label=above:{$\x$}] at (\x) {};}
        
\draw[opacity=0,name path=t5] (3, 10) to (3,6);
   \fill [name intersections={of=t5 and zigzag1, by={e,blank, m}}]
        \foreach \x in {e,m}
        {(\x) circle (2pt)};
        \foreach \x in {e,m}
        {\node[label=above:{$\x$}] at (\x) {};}
\end{tikzpicture}} $\pm \ldots \pm$ 
 { \begin{tikzpicture}[ baseline={([yshift=-.5ex]current bounding box.center)},vertex/.style={anchor=base,circle,fill=black!25,minimum size=18pt,inner sep=2pt}]
  \draw[name path=zigzag1] (0,10) to (3,9) to (0,8) to (3,7) to (0,6);
   %increase bounding box
   %\draw[opacity=0] (-0.2, -0.2) rectangle (18, 12.2);
   
   \draw[opacity=0,name path=t1] (0, 10) to (0,5);
   \fill [name intersections={of=t1 and zigzag1, by={a,blank, i,q}}]
        (a) circle (2pt)
        (i) circle (2pt)
        (q) circle (2pt);
\node[label=left:{$a$}] at (a) {};
\node[label=left:{$i$}] at (i) {};
\node[label=left:{$q$}] at (q) {};

\draw[opacity=0,name path=t2] (0.5, 10) to (0.5,6);
   \fill [name intersections={of=t2 and zigzag1, by={b,h, j, p}}]
        \foreach \x in {b, h, p}
        {(\x) circle (2pt)}
        (j) circle (2pt); 
    \foreach \x in {b, h, p}
        {\node[above=2pt] at (\x) {$\x$};}
        \node[below] at (j) {$j$};

\draw[opacity=0,name path=t3] (1.5, 10) to (1.5,6);
   \fill [name intersections={of=t3 and zigzag1, by={c,g, k, o}}]
        \foreach \x in {c,g, k, o}
        {(\x) circle (2pt)};
        \foreach \x in {c,g, k, o}
        {\node[label=above:{$\x$}] at (\x) {};}
        
%\draw[opacity=0,name path=t4] (2, 10) to (2,6);
   %\fill [name intersections={of=t4 and zigzag1, by={d,f, l, n}}]
      %  \foreach \x in {d,f, l, n}
        %{(\x) circle (2pt)};
        %\foreach \x in {d,f, l, n}
        %{\node[label=above:{$\x$}] at (\x) {};}
        
\draw[opacity=0,name path=t5] (3, 10) to (3,6);
   \fill [name intersections={of=t5 and zigzag1, by={e,blank, m}}]
        \foreach \x in {e,m}
        {(\x) circle (2pt)};      
        \node[label=above:{$def$}] at (e) {};
        \node[label=above:{$lmn$}] at (m) {};
\end{tikzpicture}}}
  \end{equation*}

\begin{definition}\label{def CH^ZZ}
Let $(\A, \cdot, d)$ be a (possibly non-commutative) associative, unital DGA over a fixed commutative ring $S$.  We define the zigzag Hochschild complex, $CH^{ZZ}(\A) = \bigoplus\limits_{n, k \ge 0, k \text{ even}}   (\A \otimes(\A^{\otimes n}\otimes \A)^{\otimes k})[n]$, where $[n]$ denotes a total shift down by $n$, with a differential defined below.  Monomials  in  $CH^{ZZ}(\A)$  will be written as
\abovedisplayskip=5pt
 \belowdisplayskip=5pt
 \begin{align*}
 \underline{x} = x^{\mathcal{L}}&\otimes (x_{(1,1)} \otimes \ldots \otimes x_{(1,n)} \otimes x^{\mathcal{R}}_1) \otimes\\
  \ldots &\otimes (x_{(i,1)} \otimes \ldots \otimes x_{(i,n)} \otimes x^{\dagger_i}_i)\otimes \\
  \ldots &\otimes (x_{(k,1)} \otimes \ldots \otimes x_{(k,n)} \otimes x^{\mathcal{L}}_k)
 \end{align*}

where $\dagger_i$ is $\mathcal{R}$ if $i$ is odd (i.e. on a ``zig'') and is equal to $\mathcal{L}$ if $i$ is even (i.e. on a ``zag'').  The differential $D: CH^{ZZ}(\A)  \to CH^{ZZ}(\A) $ is given by 
$D( \underline{x}_{k,n}):= (d + b)(  \underline{x}_{k,n})$.  The two components $d$ and $b$ are defined below: 
\begin{align*}
d( \underline{x}_{k,n}):= &(-1)^n d(x^{\mathcal{L}}) \otimes  \ldots \otimes (x_{(i,1)} \otimes \ldots \otimes x_{(i,n)} \otimes x^{\dagger_i}_i)\otimes \ldots(\dots  \otimes x^{\mathcal{L}}_k)\\
+ & \sum\limits_{p=1}^n \sum\limits_{i=1}^k (-1)^{n+ \beta_{(i,p)}} x^{\mathcal{L}} \otimes  \ldots \otimes (x_{(i,1)} \otimes \ldots \otimes d(x_{(i, p)}) \otimes \ldots \otimes x^{\dagger_i}_i) \otimes \ldots \otimes x^{\mathcal{L}}_k \\
+&  \sum\limits_{i=1}^{k} (-1)^{n+\beta_i^{\dagger_i}} x^{\mathcal{L}} \otimes  \ldots \otimes (x_{(i,1)} \otimes \ldots \otimes x_{(i,n)} \otimes d(x^{\dagger_i}_i)) \otimes \ldots  \otimes x^{\mathcal{L}}_k \end{align*}
where $\beta_i^{\dagger_i}$ and $\beta_{(i,p)}$ equal the sum of the degrees of the elements in $\A$ appearing before the element that $d$ is being applied to.  In particular $ \beta_i^{\dagger_i}:= |x^{\mathcal{L}}| + \ldots + |x_{(i,1)}| + \ldots |x_{(i,n)}|$ and $\beta_{(i,p)}:= |x^{\mathcal{L}}| + \ldots + |x_{(i,1)}| + \ldots +  |x_{(i,p-1)}|$ for $p>1$ and $\beta_{(i,1)}:= |x^{\mathcal{L}}| + \ldots + |x_{(i-1,1)}| + \ldots |x_{(i-1,n)}| + |x_i^{\dagger_i}|$. 
\begin{align*}
b( \underline{x}_{k,n}):&=  (-1)^n (x^{\mathcal{L}} \cdot x_{(1,1)}) \otimes (x_{(1,2)} \otimes \ldots \otimes x_{(1,n)} \otimes x^{\mathcal{R}}_1) \otimes \\   &\ \ldots \otimes (x_{(i,1)} \otimes \ldots \otimes x_{(i,n-1)} \otimes (x_{(i,n)} \cdot x^{\mathcal{L}}_i \cdot x_{(i+1, 1)})) \otimes \\  & \  \ldots \otimes  (x_{(k,1)} \otimes \ldots \otimes x_{(k,n-1)} \otimes (x_{(k, n)} \cdot x^{\mathcal{L}}_k))\\
			&+ \sum\limits_{p=1}^{n-1} (-1)^{n+p} x^{\mathcal{L}} \otimes (\ldots \otimes (x_{(1,p)}\cdot x_{(1, p+1)})\otimes \ldots \otimes x^{\mathcal{R}}_1) \otimes \\   	&\  \ldots \otimes (\ldots  \otimes(x_{(i,p)}\cdot x_{(i, p+1)}) \otimes \ldots \otimes x_i^{\dagger_i}) \otimes \\ 
			 &  \  \ldots \otimes  (\ldots \otimes(x_{(k,p+1)}\cdot x_{(k, p)}) \otimes \ldots \otimes x^{\mathcal{L}}_k)\\
			&+		x^{\mathcal{L}} \otimes (x_{(1,1)} \otimes \ldots \otimes x_{(1,n-1)} \otimes (x_{(1, n)} \cdot x^{\mathcal{R}}_1 \cdot x_{(2,1)})) \otimes \\ &  \  \ldots \otimes (x_{(i,1)}   \otimes \ldots \otimes x_{(i,n-1)} \otimes (x_{(i,n)} \cdot x^{\mathcal{R}}_i \cdot x_{(i+1, 1)}))  \otimes  \\ &  \ \ldots \otimes  (x_{(k,2)} \otimes \ldots \otimes x_{(k,n)} \otimes x^{\mathcal{L}}_k)
\end{align*}
For a pictorial representation of a simple monomial, see the figure on page \pageref{simple monomial pic}.  For a pictorial representation of $d$ and $b$, see the figures on page \pageref{dandbpics}.
\end{definition}

 In the definition above we are applying the following sign convention:  When $d$ is applied, each summand has a sign of $(-1)^{n+ \beta}$ where $n$ is the number of columns between endpoints in our zigzags and $\beta$ is the sum of the degrees of elements preceding the current element which $d$ is being applied to.  When $b$ is applied, each summand has a sign of $(-1)^{n+p}$ where we keep track of the fact that $b$ had to move over $p-1$ many columns and the $n-1$ is motivated by the eventual use of Stokes' formula for $\A= \Omega(M)$, $$d\int\limits_{F} \omega = (-1)^{dim{F}} \int\limits_{F} d\omega + (-1)^{dim(F) -1} \int\limits_{\partial F} \omega.$$
\begin{prop} \label{prop D^2=0 zigzag}
$D^2 = 0$
\begin{proof}

$D^2 = D \circ D = (d+ b)\circ (d + b) = d^2 + d \circ b - b \circ d  + b^2 = 0$.  We analyze each term independently:
\begin{itemize}
\item $d^2=0$ is a consequence of $\A$ being a DGA.
\item For $b^2=0$, we note that by associativity of ``$\cdot$'' in $\A$, we need only to check that we get opposite signs for the two ways in which the columns $j-1$, $j$, and $j+1$ come together.  It is easy to check that for when a collapse of the $j$-th and $j+1$-th columns is followed by a collapse of the $j-1$-th and $j$-th columns, we get a total sign of $(-1)^0$.  But when we collapse the $j-1$-th column with the $j$-th column, followed by collapsing the $j$-th and $j+1$-th columns, we obtain a sign of $(-1)^1$.  Similarly, terms in $b^2$ vanish when the columns being collapsed are separated.
\item For $d\circ b + b \circ d = 0$, we consider an element and focus on certain tensor factors\footnote{The calculation is shown along a ``zig'', the same result holds along a ``zag'' but the multiplication would be ordered differently.} of that element, namely: $$\underline{x}_{k,n} = \ldots \otimes x_{(i, p-1)}\otimes x_{(i, p)}\otimes x_{(i, p+1)}\otimes x_{(i, p+2)}\otimes \ldots. $$  Looking first at $d \circ b$ we obtain terms 
\begin{align*}
d \circ b (\underline{x}_{k,n}) &=  d \circ b (\ldots \otimes x_{(i, p-1)}\otimes x_{(i, p)}\otimes x_{(i, p+1)}\otimes x_{(i, p+2)}\otimes \ldots) \\
&= d( \ldots + (-1)^{n+p}(\ldots \otimes x_{(i, p-1)}\otimes x_{(i, p)}\cdot x_{(i, p+1)}\otimes x_{(i, p+2)}\otimes \ldots) +\ldots \\
&= \ldots + (-1)^{\beta_{(i,p-1)}+p-1}(\ldots \otimes d(x_{(i, p-1)})\otimes x_{(i, p)}\cdot x_{(i, p+1)}\otimes \ldots)+  \\
&\ldots + (-1)^{\beta_{(i,p)}+p-1}(\ldots \otimes d(x_{(i, p)}\cdot x_{(i, p+1)}) \otimes \ldots) + \\
& \ldots + (-1)^{\beta_{(i,p+2)}+p-1}(\ldots \otimes x_{(i, p)}\cdot x_{(i, p+1)}\otimes d(x_{(i, p+2)})\otimes \ldots)+ \ldots
\end{align*}
After applying the fact that $d$ is a derivation with respect to $\cdot$, along with $\beta_{(i,p)} + |x_{(i,p)}| = \beta_{(i, p+1)}$ we get
\begin{align*}
d \circ b (\underline{x}_{k,n}) &= \ldots + (-1)^{\beta_{(i,p-1)}+p-1}(\ldots \otimes d(x_{(i, p-1)})\otimes x_{(i, p)}\cdot x_{(i, p+1)}\otimes \ldots)+  \\
&\ldots + (-1)^{\beta_{(i,p)}+p-1}(\ldots \otimes d(x_{(i, p)})\cdot x_{(i, p+1)} \otimes \ldots) + \\
&\ldots + (-1)^{\beta_{(i,p+1)}+p-1 }(\ldots \otimes x_{(i, p)}\cdot d(x_{(i, p+1)}) \otimes \ldots) + \\
& \ldots + (-1)^{\beta_{(i,p+2)}+p-1}(\ldots \otimes x_{(i, p)}\cdot x_{(i, p+1)}\otimes d(x_{(i, p+2)})\otimes \ldots)+ \ldots 
\end{align*}
Next we consider $b \circ d$ and observe 
\begin{align*}
b \circ d (\underline{x}_{k,n}) &=  b \circ d (\ldots \otimes x_{(i, p-1)}\otimes x_{(i, p)}\otimes x_{(i, p+1)}\otimes x_{(i, p+2)}\otimes \ldots) \\
&= b ( \ldots + (-1)^{\beta_{(i,p-1)}+n}(\ldots \otimes d(x_{(i, p-1)})\otimes x_{(i, p)} \otimes \ldots)+  \\
&\ldots + (-1)^{\beta_{(i,p)}+ n }(\ldots \otimes d(x_{(i, p)}) \otimes x_{(i, p+1)} \otimes \ldots) + \\
&\ldots + (-1)^{\beta_{(i,p+1)}+n}(\ldots \otimes x_{(i, p)} \otimes d(x_{(i, p+1)}) \otimes \ldots) + \\
& \ldots + (-1)^{\beta_{(i,p+2)}+n}(\ldots \otimes x_{(i, p+1)}\otimes d(x_{(i, p+2)})\otimes \ldots)+ \ldots )\\
&=  \ldots + (-1)^{\beta_{(i,p-1)} + p}(\ldots \otimes d(x_{(i, p-1)}) \otimes x_{(i, p)} \cdot x_{(i, p+1)} \otimes \ldots)+  \\
&\ldots + (-1)^{\beta_{(i,p)} + p}(\ldots \otimes d(x_{(i, p)}) \cdot x_{(i, p+1)} \otimes \ldots) + \\
&\ldots + (-1)^{\beta_{(i,p+1)} + p}(\ldots \otimes x_{(i, p)} \cdot d(x_{(i, p+1)}) \otimes \ldots) + \\
& \ldots + (-1)^{\beta_{(i,p+2)}+ p}(\ldots \otimes  x_{(i, p+1)} \cdot x_{(i, p+1)} \otimes d(x_{(i, p+2)})\otimes \ldots)+ \ldots 
\end{align*}
By comparing the corresponding terms we see that $d \circ b + b\circ d=0$.  Similar arguments apply when considering other tensor-factors as well as when $b$ collapses the first or last columns; in which case one has to apply the derivation $d$ over a product of three tensor-factors.  
\end{itemize}
\end{proof}
\end{prop}

Next we define a shuffle product $\odot: CH^{ZZ}(\A)  \otimes CH^{ZZ}(\A)  \to CH^{ZZ}(\A) $.  Recall that for $CH^I(\A)$, we could define a shuffle product, but it would not be compatible with the usual Hochschild differential $D$ unless $\A$ was a commutative DGA.  For $CH^{ZZ}(\A)$, the idea is to concatenate the two monomials while shuffling in 1's like so:  \begin{equation*} 
{\resizebox{!}{3cm}
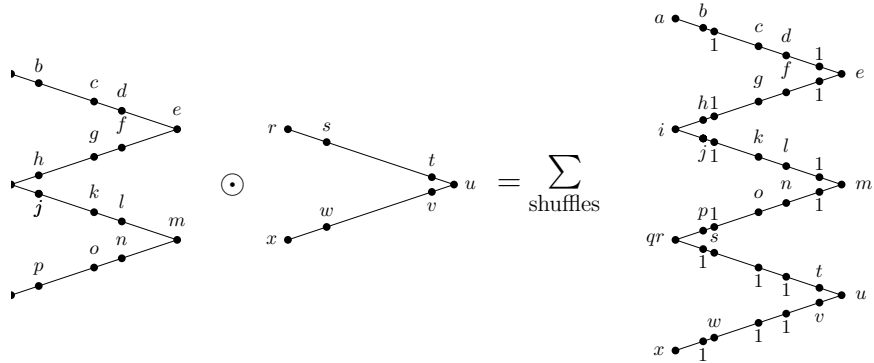
{\begin{tikzpicture}[baseline={([yshift=-.5ex]current bounding box.center)},vertex/.style={anchor=base,
    circle,fill=black!25,minimum size=18pt,inner sep=2pt}]
\clip(0,4) rectangle (16,12);
%\draw[help lines,xstep=1,ystep=1] (0,0) grid (20,20);
%\foreach \x in {0,1,...,20} { \node [anchor=north] at (\x,0) {\x}; }
%\foreach \y in {0,1,...,20} { \node [anchor=east] at (0,\y) {\y}; }
%
%
%zigzag 1
%
%
%    
   \draw[name path=zigzag1] (0,10) to (3,9) to (0,8) to (3,7) to (0,6);
   %increase bounding box
   \draw[opacity=0] (-0.2, -0.2) rectangle (18, 12.2);
   
   \draw[opacity=0,name path=t1] (0, 10) to (0,5);
   \fill [name intersections={of=t1 and zigzag1, by={a,blank, i,q}}]
        (a) circle (2pt)
        (i) circle (2pt)
        (q) circle (2pt);
\node[label=left:{$a$}] at (a) {};
\node[label=left:{$i$}] at (i) {};
\node[label=left:{$q$}] at (q) {};

\draw[opacity=0,name path=t2] (0.5, 10) to (0.5,6);
   \fill [name intersections={of=t2 and zigzag1, by={b,h, j, p}}]
        \foreach \x in {b, h, p}
        {(\x) circle (2pt)}
        (j) circle (2pt); 
    \foreach \x in {b, h, p}
        {\node[above=2pt] at (\x) {$\x$};}
        \node[below] at (j) {$j$};

\draw[opacity=0,name path=t3] (1.5, 10) to (1.5,6);
   \fill [name intersections={of=t3 and zigzag1, by={c,g, k, o}}]
        \foreach \x in {c,g, k, o}
        {(\x) circle (2pt)};
        \foreach \x in {c,g, k, o}
        {\node[label=above:{$\x$}] at (\x) {};}
        
\draw[opacity=0,name path=t4] (2, 10) to (2,6);
   \fill [name intersections={of=t4 and zigzag1, by={d,f, l, n}}]
        \foreach \x in {d,f, l, n}
        {(\x) circle (2pt)};
        \foreach \x in {d,f, l, n}
        {\node[label=above:{$\x$}] at (\x) {};}
        
\draw[opacity=0,name path=t5] (3, 10) to (3,6);
   \fill [name intersections={of=t5 and zigzag1, by={e,blank, m}}]
        \foreach \x in {e,m}
        {(\x) circle (2pt)};
        \foreach \x in {e,m}
        {\node[label=above:{$\x$}] at (\x) {};}
        
 %
%
%zigzag 2
%
%
%         
        
        \draw[name path=zigzag2] (5,9) to (8,8) to (5,7);
           \draw[opacity=0,name path=s1] (5, 10) to (5,5);
   \fill [name intersections={of=s1 and zigzag2, by={r, x}}]
        \foreach \x in {r, x}
        {(\x) circle (2pt)}; 
    \foreach \x in {r,x}
        {\node[left=2pt] at (\x) {$\x$};}

\draw[opacity=0,name path=s2] (5.7, 10) to (5.7,6);
   \fill [name intersections={of=s2 and zigzag2, by={s,w}}]
        \foreach \x in {s,w}
        {(\x) circle (2pt)}
        (j) circle (2pt); 
    \foreach \x in {s,w}
        {\node[above=2pt] at (\x) {$\x$};}
        \node[below] at (j) {$j$};
        
   \draw[opacity=0,name path=s3] (7.6, 10) to (7.6,6);
   \fill [name intersections={of=s3 and zigzag2, by={t,v}}]
        \foreach \x in {t}
        {(\x) circle (2pt)}
        (v) circle (2pt);
    \foreach \x in {t}
        {\node[above=2pt] at (\x) {$\x$};}
        \node[below=2pt] at (v) {$v$};

      \draw[opacity=0,name path=s4] (8, 10) to (8,6);
   \fill [name intersections={of=s4 and zigzag2, by={blank, u}}]
        \foreach \x in {u}
        {(\x) circle (2pt)};
    \foreach \x in {u}
        {\node[right=2pt] at (\x) {$\x$};}

%
%
%zigzag 3
%
%
%      
        
          \draw[name path=zigzag3] (12,11) to (15,10) to (12, 9) to (15, 8) to (12, 7) to (15, 6) to (12, 5);
   %increase bounding box

   \draw[opacity=0,name path=u1] (12, 12) to (12,4);
   \fill [name intersections={of=u1 and zigzag3, by={a,i, blank, qr, blank, x}}]
        \foreach \x in {a,i, qr, x}
        {(\x) circle (2pt)}; 
    \foreach \x in {a, i, qr, x}
        {\node[left=2pt] at (\x) {$\x$};}

\draw[opacity=0,name path=u2] (12.5, 12) to (12.5,4);
   \fill [name intersections={of=u2 and zigzag3, by={b,h, j, p, blank1, blank2}}]
        \foreach \x in {b, h, j, p, blank1, blank2}
        {(\x) circle (2pt)}
        (j) circle (2pt); 
    \foreach \x in {b, h, p}
        {\node[above=2pt] at (\x) {$\x$};}
        \node[below] at (j) {$j$};
        \node[below] at (blank1) {$1$};
        \node[below] at (blank2) {$1$};

\draw[opacity=0,name path=u3] (13.5, 12) to (13.5,4);
   \fill [name intersections={of=u3 and zigzag3, by={c,g, k, o, blank1, blank2}}]
        \foreach \x in {c,g, k, o, blank1, blank2}
        {(\x) circle (2pt)};
        \foreach \x in {c,g, k, o}
        {\node[label=above:{$\x$}] at (\x) {};}
         \node[below] at (blank1) {$1$};
         \node[below] at (blank2) {$1$};
        
\draw[opacity=0,name path=u4] (14, 12) to (14,4);
   \fill [name intersections={of=u4 and zigzag3, by={d,f, l, n, blank1, blank2}}]
        \foreach \x in {d,f, l, n, blank1, blank2}
        {(\x) circle (2pt)};
        \foreach \x in {d,f, l, n}
        {\node[label=above:{$\x$}] at (\x) {};}
        \node[below] at (blank1) {$1$};
        \node[below] at (blank2) {$1$};

      \draw[opacity=0,name path=u5] (12.7, 12) to (12.7,4);
   \fill [name intersections={of=u5 and zigzag3, by={blank1, blank2, blank3, blank4, s,w}}]
        \foreach \x in {blank1, blank2, blank3, blank4, s,w}
        {(\x) circle (2pt)} ;
    \foreach \x in {s,w}
        {\node[above=2pt] at (\x) {$\x$};}
       \node[below] at (blank1) {$1$};
       \node[above] at (blank2) {$1$};
       \node[below] at (blank3) {$1$};
       \node[above] at (blank4) {$1$};

        \draw[opacity=0,name path=u6] (14.6, 12) to (14.6,4);
   \fill [name intersections={of=u6 and zigzag3, by={blank1, blank2, blank3, blank4, t,v}}]
        \foreach \x in {blank1, blank2, blank3, blank4, t, v}
        {(\x) circle (2pt)};
    \foreach \x in {t}
        {\node[above=2pt] at (\x) {$\x$};}
        \node[below=2pt] at (v) {$v$};
        \node[above] at (blank1) {$1$};
        \node[below] at (blank2) {$1$};
        \node[above] at (blank3) {$1$};
        \node[below] at (blank4) {$1$};
        
\draw[opacity=0,name path=u7] (15, 12) to (15,4);
   \fill [name intersections={of=u7 and zigzag3, by={e,blank, m, blank, u}}]
        \foreach \x in {e,m, u}
        {(\x) circle (2pt)};
        \foreach \x in {e,m, u}
        {\node[label=right:{$\x$}] at (\x) {};}
   
 \node at (4,8) {\huge $\odot$};
 \node at (9,8) {\huge $=$};
 \node at (10, 8) {\huge $\sum\limits_{\text{shuffles}}$};
 \end{tikzpicture}}} 
  \end{equation*}
\begin{remark}[Some comments on shuffles and their signatures] \label{shuffle remarks}
Consider the set $S_{n,m}$ of $(n,m)$ shuffles:  $S_{n,m} = \{\sigma \in S_{n+m} | \sigma(1) < \ldots < \sigma(n), \ \text{ and, } \ \sigma(n+1) < \ldots < \sigma(n+m)\}$.  We can use these shuffles to get maps $\sigma: \A^{\otimes n} \times \A^{\otimes m} \to \A^{\otimes n+m}$ via $$\sigma((a_1 \otimes \ldots \otimes a_n),( a_{n+1} \otimes \ldots \otimes a_{n+m})) : = a_{\sigma^{-1}(1)}   \otimes \ldots \otimes a_{\sigma^{-1}(n+m)}.$$Next, since we will be shuffling in $1$'s on zigs and zags in the opposite order, we offer a formal convention to describe that process.  For an element $E = e_1 \otimes \ldots \otimes e_p \in \A^{\otimes p}$, we define $\overline{E}:= e_{p} \otimes \ldots \otimes e_1 \in \A^{\otimes p}$.  Now, for an shuffle $\sigma \in S_{n,m}$ interpreted as $\sigma: \A^{\otimes n} \times \A^{\otimes m} \to \A^{\otimes n+m}$, we can define $\sigma^{\Sh} \in S_{n,m}$ giving the induced map $\sigma^{\Sh}: \A^{\otimes n} \times \A^{\otimes m} \to \A^{\otimes n+m}$ to be $$\sigma^{\Sh}(L,R):= \overline{\sigma( \overline{L}, \overline{R} )}.$$  For a given $i\in \mathbb{N}$ we will define $\sigma^{\Sh_i}$ to be $\sigma$ if $i$ is odd (i.e. on a ``zig") and $\sigma^{\Sh}$ if $i$ is even (i.e. on a ``zag").  Finally, recall that for a shuffle (or permutation of any type) we can define the signature of that shuffle by $sgn(\sigma) = (-1)^{N(\sigma)}$ where $N(\sigma)$ is the number of transpositions needed to write $\sigma$ as a product purely of transpositions.  Note that $N(\sigma)$ is well defined $mod$ $2$.  By our convention, for reasons having nothing to do with zigs and zags, we will need to consider $sgn(\sigma^{\Sh})$.  It is straightforward to prove that $sgn(\sigma^{\Sh}) = (-1)^{nm}sgn(\sigma)$ where $\sigma$ shuffles an $n$-tuple with an $m$-tuple.
\end{remark}
\begin{definition}\label{zz shuffle}The shuffle product $\odot$ defined for $CH^{ZZ}(\A)$ is defined by 
\begin{align*}
&\underline{x}_{k, n} \odot \underline{y}_{\l, m}\\  := &\left( x^{\mathcal{L}}\otimes \ldots \otimes ((x_{(i,1)} \otimes x_{(i,2)} \otimes \ldots \otimes x_{(i,n)} )\otimes x^{\dagger_i}_i) \otimes \ldots \otimes x^{\mathcal{L}}_{k} \right)\\
&\odot \left( y^{\mathcal{L}}\otimes \ldots \otimes ((y_{(j,1)} \otimes y_{(j,2)} \otimes \ldots \otimes y_{(j,m)} )\otimes y^{\dagger_j}_j) \otimes \ldots \otimes y^{\mathcal{L}}_l \right) \\ 
:= &\sum\limits_{\sigma \in S_{n,m}} (-1)^{\epsilon_{\sigma}}x^{\mathcal{L}}\otimes \ldots \otimes \left(  \sigma^{\Sh_i} ( (x_{(i,1)} \otimes \ldots \otimes x_{(i,n)}),(1, \ldots, 1) )\otimes x^{\dagger_i}_i)  \right)\otimes \ldots  \\ &\otimes (x^{\mathcal{L}}_k \cdot y^{\mathcal{L}})\otimes  \ldots \otimes \left(\sigma^{\Sh_j} ( (1, \ldots, 1),(y_{(j,1)} \otimes \ldots \otimes y_{(j,m)})) \otimes y^{\dagger_j}_j  \right) \otimes \ldots \otimes y_l^{\mathcal{L}}
\end{align*} 
where $\epsilon_{\sigma}:= |\underline{x}_{k,n}|\cdot m + N(\sigma^{\Sh})= (|\underline{x}_{k,n}|+n)\cdot m + N(\sigma)$, using the abbreviation $|\underline{x_{k,n}}|:= |x^{\mathcal{L}}| + |x_{(1,1)}| + \ldots + |x_k^{\mathcal{L}}|$, and the shuffles are happening at all $i,j$ simultaneously for a given $\sigma$.  
\end{definition}

\begin{prop}\label{propzzDderivation}
The shuffle product $\odot$ is associative and $D$ is a derivation with respect to $\odot$.
\begin{proof}
Associativity comes from the fact that placing zigzags on top of one another is an associative operation.  We wish to show that $D(\underline{x} \odot \underline{y}) = D(\underline{x}) \odot \underline{y} + (-1)^{|\underline{x}| +n}\underline{x} \odot D(\underline{y})$.  We note that $D(\underline{x}) \odot \underline{y}$ amounts to applying the differential only to the top zigzag where as $\underline{x} \odot D(\underline{y})$ applies the differential only to the bottom zigzag.  $D(\underline{x} \odot \underline{y})$ applies first $d$ to each term, in which case you are either on the top or the bottom zigzag and so recovering those terms is straightforward.  Applying $b$ to $\underline{x} \odot \underline{y}$ has some terms which vanish and some terms which cancel with $\underline{x} \odot D(\underline{y})$ or $D(\underline{x}) \odot \underline{y}$ as we will now demonstrate.  Consider the two diagrams below:\\  
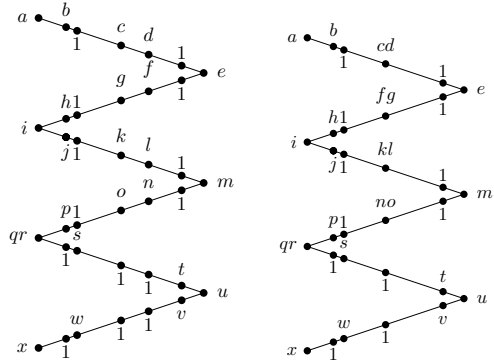
\begin{table}[H]
\centering
\begin{tabular}{cc}
\resizebox{\textwidth/4}{!}{\begin{tikzpicture}         
%
%
%zigzag 3
%
%
%      
\draw[name path=zigzag3] (12,11) to (15,10) to (12, 9) to (15, 8) to (12, 7) to (15, 6) to (12, 5);
   %increase bounding box

   \draw[opacity=0,name path=u1] (12, 12) to (12,4);
   \fill [name intersections={of=u1 and zigzag3, by={a,i, blank, qr, blank, x}}]
        \foreach \x in {a,i, qr, x}
        {(\x) circle (2pt)}; 
    \foreach \x in {a, i, qr, x}
        {\node[left=2pt] at (\x) {$\x$};}

\draw[opacity=0,name path=u2] (12.5, 12) to (12.5,4);
   \fill [name intersections={of=u2 and zigzag3, by={b,h, j, p, blank1, blank2}}]
        \foreach \x in {b, h, j, p, blank1, blank2}
        {(\x) circle (2pt)}
        (j) circle (2pt); 
    \foreach \x in {b, h, p}
        {\node[above=2pt] at (\x) {$\x$};}
        \node[below] at (j) {$j$};
        \node[below] at (blank1) {$1$};
        \node[below] at (blank2) {$1$};

\draw[opacity=0,name path=u3] (13.5, 12) to (13.5,4);
   \fill [name intersections={of=u3 and zigzag3, by={c,g, k, o, blank1, blank2}}]
        \foreach \x in {c,g, k, o, blank1, blank2}
        {(\x) circle (2pt)};
        \foreach \x in {c,g, k, o}
        {\node[label=above:{$\x$}] at (\x) {};}
         \node[below] at (blank1) {$1$};
         \node[below] at (blank2) {$1$};
        
\draw[opacity=0,name path=u4] (14, 12) to (14,4);
   \fill [name intersections={of=u4 and zigzag3, by={d,f, l, n, blank1, blank2}}]
        \foreach \x in {d,f, l, n, blank1, blank2}
        {(\x) circle (2pt)};
        \foreach \x in {d,f, l, n}
        {\node[label=above:{$\x$}] at (\x) {};}
        \node[below] at (blank1) {$1$};
        \node[below] at (blank2) {$1$};

      \draw[opacity=0,name path=u5] (12.7, 12) to (12.7,4);
   \fill [name intersections={of=u5 and zigzag3, by={blank1, blank2, blank3, blank4, s,w}}]
        \foreach \x in {blank1, blank2, blank3, blank4, s,w}
        {(\x) circle (2pt)} ;
    \foreach \x in {s,w}
        {\node[above=2pt] at (\x) {$\x$};}
       \node[below] at (blank1) {$1$};
       \node[above] at (blank2) {$1$};
       \node[below] at (blank3) {$1$};
       \node[above] at (blank4) {$1$};

        \draw[opacity=0,name path=u6] (14.6, 12) to (14.6,4);
   \fill [name intersections={of=u6 and zigzag3, by={blank1, blank2, blank3, blank4, t,v}}]
        \foreach \x in {blank1, blank2, blank3, blank4, t, v}
        {(\x) circle (2pt)};
    \foreach \x in {t}
        {\node[above=2pt] at (\x) {$\x$};}
        \node[below=2pt] at (v) {$v$};
        \node[above] at (blank1) {$1$};
        \node[below] at (blank2) {$1$};
        \node[above] at (blank3) {$1$};
        \node[below] at (blank4) {$1$};
        
\draw[opacity=0,name path=u7] (15, 12) to (15,4);
   \fill [name intersections={of=u7 and zigzag3, by={e,blank, m, blank, u}}]
        \foreach \x in {e,m, u}
        {(\x) circle (2pt)};
        \foreach \x in {e,m, u}
        {\node[label=right:{$\x$}] at (\x) {};}
 \end{tikzpicture}}
& \resizebox{\textwidth/4}{!}{ \begin{tikzpicture}           
%
%
%zigzag 3
%
%
%      
\draw[name path=zigzag3] (12,11) to (15,10) to (12, 9) to (15, 8) to (12, 7) to (15, 6) to (12, 5);
   %increase bounding box

   \draw[opacity=0,name path=u1] (12, 12) to (12,4);
   \fill [name intersections={of=u1 and zigzag3, by={a,i, blank, qr, blank, x}}]
        \foreach \x in {a,i, qr, x}
        {(\x) circle (2pt)}; 
    \foreach \x in {a, i, qr, x}
        {\node[left=2pt] at (\x) {$\x$};}

\draw[opacity=0,name path=u2] (12.5, 12) to (12.5,4);
   \fill [name intersections={of=u2 and zigzag3, by={b,h, j, p, blank1, blank2}}]
        \foreach \x in {b, h, j, p, blank1, blank2}
        {(\x) circle (2pt)}
        (j) circle (2pt); 
    \foreach \x in {b, h, p}
        {\node[above=2pt] at (\x) {$\x$};}
        \node[below] at (j) {$j$};
        \node[below] at (blank1) {$1$};
        \node[below] at (blank2) {$1$};

\draw[opacity=0,name path=u3] (13.5, 12) to (13.5,4);
   \fill [name intersections={of=u3 and zigzag3, by={cd,fg, kl, no, blank1, blank2}}]
        \foreach \x in {cd,fg, kl, no, blank1, blank2}
        {(\x) circle (2pt)};
        \foreach \x in {cd,fg, kl, no}
        {\node[label=above:{$\x$}] at (\x) {};}
         \node[below] at (blank1) {$1$};
         \node[below] at (blank2) {$1$};

      \draw[opacity=0,name path=u5] (12.7, 12) to (12.7,4);
   \fill [name intersections={of=u5 and zigzag3, by={blank1, blank2, blank3, blank4, s,w}}]
        \foreach \x in {blank1, blank2, blank3, blank4, s,w}
        {(\x) circle (2pt)} ;
    \foreach \x in {s,w}
        {\node[above=2pt] at (\x) {$\x$};}
       \node[below] at (blank1) {$1$};
       \node[above] at (blank2) {$1$};
       \node[below] at (blank3) {$1$};
       \node[above] at (blank4) {$1$};

        \draw[opacity=0,name path=u6] (14.6, 12) to (14.6,4);
   \fill [name intersections={of=u6 and zigzag3, by={blank1, blank2, blank3, blank4, t,v}}]
        \foreach \x in {blank1, blank2, blank3, blank4, t, v}
        {(\x) circle (2pt)};
    \foreach \x in {t}
        {\node[above=2pt] at (\x) {$\x$};}
        \node[below=2pt] at (v) {$v$};
        \node[above] at (blank1) {$1$};
        \node[below] at (blank2) {$1$};
        \node[above] at (blank3) {$1$};
        \node[below] at (blank4) {$1$};
        
\draw[opacity=0,name path=u7] (15, 12) to (15,4);
   \fill [name intersections={of=u7 and zigzag3, by={e,blank, m, blank, u}}]
        \foreach \x in {e,m, u}
        {(\x) circle (2pt)};
        \foreach \x in {e,m, u}
        {\node[label=right:{$\x$}] at (\x) {};}
 \end{tikzpicture} }
 \end{tabular}
 \end{table}
The diagram on the left is a term in the shuffle product of two zigzags where two columns in the top zigzag remain adjacent.  After applying the $b$-component of our differential $D$ to this new zigzag where we collapse the columns of the elements $c$ and $d$, we obtain the diagram on the right.  This term can will cancel with the term coming from $D(\underline{x}) \odot \underline{y}$ which first collapses that column in the zigzag $\underline{x}$ and then shuffles the result in with $\underline{y}$ in such a way that the  zigzag on the right (above) is obtained.  A similar argument would show were we see cancelation with terms in $\underline{x} \odot D(\underline{y})$  On the other hand, consider the three diagrams below:
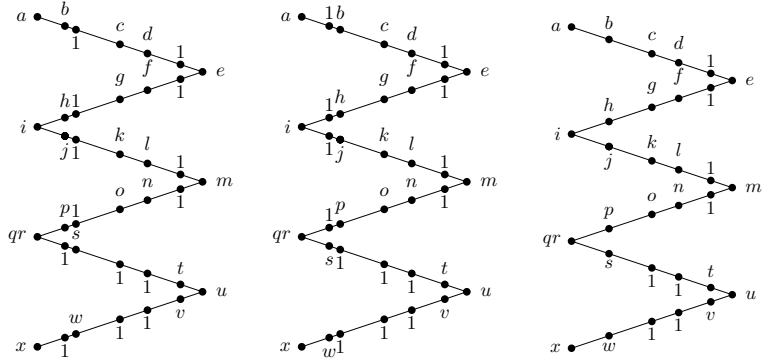
\begin{table}[H]
\centering
\begin{tabular}{ccc}
\resizebox{\textwidth/4}{!}{\begin{tikzpicture}         
%
%
%zigzag 3
%
%
%      
\draw[name path=zigzag3] (12,11) to (15,10) to (12, 9) to (15, 8) to (12, 7) to (15, 6) to (12, 5);
   %increase bounding box

   \draw[opacity=0,name path=u1] (12, 12) to (12,4);
   \fill [name intersections={of=u1 and zigzag3, by={a,i, blank, qr, blank, x}}]
        \foreach \x in {a,i, qr, x}
        {(\x) circle (2pt)}; 
    \foreach \x in {a, i, qr, x}
        {\node[left=2pt] at (\x) {$\x$};}

\draw[opacity=0,name path=u2] (12.5, 12) to (12.5,4);
   \fill [name intersections={of=u2 and zigzag3, by={b,h, j, p, blank1, blank2}}]
        \foreach \x in {b, h, j, p, blank1, blank2}
        {(\x) circle (2pt)}
        (j) circle (2pt); 
    \foreach \x in {b, h, p}
        {\node[above=2pt] at (\x) {$\x$};}
        \node[below] at (j) {$j$};
        \node[below] at (blank1) {$1$};
        \node[below] at (blank2) {$1$};

\draw[opacity=0,name path=u3] (13.5, 12) to (13.5,4);
   \fill [name intersections={of=u3 and zigzag3, by={c,g, k, o, blank1, blank2}}]
        \foreach \x in {c,g, k, o, blank1, blank2}
        {(\x) circle (2pt)};
        \foreach \x in {c,g, k, o}
        {\node[label=above:{$\x$}] at (\x) {};}
         \node[below] at (blank1) {$1$};
         \node[below] at (blank2) {$1$};
        
\draw[opacity=0,name path=u4] (14, 12) to (14,4);
   \fill [name intersections={of=u4 and zigzag3, by={d,f, l, n, blank1, blank2}}]
        \foreach \x in {d,f, l, n, blank1, blank2}
        {(\x) circle (2pt)};
        \foreach \x in {d,f, l, n}
        {\node[label=above:{$\x$}] at (\x) {};}
        \node[below] at (blank1) {$1$};
        \node[below] at (blank2) {$1$};

      \draw[opacity=0,name path=u5] (12.7, 12) to (12.7,4);
   \fill [name intersections={of=u5 and zigzag3, by={blank1, blank2, blank3, blank4, s,w}}]
        \foreach \x in {blank1, blank2, blank3, blank4, s,w}
        {(\x) circle (2pt)} ;
    \foreach \x in {s,w}
        {\node[above=2pt] at (\x) {$\x$};}
       \node[below] at (blank1) {$1$};
       \node[above] at (blank2) {$1$};
       \node[below] at (blank3) {$1$};
       \node[above] at (blank4) {$1$};

        \draw[opacity=0,name path=u6] (14.6, 12) to (14.6,4);
   \fill [name intersections={of=u6 and zigzag3, by={blank1, blank2, blank3, blank4, t,v}}]
        \foreach \x in {blank1, blank2, blank3, blank4, t, v}
        {(\x) circle (2pt)};
    \foreach \x in {t}
        {\node[above=2pt] at (\x) {$\x$};}
        \node[below=2pt] at (v) {$v$};
        \node[above] at (blank1) {$1$};
        \node[below] at (blank2) {$1$};
        \node[above] at (blank3) {$1$};
        \node[below] at (blank4) {$1$};
        
\draw[opacity=0,name path=u7] (15, 12) to (15,4);
   \fill [name intersections={of=u7 and zigzag3, by={e,blank, m, blank, u}}]
        \foreach \x in {e,m, u}
        {(\x) circle (2pt)};
        \foreach \x in {e,m, u}
        {\node[label=right:{$\x$}] at (\x) {};}
 \end{tikzpicture}}
& \resizebox{\textwidth/4}{!}{\begin{tikzpicture}         
%
%
%zigzag 3
%
%
%      
\draw[name path=zigzag3] (12,11) to (15,10) to (12, 9) to (15, 8) to (12, 7) to (15, 6) to (12, 5);
   %increase bounding box

   \draw[opacity=0,name path=u1] (12, 12) to (12,4);
   \fill [name intersections={of=u1 and zigzag3, by={a,i, blank, qr, blank, x}}]
        \foreach \x in {a,i, qr, x}
        {(\x) circle (2pt)}; 
    \foreach \x in {a, i, qr, x}
        {\node[left=2pt] at (\x) {$\x$};}

\draw[opacity=0,name path=u2] (12.7, 12) to (12.7,4);
   \fill [name intersections={of=u2 and zigzag3, by={b,h, j, p, blank1, blank2}}]
        \foreach \x in {b, h, j, p, blank1, blank2}
        {(\x) circle (2pt)}
        (j) circle (2pt); 
    \foreach \x in {b, h, p}
        {\node[above=2pt] at (\x) {$\x$};}
        \node[below] at (j) {$j$};
        \node[below] at (blank1) {$1$};
        \node[below] at (blank2) {$1$};

\draw[opacity=0,name path=u3] (13.5, 12) to (13.5,4);
   \fill [name intersections={of=u3 and zigzag3, by={c,g, k, o, blank1, blank2}}]
        \foreach \x in {c,g, k, o, blank1, blank2}
        {(\x) circle (2pt)};
        \foreach \x in {c,g, k, o}
        {\node[label=above:{$\x$}] at (\x) {};}
         \node[below] at (blank1) {$1$};
         \node[below] at (blank2) {$1$};
        
\draw[opacity=0,name path=u4] (14, 12) to (14,4);
   \fill [name intersections={of=u4 and zigzag3, by={d,f, l, n, blank1, blank2}}]
        \foreach \x in {d,f, l, n, blank1, blank2}
        {(\x) circle (2pt)};
        \foreach \x in {d,f, l, n}
        {\node[label=above:{$\x$}] at (\x) {};}
        \node[below] at (blank1) {$1$};
        \node[below] at (blank2) {$1$};

      \draw[opacity=0,name path=u5] (12.5, 12) to (12.5,4);
   \fill [name intersections={of=u5 and zigzag3, by={blank1, blank2, blank3, blank4, s,w}}]
        \foreach \x in {blank1, blank2, blank3, blank4, s,w}
        {(\x) circle (2pt)} ;
    \foreach \x in {s,w}
        {\node[below=2pt] at (\x) {$\x$};}
       \node[above] at (blank1) {$1$};
       \node[above] at (blank2) {$1$};
       \node[below] at (blank3) {$1$};
       \node[above] at (blank4) {$1$};

        \draw[opacity=0,name path=u6] (14.6, 12) to (14.6,4);
   \fill [name intersections={of=u6 and zigzag3, by={blank1, blank2, blank3, blank4, t,v}}]
        \foreach \x in {blank1, blank2, blank3, blank4, t, v}
        {(\x) circle (2pt)};
    \foreach \x in {t}
        {\node[above=2pt] at (\x) {$\x$};}
        \node[below=2pt] at (v) {$v$};
        \node[above] at (blank1) {$1$};
        \node[below] at (blank2) {$1$};
        \node[above] at (blank3) {$1$};
        \node[below] at (blank4) {$1$};
        
\draw[opacity=0,name path=u7] (15, 12) to (15,4);
   \fill [name intersections={of=u7 and zigzag3, by={e,blank, m, blank, u}}]
        \foreach \x in {e,m, u}
        {(\x) circle (2pt)};
        \foreach \x in {e,m, u}
        {\node[label=right:{$\x$}] at (\x) {};}
 \end{tikzpicture}}
& \resizebox{\textwidth/4}{!}{ \begin{tikzpicture}         
%
%
%zigzag 3
%
%
%      
\draw[name path=zigzag3] (12,11) to (15,10) to (12, 9) to (15, 8) to (12, 7) to (15, 6) to (12, 5);
   %increase bounding box

   \draw[opacity=0,name path=u1] (12, 12) to (12,4);
   \fill [name intersections={of=u1 and zigzag3, by={a,i, blank, qr, blank, x}}]
        \foreach \x in {a,i, qr, x}
        {(\x) circle (2pt)}; 
    \foreach \x in {a, i, qr, x}
        {\node[left=2pt] at (\x) {$\x$};}

\draw[opacity=0,name path=u2] (12.7, 12) to (12.7,4);
   \fill [name intersections={of=u2 and zigzag3, by={b,h, j, p, s, w}}]
        \foreach \x in {b, h, j, p, s, w}
        {(\x) circle (2pt)}
        (j) circle (2pt); 
    \foreach \x in {b, h, p}
        {\node[above=2pt] at (\x) {$\x$};}
        \node[below] at (j) {$j$};
        \node[below] at (s) {$s$};
        \node[below] at (w) {$w$};

\draw[opacity=0,name path=u3] (13.5, 12) to (13.5,4);
   \fill [name intersections={of=u3 and zigzag3, by={c,g, k, o, blank1, blank2}}]
        \foreach \x in {c,g, k, o, blank1, blank2}
        {(\x) circle (2pt)};
        \foreach \x in {c,g, k, o}
        {\node[label=above:{$\x$}] at (\x) {};}
         \node[below] at (blank1) {$1$};
         \node[below] at (blank2) {$1$};
        
\draw[opacity=0,name path=u4] (14, 12) to (14,4);
   \fill [name intersections={of=u4 and zigzag3, by={d,f, l, n, blank1, blank2}}]
        \foreach \x in {d,f, l, n, blank1, blank2}
        {(\x) circle (2pt)};
        \foreach \x in {d,f, l, n}
        {\node[label=above:{$\x$}] at (\x) {};}
        \node[below] at (blank1) {$1$};
        \node[below] at (blank2) {$1$};

        \draw[opacity=0,name path=u6] (14.6, 12) to (14.6,4);
   \fill [name intersections={of=u6 and zigzag3, by={blank1, blank2, blank3, blank4, t,v}}]
        \foreach \x in {blank1, blank2, blank3, blank4, t, v}
        {(\x) circle (2pt)};
    \foreach \x in {t}
        {\node[above=2pt] at (\x) {$\x$};}
        \node[below=2pt] at (v) {$v$};
        \node[above] at (blank1) {$1$};
        \node[below] at (blank2) {$1$};
        \node[above] at (blank3) {$1$};
        \node[below] at (blank4) {$1$};
        
\draw[opacity=0,name path=u7] (15, 12) to (15,4);
   \fill [name intersections={of=u7 and zigzag3, by={e,blank, m, blank, u}}]
        \foreach \x in {e,m, u}
        {(\x) circle (2pt)};
        \foreach \x in {e,m, u}
        {\node[label=right:{$\x$}] at (\x) {};}
 \end{tikzpicture}}
\end{tabular}
\end{table}where the left and middle diagrams are zigzags coming from different shuffles (off by a transposition and thus off by exactly $(-1)$) and the right diagram is a common term shared in the result of the $b$-component of the differential $D$ applied to these zigzags.  Thus, cancelation occurs.  
\end{proof} 
\end{prop}

\begin{remark}
We have that $(CH^{ZZ}(\A), D, \odot)$ is a DGA and thus we can iterate the $CH^{ZZ}$ functor to define the DGA $CH^{ZZ}(CH^{ZZ}(\ldots(CH^{ZZ}(\A)) \ldots))$.
\end{remark}

\subsection{Special Cases}
We begin by showing the Hochschild complex $CH^{ZZ}(\A)$ is compatible with $CH^{I}(\A)$ in the case where $\A$ is a commutative DGA.  The other special case we are interested in is when our DGA comes from a commutative DGA tensored with an associative algebra.  This case is motivated by considering matrix-valued differential forms $\Omega^{\bullet}(M, Mat) \cong \Omega^{\bullet}(M, \mathbb{R}) \otimes Mat$.  
\begin{definition}\label{def zz to I}
If $\C$ is a commutative DGA, we have a column-collapse map $Col: CH^{ZZ}(\C) \to CH^I(\C)$ defined as follows: Let $\underline{x}_{(k,n)} \in CH^{ZZ}(\C)$.  Then 
$$Col(\underline{x}_{(k,n)} ) = (-1)^{\epsilon} Col^{\mathcal{L}}(\underline{x}_{(k,n)} ) \otimes Col^1(\underline{x}_{(k,n)} ) \otimes \ldots \otimes Col^n(\underline{x}_{(k,n)} ) \otimes Col^{\mathcal{R}}(\underline{x}_{(k,n)} )$$
where 
\begin{align*}
Col^{\mathcal{L}}(\underline{x}_{(k,n)} ):= & x^{\mathcal{L}} \cdot x_2^{\mathcal{L}} \cdot \dots \cdot x_k^{\mathcal{L}}\\
Col^p(\underline{x}_{(k,n)} ):=  &\prod\limits_{\substack{ i=1 \\ i \text{ odd}}}^{k-1} x_{(i,p)} \cdot \prod\limits_{\substack{ i=2 \\ i \text{ even}}}^{k} x_{(i,n-p+1)} \\
Col^{\mathcal{R}}(\underline{x}_{(k,n)} ):= & x_1^{\mathcal{R}} \cdot x_3^{\mathcal{R}} \cdot \dots \cdot x_{k-1}^{\mathcal{R}}\\
\end{align*}
Here by $\prod$ we refer to the ordered product induced by the algebra $\C$ and $\epsilon$ comes from the usual Koszul rule of changing the order of elements $x_{(i,p)}$.
\end{definition}

\begin{prop}\label{prop zz to I}
Let $\C$ be a commutative differential graded algebra.  Then $Col: CH^{ZZ}(\C) \to CH^I(\C)$ is a chain map and an algebra map. 
\begin{proof}
It is straightforward to check that the differentials are compatible because of Leibniz and because we always collapse full rows.  Similarly the shuffle products are compatible since the insertion of $1$'s in an entire column amounts to the usual shuffle after collapsing.
\end{proof}
\end{prop}

Let $\mathcal{C}$ be a commutative DGA and $\mathcal{B}$ an associative algebra.  Recall we have the associative DGA $\mathcal{C} \otimes \mathcal{B}$ generated by monomials $c \otimes b \in \mathcal{C}^n \otimes \mathcal{B}$  with differential $d: \mathcal{C}^n \otimes \mathcal{B} \to \mathcal{C}^{n+1} \otimes \mathcal{B}$ given by $d(c\otimes b):= d(c) \otimes b$ and associative product $(c \otimes b) \cdot (c' \otimes b'):= (-1)^{|b|\cdot |c'|} (c\cdot c' \otimes b \cdot b')$ which yields the product $(\mathcal{C} \otimes \mathcal{B})^{\otimes n} \times   (\mathcal{C} \otimes \mathcal{B})^{\otimes n} \to (\mathcal{C} \otimes \mathcal{B})^{\otimes n}$.  With all of this in mind we can consider the Hochschild complex $CH^I(\C \otimes \B)$ and define a special shuffle product for it, where the idea is to shuffle the commutative part and push all of the information from $\B$ to the end, while preserving the order.

\begin{definition}\label{matrixshuffle}
For $\C$ a commutative DGA and $\B$ an associative algebra, $CH^I(\C \otimes \B)$ has the shuffle product:
\begin{align*}
(\underline{x \otimes \omega}) \odot (\underline{y \otimes \nu})
:= & \sum\limits_{\sigma \in S_{n,m}} (-1)^{\epsilon_{\sigma}}\sigma\left( (\underline{x \otimes 1}), (\underline{y \otimes 1})\right)\\
\cdot & \left((1\otimes 1) \otimes \ldots \otimes (1\otimes \omega^{\mathcal{L}} \cdot \ldots \cdot \omega_i \cdot \ldots \omega^{\mathcal{R}} \cdot \nu^{\mathcal{L}} \cdot \ldots \cdot \nu_i  \cdot \ldots \cdot \nu^{\mathcal{R}}) \right)
\end{align*}
where \begin{align*}
(\underline{x \otimes \omega})&:= (x^{\mathcal{L}} \otimes \omega^{\mathcal{L}}) \otimes (x_1 \otimes \omega_1) \otimes  \ldots \otimes (x_n \otimes \omega_n) \otimes (x^{\mathcal{R}} \otimes \omega^{\mathcal{R}}) \in (\mathcal{C} \otimes \mathcal{B})^{\otimes n}\\
\text{and} \quad (\underline{y \otimes \nu})&:= (y^{\mathcal{L}} \otimes \nu^{\mathcal{L}}) \otimes (y_1 \otimes \nu_1) \otimes  \ldots \otimes (y_n \otimes \nu_n) \otimes (y^{\mathcal{R}} \otimes \nu^{\mathcal{R}}) \in (\mathcal{C} \otimes \mathcal{B})^{\otimes n}
\end{align*}
\end{definition}

\begin{prop}
For $\C$ a commutative DGA and $\B$ an associative algebra, the differential $D:= d+b$ on $CH^I(\C \otimes \B)$ is a derivation of the shuffle product defined in Definition \ref{matrixshuffle}.
\begin{proof}
The proof is similar to the proof for Proposition \ref{propzzDderivation}. 
\end{proof}
\end{prop}

The special case for $\C \otimes \B$ we are considering can also fall into the case where we consider it as a single (non-commutative) DGA with unit, and thus we could also define $CH^{ZZ}(\C \otimes \B)$.  The column-collapse map in this case is almost the same as in Definition \ref{def zz to I} but this time we push all of the elements from $\B$ to the end, preserving order.
\begin{definition}
If $\C$ is a commutative DGA and $\B$ is an associative algebra, we have a column-collapse map $Col_{Mat}: CH^{ZZ}(\C \otimes \B) \to CH^I(\C \otimes \B)$ defined as follows:

$$Col_{Mat}(\underline{x\otimes m}_{(k,n)} ) = (-1)^{\epsilon} Col_{Mat}^{\mathcal{L}}(\underline{x\otimes m}_{(k,n)} ) \otimes  \ldots \otimes  Col_{Mat}^{\mathcal{R}}(\underline{x\otimes m}_{(k,n)} )$$
where 
\begin{align*}
Col_{Mat}^{\mathcal{L}}(\underline{x\otimes m}_{(k,n)} ):= & (x^{\mathcal{L}} \cdot x_2^{\mathcal{L}} \cdot \dots \cdot x_k^{\mathcal{L}})\otimes 1\\
Col_{Mat}^p(\underline{x\otimes m}_{(k,n)} ):=  &\left( \prod\limits_{\substack{ i=1 \\ i \text{ odd}}}^{k-1} x_{(i,p)} \cdot \prod\limits_{\substack{ i=2 \\ i \text{ even}}}^{k} x_{(i,n-p+1)}\right) \otimes 1 \\
Col_{Mat}^{\mathcal{R}}(\underline{x\otimes m}_{(k,n)} ):= & \left( x_1^{\mathcal{R}} \cdot x_3^{\mathcal{R}} \cdot \dots \cdot x_{k-1}^{\mathcal{R}}\right) \otimes ( m^{\mathcal{L}} \cdot \prod m_{(i,p)}) \\
\end{align*}
where again by $\prod$ we mean the ordered product induced by $\C$ or $\B$ and $\epsilon$ again comes from the Koszul rule.
\end{definition}

\begin{prop}
Given a commutative DGA, $\C$, and an associative algebra $\B$, the column-collapse map $Col_{Mat}: CH^{ZZ}(\C \otimes \B) \to CH^I(\C \otimes \B)$ is a chain map and an algebra map.
\begin{proof}
The proof is similar to the one for Proposition \ref{prop zz to I}.
\end{proof}
\end{prop}

 \section{A Chen Map out of the zigzag Hochschild complex}
We use $CH^{ZZ}(\A)$ to model non-abelian differential forms on the path space and so for the remainder of the paper, when we are considering matrix-valued differential forms, we write $\OmegaMat(M):=  \Omega(M, Mat)$, so as to distinguish from our real-valued forms $\Omega(M):= \Omega(M, \mathbb{R})$.  First let us recall the usual interval Hochschild model.  Let $\underline{\omega }_n \in CH^I(\Omega(M))$.  We have an evaluation map 
\begin{align*}
\Delta^n \times PM &\xrightarrow{ev} M^{n+2}, \\
((t_1, \ldots, t_n), \gamma) &\xmapsto{ev} (\gamma(0), \gamma(t_1), \ldots, \gamma(t_n), \gamma(1))
\end{align*}
 and so we use the composition
 $$\Omega(M)^{\otimes n+2} \xrightarrow{} \Omega(M^{n+2}) \xrightarrow{ev^*} \Omega(\Delta^n \times PM) \xrightarrow{\int_{\Delta^n}} \Omega(PM)$$
 to define $It(\underline{\omega}) := \int_{\Delta^n} ev^*(\underline{\omega})$.  Similarly we can use the evaluation map governed by placing differential forms along a zigzag at each time-slot and can define an iterated integral map out of $CH^{ZZ}(\OmegaMat(M))$:
\begin{definition}\label{zigzag ev}
Let $\underline{\omega }_{(n,k)} \in CH^{ZZ}(\OmegaMat(M))$.  We have an evaluation map 
\begin{align*}
\Delta^n \times PM \xrightarrow{ev} &M^{nk + k+1}, \\
((t_1, \ldots, t_n), \gamma) \xmapsto{ev} &(\gamma(0), \gamma(t_1), \ldots, \gamma(t_n), \gamma(1), \\
& \gamma(t_n), \ldots \gamma(t_1), \gamma(0), \\
\ldots & \gamma(t_n), \ldots \gamma(t_1), \gamma(0))
\end{align*}
   and so we use the composition
 $$\OmegaMat(M)^{\otimes nk + k+1} \rightarrow \OmegaMat(M^{nk +k+1}) \xrightarrow{ev^*} \OmegaMat(\Delta^n \times PM) \xrightarrow{\int_{\Delta^n}} \OmegaMat(PM)$$
 to define $It(\underline{\omega}) := \int_{\Delta^n} ev^*(\underline{\omega})$.  The evaluation map here can be clarified by the figure on page \pageref{general zigzag}.
 \end{definition}
When applying the differential in the image of $CH^{ZZ}(\OmegaMat(M)) \xrightarrow{It} \OmegaMat(PM)$, we encounter the situation where two time-slots in $\Delta^n$ come together.  For this reason, we recall the following Lemma which will be applied below without further reference:
 \begin{lemma}\label{lemma:diagonalwedge}
 The map $\OmegaMat(M) \otimes \OmegaMat(M) \xrightarrow{EZ} \OmegaMat(M \times M) \xrightarrow{\Delta^*} \OmegaMat(M)$, where $EZ$ is the Eilenberg-Zilber map and $\Delta:M\times M \to M$ is the diagonal, is given by the wedge product of forms. 
 \end{lemma}
 
We now arrive at our first important result and continue by establishing the relationship between $CH^I(\Omega(M))$, $CH^{ZZ}(\OmegaMat(M))$, and their Chen maps.
\begin{prop}\label{prop zz chain map}
The Iterated Integral $It: CH^{ZZ}(\OmegaMat(M)) \to \OmegaMat(PM)$ is a chain map.
\begin{proof}
We have $d_{DR} \int\limits_{\Delta^n} ev^*(\underline{\omega}) = (-1)^n \int\limits_{\Delta^n} d_{DR} ev^*(\underline{\omega})+ (-1)^{n-1} \int\limits_{\partial \Delta^n} ev^*(\underline{\omega})$.  Now we use: (a) the chain map $\OmegaMat(M)^{\otimes k} \to \OmegaMat(M^k)$ along with the exterior derivative acting as a derivation and (b) the commutative diagram 
\begin{equation*}
\xymatrix{\Delta^{n-1} \times PM \ar[d]_{d_i \times id} \ar[r]^{ev} & M^{(n-1)k + k+1 } \ar[d]^{d_i} \\
 \Delta^n \times PM  \ar[r]^{ev} & M^{nk+k+1} \\  }
\end{equation*}
where $d_i: \Delta^{n-1} \to \Delta$ is the map  $(t_1, \ldots , t_i, \ldots t_{n-1}) \mapsto (t_1, \ldots , t_i, t_i, \ldots t_{n-1})$ and $d_i$ on the right arrow is, by abuse of notation, the diagonal making the diagram commute.  So then the above equation can be expressed as:
\begin{align*}
d_{DR} \int\limits_{\Delta^n} ev^*(\underline{\omega}) 
&=   (-1)^n \int\limits_{\Delta^n} d_{DR} ev^*(\underline{\omega})   +  (-1)^{n-1} \int\limits_{\partial\Delta^{n}} ev^*(\underline{\omega}) \\
&=   (-1)^n \int\limits_{\Delta^n}  ev^*(d_{DR}\underline{\omega})   + (-1)^{n-1}  \sum\limits_i \int\limits_{\Delta^{n-1}} (d_i \times id)^*(ev^*(\underline{\omega}) )\\
&= (-1)^n \sum\limits_{p} (-1)^{\beta_p} \int\limits_{\Delta^n}  ev^*(\ldots \otimes d_{DR}(\omega)\otimes \ldots) \\ &+(-1)^{n-1}  \sum\limits_i (-1)^{i-1}  \int\limits_{\Delta^{n-1}} ev^*(\ldots \otimes (\omega \cdot \omega') \otimes \ldots  )\\
&= \int\limits_{\Delta^n}  ev^*(d (\underline{\omega})) +  \int\limits_{\Delta^{n-1}} ev^*(b (\underline{\omega}) )\\
&= It (D(\underline{\omega}))
\end{align*}
\end{proof}
\end{prop}

\begin{remark}
By a similar calculation, we have that $It: CH^I(\Omega(M)) \to \Omega(PM)$ is a chain map for any coefficients since we never have to commute forms in the differential of $CH^I(\Omega(M))$.
\end{remark}
We want to show that $It: CH^{ZZ}(\OmegaMat(M) \to \OmegaMat(PM)$ is also an algebra map.  As a warm-up, we first recall why in the abelian case, $It: CH^{I}(\Omega(M, \mathbb{R})) \to \Omega(PM, \mathbb{R})$ is an algebra map.  
\begin{prop}\label{Prop interval algebra map}
The Iterated Integral $CH^I(\Omega(M)) \to \Omega(PM)$ is a map of algebras.
\begin{proof}
We have degeneracy maps $s_i: \Delta^{r+1} \to \Delta^r, (t_1, \ldots, t_r) \mapsto (t_1, \ldots \hat{t_i}, \ldots, t_r)$ yielding 
\begin{equation*}
\xymatrix{
\Delta^{r+1} \ar[d]^{s_i} \ar[r]^{ev} & M^{r+1+2}\ar[d]^{s_i} \\
\Delta^r  \ar[r]^{ev}       &M^{r+2} }
\end{equation*}
where $s_i$ on the right arrow is the induced projection making the diagram commute.  So by composing these degeneracy maps (and again abusing notation) we obtain for a fixed shuffle $\sigma \in S_{n,m}$ the commutative diagram:
\begin{equation*}
\xymatrix{\Delta^{n+m} \times PM \ar[d]^{\beta^{\sigma} \times id} \ar[rr]^{ev_{n+m}} &&  M^{m+n+2} \ar[d]^{\rho^{\sigma}} \\
\Delta^m \times \Delta^n \times PM \ar[rr]^{(ev_n, ev_m)} & &M^{m+2} \times M^{n+2}  }
\end{equation*}
Here $\beta^{\sigma} : \Delta^{n+m} \to \Delta^n \times \Delta^m$ is the unshuffle map $$(t_1, \ldots , t_{n+m}) \xmapsto{\beta^{\sigma}} ((t_{\sigma(1)}, \ldots , t_{\sigma(n)}), (t_{\sigma_{(n+1)}}, \ldots , t_{\sigma_{(n+m)}}))$$ and $\rho^{\sigma}$ on the right is the map which makes the diagram commute.  Now we claim that $$\int\limits_{\Delta^m} ev^*(\underline{\omega}) \wedge \int\limits_{\Delta^n} ev^*(\underline{\nu}) = 
\int\limits_{\Delta^{m+n}} ev^*( \underline{\omega} \odot \underline{\nu})$$
We begin proving this claim by noting that we can evaluate the left hand side 
\begin{align*}
&\int\limits_{\Delta^n} ev_n^*(\underline{\omega}) \wedge \int\limits_{\Delta^m} ev_m^*(\underline{\nu}) \\  & = (-1)^{(|\underline{\omega}| + n)\cdot m} \int\limits_{\Delta^m \times \Delta^n} (ev_n, ev_m)^*(\underline{\omega} \otimes \underline{\nu})\\
&= (-1)^{(|\underline{\omega}| + n)\cdot m} \int\limits_{\coprod\limits_{\sigma \in S_{n,m}} \beta^{\sigma}(\Delta^{m+n})} (ev_n^*, ev_m^*)(\underline{\omega} \otimes \underline{\nu})\\
&=  \sum\limits_{\sigma \in S_{n,m}}  (-1)^{(|\underline{\omega}| + n)\cdot m}\cdot sgn(\sigma) \int\limits_{\Delta^{m+n}} (\beta^{\sigma} \times id)^* (ev_n^*, ev_m^*)(\underline{\omega} \otimes \underline{\nu})\\
&= \int\limits_{\Delta^{m+n}} ev_{n+m}^*( \sum\limits_{\sigma \in S_{n,m}}  (-1)^{(|\underline{\omega}| + n)\cdot m}\cdot sgn(\sigma) (\rho^{\sigma})^* (\underline{\omega} \otimes \underline{\nu}))\\
&= \int\limits_{\Delta^{n+m}} ev_{n+m}^* (\underline{\omega} \odot \underline{\nu})
\end{align*}
\end{proof}
\end{prop}
\begin{remark}
In the above proof, the reason why we need abelian coefficients can be seen by the fact that the map $\rho^{\sigma}: M^{m+n+2} \to M^{m+2} \times M^{n+2}$ switches the order of coordinates.  In the following proposition and proof, no coordinates are switched, which is exactly why we can correctly model the wedge product of two Iterated Integrals with our shuffle product. 
\end{remark}
\begin{prop}\label{zz It algebra map}
The Iterated Integral $CH^{ZZ}(\OmegaMat(M)) \to \OmegaMat(PM)$ is a map of algebras.
\begin{proof}
The proof is analogous to the above proposition, so we include the commutative diagram in this case for clarity: 
\begin{equation*}
\xymatrix{\Delta^{n+m} \times PM \ar[d]^{\beta^{\sigma} \times id} \ar[rrr]^{ev_{n+m, k}} & &  & M^{(m+n+1)(k_n + k_m)+1} \ar[d]^{\pi^{\sigma}} \\
\Delta^m \times \Delta^n \times PM \ar[rrr]^{(ev_{m,k_m}, ev_{n,k_n})} & & &M^{(m+1)k_m +1} \times M^{(m+1)k_m +1}  }
\end{equation*}
where $\beta^{\sigma}$ is the same as in Proposition \ref{Prop interval algebra map}, $\pi^{\sigma}$ is the map which makes the diagram commute, and the evaluation maps $ev_{\bullet, \bullet}$ are the ones defined in Definition \ref{zigzag ev}.  Then 
\begin{align*}
&\int\limits_{\Delta^n} ev_{n, k_n}^*(\underline{\omega}) \wedge \int\limits_{\Delta^m} ev_{m, k_m}^*(\underline{\nu})\\ &=  \sum\limits_{\sigma \in S_{n,m}} (-1)^{(|\underline{\omega}| + n)\cdot m}\cdot sgn(\sigma) \int\limits_{\Delta^{n+m}} (\beta^{\sigma} \times id)^*(ev_{n, k_n}, ev_{m, k_m})^*(\underline{\omega} \otimes \underline{\nu})\\
&=  \int\limits_{\Delta^{n+m}} ev_{n+m, k_n+k_m}^* (\sum\limits_{\sigma}  (-1)^{(|\underline{\omega}| + n)\cdot m}\cdot sgn(\sigma) (\pi^{\sigma})^*(\underline{\omega} \otimes \underline{\nu}))\\
&= \int\limits_{\Delta^{n+m}} ev_{n+m, k_n+k_m}^*(\underline{\omega} \odot \underline{\nu})
\end{align*}
\end{proof}
\end{prop}

\begin{lemma}
The collapse map $Col: CH^{ZZ}(\Omega(M, \mathbb{R})) \to CH^I(\Omega(M, \mathbb{R}))$ from Definition \ref{def zz to I}  is given by $zz^*$ where $zz:M^{n+2} \to M^{nk+k+1}$ is a particular diagonal map.  In other words, we have the commutative diagram: 
\begin{equation*}
\xymatrix{CH^{ZZ}(\Omega(M)) \ar[d]_{Col} \ar[r]& \Omega(M^{nk+k+1}) \ar[d]^{zz^*} \\
CH^I(\Omega(M)) \ar[r]& \Omega(M^{n+2})}
\end{equation*}
\begin{proof}
This fact is another consequence of Lemma \ref{lemma:diagonalwedge}.
\end{proof}
\end{lemma}

\begin{prop}\label{prop Col It commute}
If $\Omega(M)= \Omega(M, \mathbb{R})$ then we have the following commutative diagram of DGAs.
\begin{equation*}
\xymatrix{CH^{ZZ}(\Omega(M)) \ar[d] \ar[r] & \Omega(PM)\\
CH^I(\Omega(M)) \ar[ru] &  }
\end{equation*}
\begin{proof}
We observe the following commutative diagram
\begin{equation*}
\xymatrix{M^{nk + k+1} & \Delta^n \times PM   \ar[l]_{ev^{ZZ}} \ar[ld]^{ev^I} \\
M^{n+2} \ar[u]^{zz} &  }
\end{equation*}
combined with the lemma above we obtain our result from the diagram
 \begin{equation*}
\xymatrix{CH^{ZZ}(\Omega(M)) \ar[d]_{Col} \ar[r] & \Omega(M^{nk+k+1})  \ar[d]^{zz^*} \ar[r]^{(ev^{ZZ})^*} &\Omega(\Delta^n \times PM) \ar[r]^{\int_{\Delta^n}} & \Omega(PM)\\
CH^I(\Omega(M)) \ar[r]& \Omega(M^{n+2}) \ar[ru]_{(ev^{I})^*} & & }
\end{equation*}
\end{proof}
\end{prop}

\begin{remark}
In the case of non-abelian coefficients, $\OmegaMat(M) = \Omega(M, Mat)$, both $It: CH^{ZZ}(\OmegaMat(M)) \to \OmegaMat(PM)$ and $It: CH^{I}(\OmegaMat(M)) \to \OmegaMat(PM)$ are chain maps, while $Col: CH^{ZZ}(\OmegaMat(M)) \to CH^I(\OmegaMat(M))$ is not a chain map.  However, in this case the diagram from Propostition \ref{prop Col It commute} does not commute either.
\end{remark}

\section{The two-dimensional zigzag Hochschild complex(es)}
Given a commutative DGA, $(\A, d, \cdot)$, we have the Hochschild complex of the Hochschild complex, the rectangular Hochschild complex, and the square Hochschild complex, denoted $CH^I(CH^I(\A))$, $CH^{I \times I}_{Rec}(\A)$, and $CH^{I\times I}_{Sq}(\A)$, respectively.  Note that for a non-commutative DGA, $(\A, d, \cdot)$, the 2-dimensional (both simplicial and bisimplicial) higher Hochschild structures do not form a complex \cite{GTZ}.  In this section we give two related 2-d Hochschild complexes using the ideas from $CH^{ZZ}(\A)$.\\

Given an associative DGA, $(\A, d, \cdot)$, we have the Hochschild complex of the zigzag Hochschild complex, denoted $CH^I(CH^{ZZ}(\A))$, and below we define the rectangular zigzag Hochschild complex and the square zigzag Hochschild complex, denoted $CH^{ZZ}_{Rec}(\A)$ and $CH^{ZZ}_{Sq}$, respectively.  
\begin{definition}\label{def rect zig zag}
Let $(\A, d, \cdot)$ be a DGA.  The rectangular zigzag Hochschild complex has underlying vector space

$$CH^{ZZ}_{Rec}(\A): = \bigoplus\limits_{\substack{m,n \ge 0 \\ k_0, \ldots, k_{m+1} \ge 0 \\ k_i \text{ even}}} (\A \otimes (\A^{\otimes n} \otimes \A)^{\otimes {k_0}} \otimes \ldots \otimes \A \otimes (\A^{\otimes n} \otimes \A)^{\otimes {k_{m+1}}})[n+m] $$ Monomials, $(\underline{x}^j_{k_j, n})_{j=0}^{m+1}$, are to be thought of as $m+2$-many rows of elements $\underline{x}^j_{k_j, n} \in CH^{ZZ}(\A)$, 
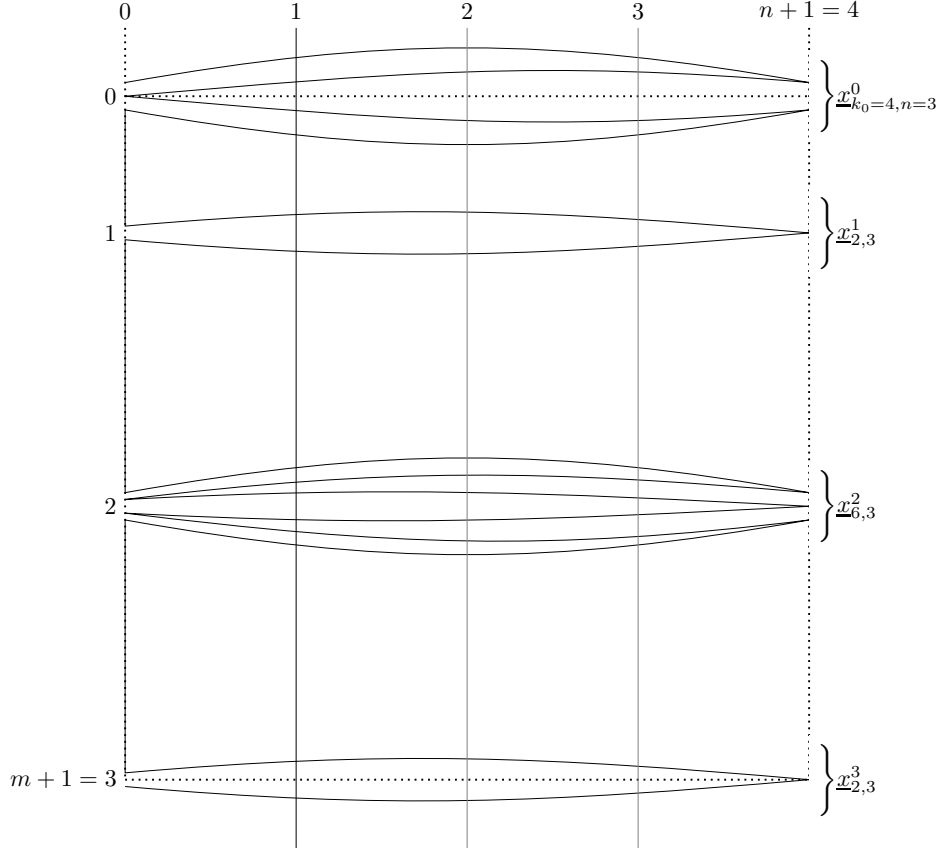
\begin{equation*}\resizebox{\textwidth}{!}{
\begin{tikzpicture}
%help lines grid
%\draw[help lines] (0.1,0.1) grid (9.9,9.9);

% s axis
\draw[thick][dotted] (0,11) -- (0,0);

% t axis
\draw[thick][dotted](0,10) -- (10,10);

%label axes
\node[above] at (0,11) {$0$};
\draw[thick][dotted] (10,11) -- (10,0);
\node[above] at(10,11) {$n+1=4$};
\draw[thick][dotted] (0,0) -- (10,0);
\node[left] at(0,0) {$m+1=3$};
\node[left] at (0,10) {$0$};

% s=0 zigzags
\draw[name path=zigzag] (0,10.2) to [out=10, in=170] (10,10.2) to [out=175, in=5](0,10)
to [out=-5, in=185] (10,9.8) to [out=190, in=-10] (0,9.8)
% s=s_1 zigzazags
-- (0,8.1) to [out=5, in=175] (10,8) to [out=185, in=-5](0,7.9)

%s=s_2 zigzags
--(0, 4.2) to  [out=10, in=170] (10,4.2) to [out=175, in=7](0,4.1)
to [out=3, in=177] (10,4) to [out=183, in=-3] (0,3.9)
to [out=-7, in=187] (10,3.8) to [out=190, in=-10](0,3.8)

%s=1 zigzags
-- (0,0.1) to [out=5, in=175] (10,0) to [out=185, in=-5](0,-0.1);

%label si's
\node[left] at (0,8) {$1$};
\node[left] at (0,4) {$2$};

%label ki's
\node[right, fill=white] at (10,10) { $\Bigg\} \underline{x}^0_{k_0=4,n=3}$};
\node[right, fill=white] at (10,8) { $\Bigg\} \underline{x}^1_{2,3}$};
\node[right, fill=white] at (10,4) { $\Bigg\} \underline{x}^2_{6,3}$};
\node[right, fill=white] at (10,0) { $\Bigg\} \underline{x}^3_{2,3}$};

%t-lines

%t1
\draw[name path=t1] (2.5,-1) -- (2.5,11);
\node[above] at (2.5,11) {$1$};

%t2
\draw[thin][gray] (5,-1) -- (5,11);
\node[above] at (5,11) {$2$};

%t3
\draw[thin][gray] (7.5,-1) -- (7.5,11);
\node[above] at (7.5,11) {$3$};

%\fill [name intersections={of=t1 and zigzag, by={a,b}}]
   %     (a) circle (2pt)
      %  (b) circle (2pt);
%\node[label=above:$r$] at (a) {};

\end{tikzpicture}}
\end{equation*}
with differential 
\begin{align*}
&D((\underline{x}^0_{k_0, n})\otimes \ldots \otimes (\underline{x}^{m+1}_{k_{m+1}, n}))\\
:=& \sum\limits_{r=0, p=0}^{m+1,n+1} (-1)^{n + m + \beta_{r,p}} (\underline{x}^0_{k_0, n})\otimes \ldots \otimes d_p((\underline{x}^r_{k_r, n})) \otimes \ldots \otimes (\underline{x}^{m+1}_{k_{m+1}, n}) \\
+ & \sum\limits_{p=0}^{n} (-1)^{m+n+p}\  b_p((\underline{x}^0_{k_0, n})\otimes \ldots \otimes (\underline{x}^{m+1}_{k_{m+1}, n}))\\
+&\sum\limits_{r=0}^{m} (-1)^{m+r} (\underline{x}^0_{k_0, n})\otimes \ldots \otimes((\underline{x}^r_{k_r, n})\star(\underline{x}^{r+1}_{k_{r+1}, n})) \otimes \ldots \otimes (\underline{x}^{m+1}_{k_{m+1}, n})
\end{align*}
where $d_p$ is the differential coming from $\A$ applied to exactly the $p$-th slot of 
$(\underline{x}_{k_i, n})$ and $b_p$ is collapsing/multiplying the slots in the $p$-th and $(p+1)$-th columns of $(\underline{x}_{k_i, n})$.  The operation $\star$ is simply a concatenation of the two zigzags:
\begin{align*}
\underline{a}_{k, n}\star \underline{b}_{l, n}:= &a^{\mathcal{L}}\otimes (a_{(1,1)} \otimes \ldots \otimes a_{(1,n)} \otimes a^{\mathcal{R}}_1) \otimes\\
  \ldots &\otimes (a_{(i,1)} \otimes \ldots \otimes a_{(i,n)} \otimes a^{\dagger_i}_i)\otimes \\
  \ldots &\otimes a_{(k,1)} \otimes \ldots \otimes a_{(k,n)} \\
  &\otimes (a^{\mathcal{L}}_k\cdot b^{\mathcal{L}})\otimes (b_{(1,1)} \otimes \ldots \otimes b_{(1,n)} \otimes b^{\mathcal{R}}_1) \otimes\\
  \ldots &\otimes (b_{(l,1)} \otimes \ldots \otimes b_{(l,n)} \otimes b_{l}^{\mathcal{L}})
\end{align*}

\begin{prop}
Let $(\A, d, \cdot)$ be an associative $DGA$, then $CH^{ZZ}_{Rec}(\A)$ forms a complex.
\begin{proof}
The components $d$ and $b$ of $D$ are just as in Proposition \ref{prop D^2=0 zigzag}.  For the new $\star$ part, it is straightforward to check that $\star d + d \star =0$, $\star^2=0$, and $\star b + b\star=0$.
\end{proof}
\end{prop}

\end{definition}
In 1-d holonomy, we integrate over a path in $M$, and now we would like to define our surfaces of integration for 2-d holonomy.  While we are not interested in this paper in many of the properties one normally wants in their space of bigons for 2-d holonomy, we acknowledge that notation developed in \cite{BaSc}. 
\begin{definition}
Given a manifold, $M$, we define the space of smooth bigons $BM$ to be the space of smooth maps $\Gamma: [0,1]^2 \to M$.
\end{definition}
\begin{remark}
While bigons usually require that the squares being mapped in are constant along the vertical edges, we will not require this until the last chapter, even though we will continue to refer to them as bigons in both contexts.
\end{remark}
%: iterated integral for CH^ZZ_Rec
\begin{definition}\label{It for CH^ZZ_Rec}
Let $(\underline{\omega}^j_{k_j, n})_{j=0}^{m+1} \in CH^{ZZ}_{Rec}(\OmegaMat(M))$.  We define the iterated integral, $CH^{ZZ}(\OmegaMat(M)) \xrightarrow{It} \OmegaMat(BM)$, by
$$It((\underline{\omega}^j_{k_j,n})_{j=0}^{m+1} ) := \int\limits_{\Delta^n \times \Delta^m} ev_{\underline{k},n}^*((\underline{\omega}^j_{k_j, n})_{i=0}^{m+1} )$$
where we use the evaluation map, with $\underline{k}:=( k_0, k_1, \ldots k_{m+1})$, 
\begin{align*}
&\Delta^n \times \Delta^n \times B(M) \xrightarrow{ev_{\underline{k},n}} M^{\sum\limits_{i=0}^{m+1} nk_i + k_i +1}\\
&(t_1, \ldots, t_n, s_1, \ldots, s_m, \Sigma)\\
 \mapsto &(\Sigma(0,0), \Sigma(0,t_1), \ldots, \Sigma(0, t_n), \Sigma(0, 1), \Sigma(0, t_n), \ldots , \Sigma(0,0),\ldots,  \\
              & \Sigma(s_1,0), \Sigma(s_1,t_1), \ldots, \Sigma(s_1, t_n), \Sigma(s_1, 1) \ldots, \Sigma(s_1,0), \ldots ,\\
									&\vdots\\
									& \Sigma(s_m,0), \Sigma(s_m,t_1), \ldots, \Sigma(s_m, t_n), \Sigma(s_m, 1),\ldots, \Sigma(s_m,0), \ldots, \\
									& \Sigma(1,0), \Sigma(1,t_1), \ldots, \Sigma(1, t_n), \Sigma(1, 1), \ldots, \Sigma(1,0))\\						
\end{align*}
\end{definition}

While the $d_p$ and $b_p$ above will be similar to the pieces of the 1-d zigzag story, something should be said about the $\star$ component of the differential.  Consider a very simple monomial (pictured below) on which we'd like to observe the star operation, say two zig zags with $n=2$, each having $k=1$: i.e. $1 \otimes \underline{x}^1_{(1,2)} \otimes \underline{x}^2_{(1,2)} \otimes 1$ with $a =1$, $\underline{x}_{(1,2)} = b \otimes c \otimes \ldots \otimes h$, $\underline{x'}_{(1,2)} = i \otimes j \otimes \ldots \otimes o$, and $p=1$.

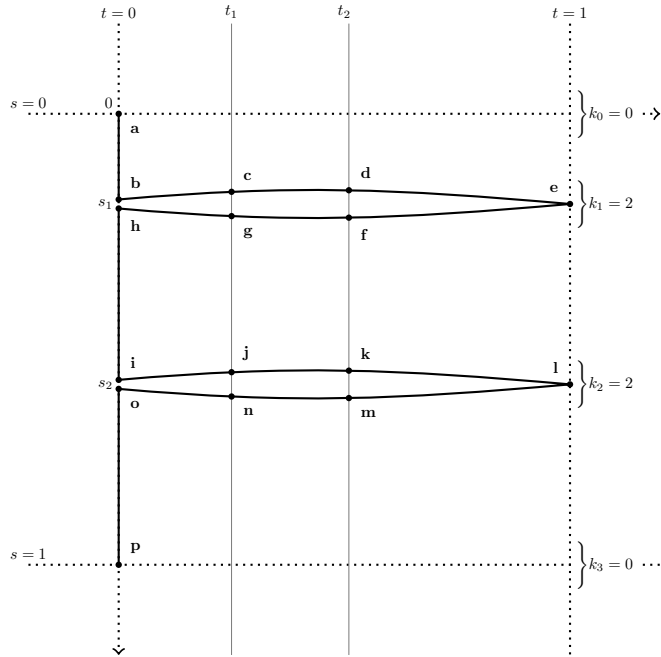
\begin{equation*}
\begin{tikzpicture}[thick,scale=0.6, every node/.style={transform shape}]
%help lines grid
%\draw[help lines] (0.1,0.1) grid (9.9,9.9);

% s axis
\draw[thick][dotted][->] (0,12) -- (0,-2);

% t axis
\draw[thick][dotted][->] (-2,10) -- (12,10);

%label axes
\node [above left] at (0,10) {0};
\node[above] at (0,12) {$t=0$};
\draw[thick][dotted] (10,12) -- (10,-2);
\node[above] at(10,12) {$t=1$};
\draw[thick][dotted] (-2,0) -- (12,0);
\node[above] at(-2,0) {$s=1$};
\node[above] at (-2,10) {$s=0$};

% s=0 zigzags
\draw[name path=zigzag] (0,10)
% s=s_1 zigzazags
-- (0,8.1) to [out=5, in=175] (10,8) to [out=185, in=-5](0,7.9)

%s=s_2 zigzags
-- (0,4.1) to [out=5, in=175] (10,4) to [out=185, in=-5](0,3.9)

%s=1 zigzags
-- (0,0);% to [out=5, in=175] (10,0) to [out=185, in=-5](0,-0.1);

%label si's
\node[left] at (0,8) {$s_1$};
\node[left] at (0,4) {$s_2$};

%label ki's
\node[right, fill=white] at (10,10) { $\Bigg\} k_0 = 0$};
\node[right, fill=white] at (10,8) { $\Bigg\} k_1 = 2$};
\node[right, fill=white] at (10,4) { $\Bigg\} k_2 = 2$};
\node[right, fill=white] at (10,0) { $\Bigg\} k_3 = 0$};

%t-lines

%t1
\draw[name path=t1, thin, gray] (2.5,-2) -- (2.5,12);
\node[above] at (2.5,12) {$t_1$};

%t2
\draw[name path=t2, thin, gray] (5.1,-2) -- (5.1,12);
\node[above] at (5,12) {$t_2$};

%t3
%\draw[name path=t3, thin, gray] (7.5,-2) -- (7.5,12);
%\node[above] at (7.5,12) {$t_3$};

\fill [name intersections={of=t1 and zigzag, by={a1, a2, a3, a4}}]
        \foreach \y in {1, 2, ...,4}
        {(a\y) circle (2pt)};
        \node[label=above right:$\bold{c}$] at (a1) {};
        \node[label=above right:$\bold{j}$] at (a3) {};
        \node[label=below right:$\bold{n}$] at (a4) {};
        \node[label=below right:$\bold{g}$] at (a2) {};

% \fill [name intersections={of=t3 and zigzag, by={a1, a2, a3, a4, a5, a6, a7, a8, a9, a10, a11, a12, a13, a14}}]
   %     \foreach \y in {1, 2, ..., 14}
      %  {(a\y) circle (2pt)};
        
\fill [name intersections={of=t2 and zigzag, by={a1, a2, a3, a4}}]
        \foreach \y in {1, 2, ...,4}
        {(a\y) circle (2pt)};
        \node[label=above right:$\bold{d}$] at (a1) {};
        \node[label=below right:$\bold{f}$] at (a2) {};
        \node[label=above right:$\bold{k}$] at (a3) {};
        \node[label=below right:$\bold{m}$] at (a4) {};

\fill (0,10) circle (2pt);
\node[label=below right:$\bold{a}$] at (0,10) {};
\fill (0,8.1) circle (2pt);
\node[label=above right:$\bold{b}$] at (0,8.1) {};
\fill (10,8) circle (2pt);
\node[label=above left:$\bold{e}$] at (10,8) {};
\fill (0,7.9) circle (2pt);
\node[label=below right:$\bold{h}$] at (0,7.9) {};
\fill (0,4.1) circle (2pt);
\node[label=above right:$\bold{i}$] at (0,4.1) {};
\fill (10,4) circle (2pt);
\node[label=above left:$\bold{l}$] at (10,4) {};
\fill (0,3.9) circle (2pt);
\node[label=below right:$\bold{o}$] at (0,3.9) {};
\fill (0,0) circle (2pt);
\node[label=above right:$\bold{p}$] at (0,0) {};

\end{tikzpicture}
\end{equation*}

When applying star to the two no-trivial zigzags above,$1 \otimes (\underline{x}_{(1,2)} \star \underline{x'}_{(1,2)}) \otimes 1$, the picture above simply becomes

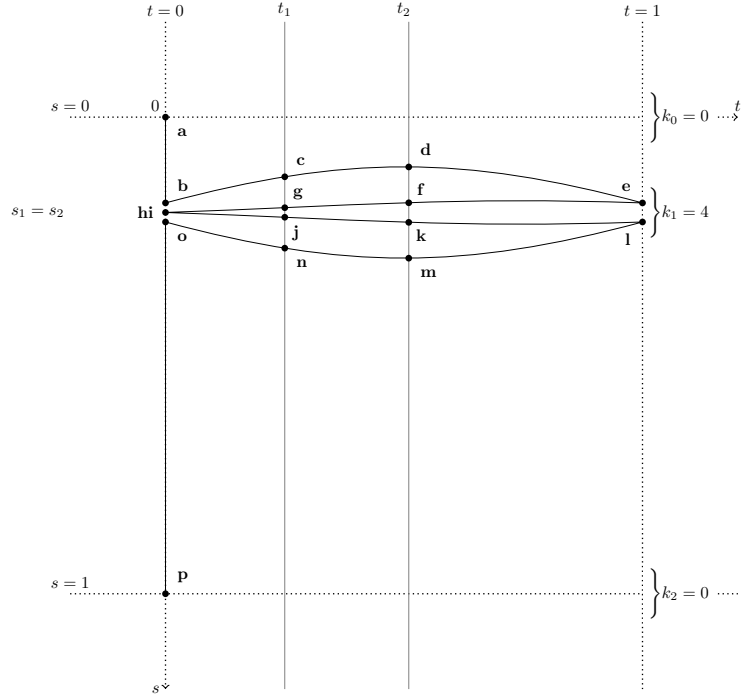
\begin{equation*}\resizebox{10cm}{!}{
\begin{tikzpicture}%[thick,scale=0.6, every node/.style={transform shape}]
%help lines grid
%\draw[help lines] (0.1,0.1) grid (9.9,9.9);

% s axis
\draw[thick][dotted][->] (0,12) -- (0,-2);
\node [left] at (0,-2) {$s$};
% t axis
\draw[thick][dotted][->] (-2,10) -- (12,10);
\node [above] at (12,10) {$t$};
%label axes
\node [above left] at (0,10) {0};
\node[above] at (0,12) {$t=0$};
\draw[thick][dotted] (10,12) -- (10,-2);
\node[above] at(10,12) {$t=1$};
\draw[thick][dotted] (-2,0) -- (12,0);
\node[above] at(-2,0) {$s=1$};
\node[above] at (-2,10) {$s=0$};

% s=0 zigzags
\draw[name path=zigzag] (0,10)
% s=s_1 zigzazags
-- (0,8.2) to [out=15, in=165] (10,8.2) to [out=178, in=2](0,8)
to [out=-2, in=182] (10,7.8) to [out=195, in=-15] (0,7.8)

%s=1 zigzags
-- (0,0);% to [out=5, in=175] (10,0) to [out=185, in=-5](0,-0.1);

%label si's
\node[left] at (-2,8) {$s_1= s_2$};

%label ki's
\node[right, fill=white] at (10,10) { $\Bigg\} k_0 = 0$};
\node[right, fill=white] at (10,8) { $\Bigg\} k_1 = 4$};
\node[right, fill=white] at (10,0) { $\Bigg\} k_2 = 0$};

%t-lines

%t1
\draw[name path=t1, thin, gray] (2.5,-2) -- (2.5,12);
\node[above] at (2.5,12) {$t_1$};

%t2
\draw[name path=t2, thin, gray] (5.1,-2) -- (5.1,12);
\node[above] at (5,12) {$t_2$};

%t3
%\draw[name path=t3, thin, gray] (7.5,-2) -- (7.5,12);
%\node[above] at (7.5,12) {$t_3$};

\fill [name intersections={of=t1 and zigzag, by={a1, a2, a3, a4}}]
        \foreach \y in {1, 2, ...,4}
        {(a\y) circle (2pt)};
        \node[label=above right:$\bold{c}$] at (a1) {};
        \node [inner sep=1pt, label={[label distance=0.01cm]-45:$\bold{j}$}] at (a3){};
        %\node[label=below right = .3cm:$\bold{j}$] at (a3) {};
        \node[label=below right:$\bold{n}$] at (a4) {};
        \node [inner sep=1pt, label={[label distance=0.01cm]45:$\bold{g}$}] at (a2){};

% \fill [name intersections={of=t3 and zigzag, by={a1, a2, a3, a4, a5, a6, a7, a8, a9, a10, a11, a12, a13, a14}}]
   %     \foreach \y in {1, 2, ..., 14}
      %  {(a\y) circle (2pt)};
        
\fill [name intersections={of=t2 and zigzag, by={a1, a2, a3, a4}}]
        \foreach \y in {1, 2, ...,4}
        {(a\y) circle (2pt)};
        \node[label=above right:$\bold{d}$] at (a1) {};
        \node [inner sep=1pt, label={[label distance=0.01cm]45:$\bold{f}$}] at (a2){};
        \node [inner sep=0pt, label={[label distance=0.01cm]-45:$\bold{k}$}] at (a3){};
        \node[label=below right:$\bold{m}$] at (a4) {};

\fill (0,10) circle (2pt);
\node[label=below right:$\bold{a}$] at (0,10) {};
\fill (0,8.2) circle (2pt);
\node[label=above right:$\bold{b}$] at (0,8.2) {};
\fill (10,8.2) circle (2pt);
\node[label=above left:$\bold{e}$] at (10,8.2) {};
\fill (0,8) circle (2pt);
\node[label=left:$\bold{hi}$] at (0,8) {};
\fill (10,7.8) circle (2pt);
\node[label=below left:$\bold{l}$] at (10,7.8) {};
\fill (0,7.8) circle (2pt);
\node[label=below right:$\bold{o}$] at (0,7.8) {};
\fill (0,0) circle (2pt);
\node[label=above right:$\bold{p}$] at (0,0) {};

\end{tikzpicture}
}\end{equation*}
The maps which describe this $\star$ operation include the diagonal
\begin{align*}
				&\Delta^1 \times \Delta^2 \times B(M) \to \Delta^2 \times \Delta^2 \times B(M)\\
\text{given by} \quad 	&(s, (t_1, t_2), \Sigma) 		     \mapsto  ((s, s), (t_1, t_2), \Sigma)
\end{align*}
and 
\begin{align*}
				M^{15} &\to M^{16}\\
\text{given by} \quad 	(x_1, \ldots , x_{15}) 		     &\mapsto  (x_1, \ldots, x_7, x_8, x_8, x_9, \ldots, x_{15})
\end{align*}
yielding the commutative diagram
\begin{equation*}
\xymatrix{\Delta^{1} \times \Delta^2 \times B(M) \ar[d] \ar[rrr]^{ev_{(0,2,2,0),2}} & & & M^{15} \ar[d] \\
\Delta^2 \times \Delta^2 \times B(M)  \ar[rrr]_{ev_{(0,4,0),2}} & & & M^{16}  }
\end{equation*}
So in our figures above, the multiplication of exactly $h$ and $i$ comes from the fact that the map $M^{15} \to M^{16}$ had exactly one diagonal map built into it.  Also recall that before the star operation, we had $k_1=2$ and $k_2=2$, and after the star operation we had $m=1$ and $k_1=4$.  In general, the evaluation maps in the commutative diagram need to keep track of the $k_j$'s in order for $CH^{ZZ}_{Rec}(\OmegaMat(M))$ to be relevant to, let alone be a model for, $\OmegaMat(BM)$. 
In general, when we want to investigate
 $$It\left( (\underline{x}^0_{(k_0, n)}) \otimes \ldots \otimes( (\underline{x}^r_{(k_r, n)})\star (\underline{x}^{r+1}_{(k_{r+1}, n)})) \otimes \ldots (\underline{x}^{m+1}_{(k_{m+1}, n)}) \right)$$
 we use the commutative diagram with $\phi_j = nk_j +k_j +1$,

\begin{equation*}
\xymatrix{\Delta^{m-1} \times \Delta^n \times B(M) \ar[d] \ar[rrrrrr]^{ev_{(k_0, \ldots, k_{r-1}, k_r + k_{r+1}, k_{r+2}, \ldots, k_{m+1}),n}} & &  & & & & M^{\left(\sum\limits_{j=0}^{m+1} \phi_j \right)-1}  \ar[d] \\
 \Delta^m \times \Delta^n \times B(M) \ar[rrrrrr]_{ev_{(k_0, \ldots , k_{m+1}),n}} & &  & & & & M^{\sum\limits_{j=0}^{m+1} \phi_j} \\  }
\end{equation*}
and see that the $\star$ operation corresponds to the wedge product during concatenation as shown in the figures above.  In the calculation of $d(\int\limits_{\Delta^n \times \Delta^m} ev^*(\underline{\omega}))$ we then obtain precisely the new term $\int\limits_{\Delta^n \times \partial\Delta^m} ev^*(\underline{\omega})$ (see Proposition \ref{prop zz chain map}).  Now we have:
\begin{prop}
The iterated integral $CH^{ZZ}_{Rec}(\OmegaMat(M)) \xrightarrow{It} \OmegaMat(B(M))$ is a chain map.
\end{prop}
Below we briefly mention a variation of the 2-d Rectangular Hochschild model and use the fact that since $CH^{ZZ}(\A)$ is a DGA, we can take $CH^I(CH^{ZZ}(\A))$.
\begin{definition}
$$CH^{ZZ}_{Sq}(\A): = \bigoplus\limits_{\substack{n \ge 0 \\ k_0, \ldots, k_{n+1} \ge 0 \\ k_i \text{ even}}} (\A \otimes (\A^{\otimes n} \otimes \A)^{\otimes {k_0}} \otimes \ldots \otimes \A \otimes (\A^{\otimes n} \otimes \A)^{\otimes {k_{n+1}}})[n] $$  is to be thought of as $(n+2)$-many rows of elements $(\underline{\omega}^i_{(k_i, n)}) \in CH^{ZZ}(\A)$, with differential 
\begin{align*}
&D((\underline{x}^0_{(k_0, n)})\otimes \ldots \otimes (\underline{x}^{n+1}_{(k_{n+1}, n)}))\\
:=& \sum\limits_{r=0}^{n} (-1)^{n+r} b_r \left( (\underline{x}^0_{(k_0, n)})\otimes \ldots \otimes(\underline{x}^r_{(k_r, n)})\star(\underline{x}^{r+1}_{(k_{r+1}, n)}) \otimes \ldots \otimes (\underline{x}^{n+1}_{(k_{n+1}, n)}) \right)\\
+ & \sum\limits_{r=1,p=1}^{n,n} (-1)^{n + \beta_{r,p}} (\underline{x^0}_{(k_0, n)})\otimes \ldots \otimes d_p((\underline{x}^r_{(k_r, n)})) \otimes \ldots \otimes (\underline{x}^{n+1}_{(k_{n+1}, n)})
\end{align*}
where $d_p$ is the differential coming from $\A$ applied to exactly the $p$-th slot of 
$(\underline{x}_{(k_i, n)})$ and $b_r$ is collapsing/multiplying the $r$-th and $(r+1)$-th slots of $(\underline{x}_{(k_i, n)})$.  Note that we use the fact that collapsing two columns commutes with the $\star$ product of two rows coming from elements in $CH^{ZZ}(\A)$. 
\end{definition}

\begin{prop}
If $\A$ is a commutative DGA, then using our map from Proposition \ref{prop zz to I}, there exists a commutative diagram of chain complexes:
\begin{equation*}
\xymatrix{
CH^I(CH^{ZZ}(\A))\ar[d] \ar[r] &CH^I(CH^I(\A))\ar[d] \\
CH^{ZZ}_{Rec}(\A)\ar[d]  \ar[r]        &CH^{I \times I}_{Rec}(\A) \ar[d] \\
CH^{ZZ}_{Sq}(\A) \ar[r] &CH^{I \times I}_{Sq}(\A)}
\end{equation*}
where the horizontal arrows come from collapsing zigzags and the vertical arrows come from adding in degeneracies (see for example \cite{GTZ}, Corollary 2.4.4, for the maps on the right). 
\end{prop}

\section{The zigzag Hochschild complex for a curved DGA}
\subsection{The one-dimensional case}
We have in mind integrating differential forms along our zigzags as before, but now we would like to apply parallel transport between the variables at which the differential forms are sitting.  When we take the De Rham differential after integration, we will have some extra terms show up.  For this reason we define the curved zigzag Hochschild complex, using the same underlying vector space as in the ``1-d'' case, $CH^{ZZ}(\A)$, but with an additional component added to its differential.  
\begin{definition}\label{def of curved CH^ZZ}
Let $(\A, \cdot, d)$ be a DGA and let $A \in \A$ be an element of degree 1.  Then denote by $\nabla: \A \to \A$, $\nabla(x):= dx + [A,x]$, with $\nabla_R(x):= dx + (-1)^{|x|} x\cdot A$ and $\nabla_L(x):= dx + A \cdot x$.  The curved zigzag Hochschild complex is defined as
 $$CH^{ZZ}(\A) = \bigoplus\limits_{n, k \ge 0, k \text{ even}}   (\A \otimes(\A^{\otimes n}\otimes \A)^{\otimes k})[n],$$ with differential $D: CH^{ZZ}(\A)  \to CH^{ZZ}(\A) $ given by $D( \underline{x}_{k,n}):=( \nabla + b + c) ( \underline{x}_{k,n})$.  We define these three components below: 
\begin{align*}
\nabla( \underline{x}_{k,n}):= &(-1)^n \nabla_R(x^{\mathcal{L}}) \otimes  \ldots \otimes  x^{\mathcal{L}}_k)\\
+ & \sum\limits_{p=1}^n \sum\limits_{i=1}^k (-1)^{n+\beta_{(i,p)}} x^{\mathcal{L}} \otimes  \ldots \otimes \nabla(x_{(i, p)}) \otimes \ldots  x_k^{\mathcal{L}}\\
+& \sum\limits_{i=1}^{k-1} (-1)^{n+\beta_i^{\dagger_i}}  x^{\mathcal{L}}  \otimes \ldots \otimes \nabla(x^{\dagger_i}_i)\otimes  \ldots \otimes x_k^{\mathcal{L}}\\
+&(-1)^{n+ \beta_k^{\mathcal{L}}} x^{\mathcal{L}} \otimes  \ldots \otimes  \nabla_L(x^{\mathcal{L}}_k)\\
\intertext{$b$ is defined in exactly the same way as in Definition \ref{def CH^ZZ}, and }
c( \underline{x}_{k,n}):=  
& \sum\limits_{\sigma \in S_{n,1}}
c_{\sigma}^R(\underline{x}_{k,n}) \\
= & \sum\limits_{\sigma \in S_{n,1}} \sum\limits_{i=1}^k 
(-1)^{n+ \sigma(n+1)} x^{\mathcal{L}} \otimes \ldots \\
& \otimes \ldots \otimes \left( \sigma ((x_{(i,1)} \otimes \ldots \otimes x_{(i,n)}), R) \otimes x^{\dagger_i}_i  \right)\\
& \otimes \ldots \otimes \left( \sigma^{\Sh_i} ((x_{(j,1)} \otimes \ldots \otimes x_{(j,n)}), 1) \otimes x^{\dagger_j}_j \right) \ldots 
\end{align*}
The term $c_{\sigma}^R$ shuffles in exactly one new column in $\underline{x}_{k,n}$ which consists of all 1's except for exactly one $R$, and whose placement is determined by $\sigma(n+1)$ (see Remark \ref{shuffle remarks} for the definition of $\sigma$ and $\sigma^{\Sh_i}$):  
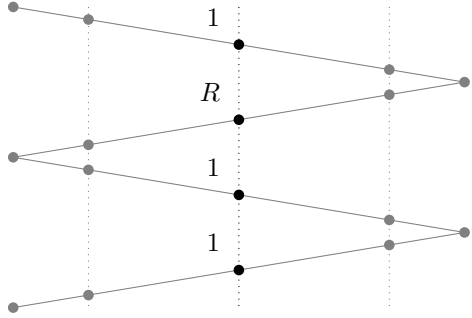
\begin{equation*}
\begin{tikzpicture}
\draw[name path=zigzagtop, gray] (0,10) to (6,9) to (0,8) to (6,7) to (0,6); 
\fill[gray] (0,10) circle (2pt) ;
\fill[gray] (0,8) circle (2pt) ;
\fill[gray] (0,6) circle (2pt) ;
\fill[gray] (6,9) circle (2pt) ;
\fill[gray] (6,7) circle (2pt) ;
\draw[gray, dotted,name path=t1] (1, 10) to (1,6);
\fill [gray, name intersections={of=t1 and zigzagtop, by={a1,b1, c1, d1}}]
        (a1) circle (2pt)
        (b1) circle (2pt)
        (c1) circle (2pt)
        (d1) circle (2pt);

\draw[dotted, name path=t2] (3, 10) to (3,6);
\fill [name intersections={of=t2 and zigzagtop, by={a2,b2, c2, d2}}]
        (a2) circle (2pt)
        (b2) circle (2pt)
        (c2) circle (2pt)
        (d2) circle (2pt);
\node[label=above left:{$1$}] at (a2) {};
\node[label=above left:$R$] at (b2) {};
\node[label=above left:$1$] at (c2) {};
\node[label=above left:$1$] at (d2) {};

%\draw[opacity=0, name path=t2] (2, 10) to (2,0);

%\draw[opacity=0, name path=t3] (8, 10) to (8,0);

\draw[dotted, gray, name path=t3] (5, 10) to (5,6);
\fill [gray, name intersections={of=t3 and zigzagtop, by={a3,b3, c3, d3}}]
        (a3) circle (2pt)
        (b3) circle (2pt)
        (c3) circle (2pt)
        (d3) circle (2pt);

\end{tikzpicture}
\end{equation*}
Moreover, we have the same shuffle product, $\odot$, for the curved zigzag Hochschild complex as we do in the zigzag Hochschild complex, defined in Definition \ref{zz shuffle}.
\end{definition} 
\begin{prop}
For a unital DGA, $(A, \cdot, d)$ with $A \in \A^1$, the curved differential $D$ on $CH^{ZZ}(\A)$ satisfies $D^2=0$.  Furthermore, $D$ is a (graded) derivation of $\odot$.  
\begin{proof}
To see that $D^2 =0$, we observe $\nabla^2 + cb + bc =0$, $b^2=c^2 =0$, $\nabla b + b \nabla=0$, and $\nabla c + c \nabla = 0$.  The proof is straightforward but the reader should note that the terms $\nabla_L$ and $\nabla_R$ correspond to $R$ being multiplied to the left or the right of the endpoints, respectively, when considering $b \circ c$.  The proof that $D$ is a derivation of $\odot$ is comparable to that of Proposition \ref{propzzDderivation} if one replaces all instances of $d$ with $\nabla$.
\end{proof}
\end{prop}
The setup we actually have in mind for this curved algebra case is when $\A = \OmegaMat(M) := \Omega(M, Mat)$ is the DGA of matrix-valued forms on $M$, or more generally, $\OmegaMat(M) = \Omega(M, \mathfrak{g})$ where $\mathfrak{g}$ is a Lie-Algebra whose bracket comes from an underlying product.  Now we fix a 1-form $A \in \OmegaMat^1(M)$. 
We wish to define a map $It^A: CH^{ZZ}(\OmegaMat(M)) \to \OmegaMat(PM)$.  We illustrate the map on a simple monomial in $CH^{ZZ}(\OmegaMat(M))$.  So consider the following element of $\OmegaMat(M)^{\otimes 7}[2] \subset CH^{ZZ}(\OmegaMat(M))$, where the evaluation map normally used in the non-curved case $\Delta^2 \times PM  \to M^{7}$ should be evident from the picture, 
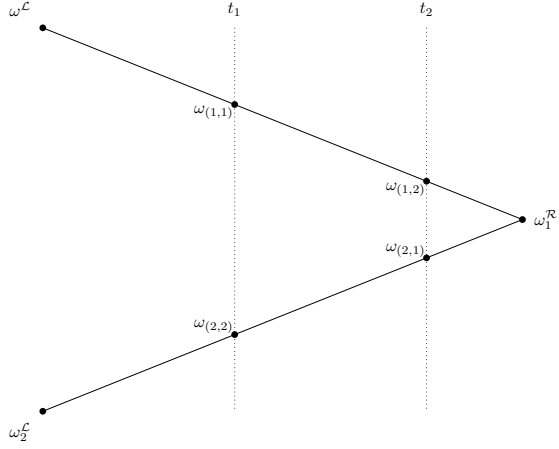
\begin{equation*}\resizebox{!}{6cm}{
\begin{tikzpicture}
\draw[name path=zigzag] (0,12) -- (10,8) -- (0,4);
\draw[dotted, opacity=50,name path=t1] (4, 12) to (4,4);
   \fill [name intersections={of=t1 and zigzag, by={a,b}}]
        \foreach \x in {a,b}
        {(\x) circle (2pt)};
        \node[below left = -0.1cm] at (a) {$\omega_{(1,1)}$};
        \node[above left = -0.1cm] at (b) {$\omega_{(2,2)}$};
    \node[label=above:{$t_1$}] at (4,12) {};

\draw[dotted, opacity=50,name path=t2] (8, 12) to (8,4);
   \fill [name intersections={of=t2 and zigzag, by={a,b}}]
        \foreach \x in {a,b}
        {(\x) circle (2pt)};
        \node[below left = -0.1cm] at (a) {$\omega_{(1,2)}$};
        \node[above left = -0.1cm] at (b) {$\omega_{(2,1)}$};
     \node[label=above:{$t_2$}] at (8,12) {};

\fill (0,12) circle (2pt);
\node[label=above left:{$\omega^{\mathcal{L}}$}] at (0,12) {};
\fill (10,8) circle (2pt);
\node[label=right:{$\omega^{\mathcal{R}}_1$}] at (10,8) {};
\fill (0,4) circle (2pt);
\node[label=below left:{$\omega^{\mathcal{L}}_2$}] at (0,4) {};
\end{tikzpicture}}
\end{equation*}
 which represents an element 
$$\omega^{\mathcal{L}} \otimes \omega_{(1,1)} \otimes \omega_{(1,2)} \otimes \omega^{\mathcal{R}}_1 \otimes \omega_{(2,1)} \otimes \omega_{(2,2)} \otimes \omega_2^{\mathcal{L}} \in \OmegaMat(M)^{\otimes 7}[2]$$
Next, we apply a map of degree zero\footnote{While we are inserting arbitrarily many 1-forms $A$, they will all be shifted down by one-degree, since they each will be integrated along some interval at some $\tau$.  This is consistent with the rest of the shifts in the paper where we always shift a monomial by the dimension of the fiber we integrate over.} , $$Ins^A_{(q_1, q_2, \ldots, q_6)}: \OmegaMat(M)^{\otimes 7}[2] \to \OmegaMat(M)^{\otimes {7 + \sum q_i}}[2+ \sum q_i]$$ which will insert some $A$'s as prescribed by the $q_i \in \mathbb{Z}_{\ge 0}$.  So applying $Ins^A_{(2, 3, 1, 2, 0, 5)}:  \OmegaMat(M)^{\otimes 7}[2] \to \OmegaMat(M)^{20}[15]$ to our element in $CH^{ZZ}(\OmegaMat(M))$ above yields
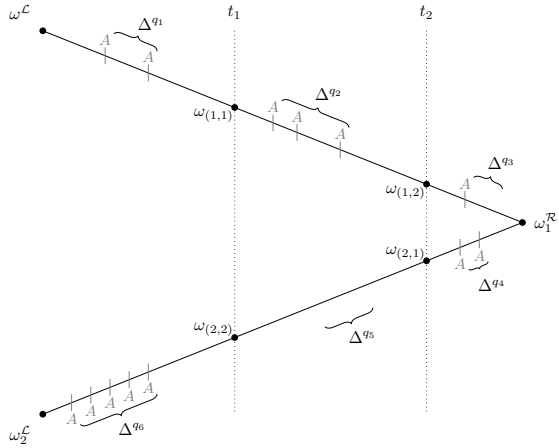
\begin{equation*}\resizebox{!}{6cm}{
\begin{tikzpicture}
\draw[name path=zigzag] (0,12) --  coordinate[pos=.13] (A1)  coordinate[pos=.18] (D1) coordinate[pos=.22] (A2) coordinate[pos=.48] (A3) coordinate[pos=.53] (A4)  coordinate[pos=.55] (D2) coordinate[pos=.62] (A5) coordinate[pos=.88] (A6)coordinate[pos=.91] (D3) coordinate[pos=.94] (B6) (10,8) -- coordinate[pos=.09] (A7) coordinate[pos=.11] (D4)   coordinate[pos=.13] (A8)  coordinate[pos=.33] (B8) coordinate[pos=.38] (D5)  coordinate[pos=.43] (B9)   coordinate[pos=.78] (A9) coordinate[pos=.82] (A10)  coordinate[pos=.86] (A11) coordinate[pos=.90] (A12) coordinate[pos=.94] (A13)(0,4);
\draw[dotted, opacity=50,name path=t1] (4, 12) to (4,4);
   \fill [name intersections={of=t1 and zigzag, by={a,b}}]
        \foreach \x in {a,b}
        {(\x) circle (2pt)};
        \node[below left = -0.1cm] at (a) {$\omega_{(1,1)}$};
        \node[above left = -0.1cm] at (b) {$\omega_{(2,2)}$};
    \node[label=above:{$t_1$}] at (4,12) {};

\draw[dotted, opacity=50,name path=t2] (8, 12) to (8,4);
   \fill [name intersections={of=t2 and zigzag, by={a,b}}]
        \foreach \x in {a,b}
        {(\x) circle (2pt)};
        \node[below left = -0.1cm] at (a) {$\omega_{(1,2)}$};
        \node[above left = -0.1cm] at (b) {$\omega_{(2,1)}$};
      \node[label=above:{$t_2$}] at (8,12) {};

\fill (0,12) circle (2pt);
\node[label=above left:{$\omega^{\mathcal{L}}$}] at (0,12) {};
\fill (10,8) circle (2pt);
\node[label=right:{$\omega^{\mathcal{R}}_1$}] at (10,8) {};
\fill (0,4) circle (2pt);
\node[label=below left:{$\omega^{\mathcal{L}}_2$}] at (0,4) {};

\foreach \x in {1,2,...,6}
{  \draw[gray] ($(A\x)+(0,5pt)$) -- ($(A\x)-(0,5pt)$);

\node[above= 0.1cm, gray] at (A\x) {\small$A$};}

\foreach \x in {7,8,...,13}
{  \draw[gray] ($(A\x)+(0,5pt)$) -- ($(A\x)-(0,5pt)$);

\node[below= 0.1cm, gray] at (A\x) {\small$A$};}

\draw [decorate,decoration={brace,amplitude=4pt,raise=14pt}]
(A1) -- (A2);
\node [label={[label distance=0.5cm]85:$\Delta^{q_1}$}]  at (D1) {};

% \node[below  = 0.4cm] at (D1) {$\Delta^{q_1}$};

\draw [decorate,decoration={brace,amplitude=4pt,raise=14pt}]
(A3) -- (A5);

\node [label={[label distance=0.5cm]85:$\Delta^{q_2}$}]  at (D2) {};
 %\node[below  = 0.4cm] at (D2) {$\Delta^{q_2}$};
\draw [decorate,decoration={brace,amplitude=4pt,raise=14pt}]
(A6) -- (B6);

\node [label={[label distance=0.5cm]85:$\Delta^{q_3}$}]  at (D3) {};

\draw [decorate,decoration={brace,amplitude=2pt,raise=14pt}]
(A7) -- (A8);

\node [label={[label distance=0.5cm]-85:$\Delta^{q_4}$}]  at (D4) {};

\draw [decorate,decoration={brace,amplitude=4pt,raise=14pt}]
(B8) -- (B9);

\node [label={[label distance=0.5cm]-85:$\Delta^{q_5}$}]  at (D5) {};

\draw [decorate,decoration={brace,amplitude=4pt,raise=14pt}]
(A9) -- (A13);

\node [label={[label distance=0.5cm]-85:$\Delta^{q_6}$}]  at (A11) {};

\end{tikzpicture}}\end{equation*}
resulting in an element $$\omega^{\mathcal{L}} \otimes A^{\otimes 2} \otimes \omega_{(1,1)} \otimes  A^{\otimes 3} \otimes \omega_{(1,2)} \otimes  A \otimes \omega^{\mathcal{R}}_1 \otimes  A^{\otimes 2} \otimes  \omega_{(2,1)} \otimes \omega_{(2,2)} \otimes  A^{\otimes 5} \otimes \omega_2^{\mathcal{L}} \in \OmegaMat(M)^{\otimes {20}}[15]$$
where, for example $\Delta^{q_5} = {*}$ and $$\Delta^{q_6} = \{( \tau_1, \ldots, \tau_5) \in \mathbb{R}^5 | t_1 \ge \tau_1 \ge \tau_2 \ge \ldots \ge \tau_5 \ge 0 \}.$$
We consider the bounded convex polytope:
\begin{align*}
E = \{ (&\tau^1_1, \tau^1_2, \tau^2_1, \tau^2_2 , \tau^2_3,\tau^3_1,\tau^4_1, \tau^4_2, \tau^6_1, \ldots \tau^6_5, t_1, t_2) | \\
&0 \le t_1 \le t_2 \le 1,\\
&0 \le \tau^1_1 \le  \tau^1_2 \le t_1 \le \tau^2_1 \le  \tau^2_2 \le \tau^2_3 \le  t_2 \le  \tau^3_1 \le 1,\\
&1 \ge \tau^4_1 \ge  \tau^4_2 \ge t_2 \ge t_1 \ge  \tau^6_1 \ge \ldots  \ge \tau^6_5 \ge 0\}.
\end{align*}
Once again, we use the picture as our guide for the evaluation map, and define a summand of our iterated integral using the following diagram (ignoring degree-shifts): 
\begin{equation*}
\xymatrix{\OmegaMat(PM \times E)    \ar[d]^{\int_E } &\ar[l]_{ev^*}\OmegaMat(M^{20}) &\ar[l]  \OmegaMat(M)^{\otimes {20}}  & \ar[llld]^{It^A_{\underline{q}}}  \ar[l]_{ Ins^A_{\underline{q}}}  \OmegaMat(M)^{\otimes {7}} \\
\OmegaMat(PM) &  & &}
\end{equation*}
where $\underline{q} = (2, 3, 1, 2, 0, 5)$ in this case and then we can finally define the iterated integral map $\OmegaMat(M)^{\otimes {7}}[2] \xrightarrow{It^A} \OmegaMat(PM)$ by $It^A (\underline{\omega}) := \sum\limits_{q_1, \ldots, q_6 \ge 0} (-1)^{\sum |\omega_i| Q_i}It^A_{(q_1, \ldots, q_6)}(\underline{\omega})$, where here by $\omega_i$ we mean the $i$-th tensor-factor in $\underline{\omega} = \omega_0\otimes \dots \otimes \omega_{(n+1)k +1}$ and we define $Q_i := q_{i+1} + \ldots + q_{(n+1)k }$.  In the same way we can now define the iterated integral in general.

\begin{definition}   Consider some monomial $\underline{\omega}_{k,n} \in \OmegaMat(M)^{\otimes \phi}[n] \subset CH^{ZZ}(\OmegaMat(M))$ where $\phi = nk + k +1$. For each choice $q_i \ge 0$, for $i=1, \ldots, (n+1)\cdot k$ we have the corresponding diagram
\begin{equation*}
\xymatrix{\OmegaMat(PM \times E)    \ar[d]^{\int_{E}} &\ar[l]_{ev^*}\OmegaMat(M^{\phi^A}) &\ar[l]  \OmegaMat(M)^{\otimes {\phi^A}}  & \ar[llld]^{It^A_{\underline{q}}}  \ar[l]_{ Ins^A_{\underline{q}}}  \OmegaMat(M)^{\otimes {\phi}} \\
\OmegaMat(PM) &  & &}
\end{equation*}
where $\phi^A = \phi + \sum q_i$ and $\underline{q} = (q_1 , \ldots, q_{\phi -1})$.  Then $\OmegaMat(M)^{\otimes {\phi}}[n] \xrightarrow{It^A} \OmegaMat(PM)$ is defined by $It^A := \sum\limits_{\underline{q}} It^A_{\underline{q}}$.
\end{definition}
A point should be made about this infinite sum.  In the case when we were considering $It :\OmegaMat(M)^{\otimes 7} \to \OmegaMat(PM)$ we were considering a cartesian product of 6 different infinite-sums of choices: $q_1, \ldots, q_6$.  However, for a fixed path $\gamma$ and $(t_1, t_2)$, we want to simply compute the parallel transport between two points on $\gamma$ using a connection 1-form, $A$.  In particular, since the image of $\gamma$ is compact in $M$, $|A|$ is bounded on $\gamma$ by some element $\rho \in \mathbb{R}$ and so we consider
\begin{align*}
&|1 + \int_{\Delta^1} A + \int_{\Delta^2} AA + \ldots|\\
\le &|1| + |\int_{\Delta^1} A| + |\int_{\Delta^2} AA| + \ldots\\
\le &1 + {\rho} + \frac{\rho^2}{2!} + \frac{\rho^3}{3!} + \ldots \\
= &\exp(\rho) \in \mathbb{R}.
\end{align*}
This shows that the infinite sum in $It^A$ indeed converges.  
\begin{prop}\label{It^A chain and algebra}
The map $It^A: CH^{ZZ}(\OmegaMat(M)) \to \OmegaMat(PM)$ as defined above is a chain map and an algebra map.
\begin{proof}
Contrasting with Proposition \ref{zz It algebra map}, we now have inserted $A$'s.  However, the insertion of $A$'s produces parallel transport functions which will not interact significantly with the wedge product.  Hence, up to sign, the proof of $It^A$ being an algebra map similar to Proposition \ref{zz It algebra map}, since we have a diagram using the evaluation map $PM \times E \xrightarrow{ev} M^{\phi^A}$ where no factors of $M$ have to be switched.  The computations involving keeping track of the choices of $\underline{q}_i$, although not trivial, are more straightforward than those for the proof that we have a chain map and so we leave those details to the reader.

We now show that $It^A$ is a chain map by applying $d_{DR}$ and focusing on what occurs in two scenarios: at some $\omega_{(-,-)}$ and at some inserted $A$, respectively.  We have:
\begin{align}
d\circ It^A(\underline{\omega}) &= \sum\limits_{\underline{q}} d \circ It^A_{\underline{q}}(\underline{\omega}) \nonumber \\
&=\sum\limits_{\underline{q}}  (-1)^{dimE} \int\limits_{E} d(ev^*(Ins^A_{\underline{q}}(\underline{\omega})))+ 
\sum\limits_{\underline{q}}(-1)^{dimE-1}  \int\limits_{\partial E} ev^*(Ins^A_{\underline{q}}(\underline{\omega})).\label{disdplusboundary}
\end{align}
First recall that for our bounded polytope $$E= \{ \ldots \tau^1_a \le \ldots \le \tau^{q_a}_a \le t_i \le \tau^1_b \le \dots \le \tau^{q_b}_b \le t_{i+1} \ldots\},$$taking the place of $\Delta^n$, is a subspace of $( \prod\limits_{q_i}  \Delta^{q_i}) \times \Delta^n $ given by elements that satisfy conditions such as $t_i \le \tau_b^1$ and $\tau_b^{q_b} \le t_{i+1}$, etc.  Using this as a guide, we note that $\partial E$ has components of the form 
$$\restr{E}{t_i = t_{i+1}} \quad \restr{E}{\tau_a^j = \tau_a^{j+1}} \quad \restr{E}{t_i = \tau_b^1} \quad \restr{E}{\tau_a^{q_a} = t_{i}}.$$  Each component of the boundary comes with a well-defined induced orientation coming from Stokes' Theorem (using outward pointing normal-vectors of each component).  We have boundary maps $\partial_{(-,-)}$ and $\partial_{-}$ which take {\it adjacent} coordinates and identify them.  For example, we consider the adjacent coordinates $$\tau^1_a \le \ldots \le \tau^{q_a}_a \le t_i \le \tau^1_b \le \dots \le \tau^{q_b}_b \le t_{i+1} \ldots$$ in $E$.  Then we have maps of the form
\begin{align*}
\partial_{(\tau_a^{q_a}, t_i )}: \restr{E}{\tau^{q_a}_a =  t_i} &\hookrightarrow E \\
\partial_{(\tau_a^{j}, \tau_a^{j+1})}: \restr{E}{\tau^{j}_a=  \tau^{j+1}_a} &\hookrightarrow E\\
\partial_{(t_i, \tau_b^1)}: \restr{E}{t_i = \tau^1_b } &\hookrightarrow E \\
\partial_{i}: \restr{E}{t_i = \tau^1_b = \dots = \tau^{q_b}_b = t_{i+1}} &\hookrightarrow E 
\end{align*} Note that for any zig-zig diagram having $n$ columns and $k$ zigzags, consisting of information $(n, k, \underline{q})$, we can associate its corresponding information $\partial_{(-,-)}(n, k, \underline{q})$ and 
$\partial_{i}(n, k, \underline{q})$, coming from the zigzag diagrams of the boundaries $\partial_{(-,-)}$ and $\partial_i$, respectively.  We have corresponding maps $b_{(-,-)}$ and $b_i$ which multiply the appropriate differential-forms in our monomial $Ins^A_{\underline{q}} (\underline{\omega})$.  We can rewrite the term on the right in equation \eqref{disdplusboundary} above as 
\begin{align}
&\sum\limits_{\underline{q}} (-1)^{\epsilon_{\underline{q}} + \epsilon_{\partial E}} \int\limits_{\partial E} ev^*(Ins^A_{\underline{q}}(\underline{\omega}))\\  &= \sum\limits_{\underline{q}} (-1)^{\epsilon_{\underline{q}} + \epsilon_{\partial E}}  \sum\limits_{i}  (-1)^{\epsilon_{\partial_i E}} \int_{\partial_{i} E}(\partial_{i} \times id)^*(ev_{(n,k, \underline{q})}^*(Ins^A_{\underline{q}}(\underline{\omega}))) \nonumber \\
&+\sum\limits_{\underline{q}} (-1)^{\epsilon_{\underline{q}} + \epsilon_{\partial E}} \sum\limits_{\tiny{(-,-)}}  (-1)^{\epsilon_{\partial_{\tiny{(-,-)}}E}} \int_{\partial_{\tiny{(-,-)}} E}( \partial_{\tiny{(-,-)}} \times id)^*(ev_{(n,k, \underline{q})}^*(Ins^A_{\underline{q}}(\underline{\omega}))) \nonumber \\
&= \sum\limits_{\underline{q}} (-1)^{\epsilon_{\underline{q}} + \epsilon_{\partial E}}\sum\limits_{i} (-1)^{\epsilon_{\partial_i E}} \int_{\partial_{i} E} (ev_{\partial_{i}(n,k, \underline{q})}^*(b_{i}(Ins^A_{\partial_i(\underline{q})}(\underline{\omega}))) \label{boundarytsandtaus}\\
&+ \sum\limits_{\underline{q}} (-1)^{\epsilon_{\underline{q}} + \epsilon_{\partial E}}  \sum\limits_{\tiny{(-,-)}}  (-1)^{\epsilon_{\partial_{\tiny{(-,-)}}E}}  \int_{\partial_{\tiny{(-,-)}} E} (ev_{\partial_{\tiny{(-,-)}}(n,k, \underline{q})}^*(b_{\tiny{(-,-)}}(Ins^A_{\underline{q}}(\underline{\omega}))) \nonumber
\end{align} 
where $\epsilon_{\partial E}: =n + 1 + \sum q_r$, $\epsilon_{\underline{q}}:= \sum |\omega_r| Q_r$, $\epsilon_{\partial_i E}: = i +1 + \sum q_i$, and $\epsilon_{\partial_{(-,-)}E}$ is defined as well by the orientation of $E$.  We also use the fact that the pullback along $\partial_i$  amounts to pulling back along a diagonal and so we wedge the adjacent\footnote{or $\omega_{(i, n-p)}$ and $\omega_{(i, n-p+1)}$ depending upon whether or not the $i$-th level is a ``zig'' or a ``zag''}  forms $\omega_{(i, p)}$ and $\omega_{(i, p+1)}$.  Note, however, that we have essentially dropped all of the inserted $A$'s between the two adjacent forms.  This is only because if we were to follow through with the iterated integral, we would be integrating along a 0-dimensional subspace (i.e. a point) and so the integral over that point of $1 + A + A\wedge A + \ldots$ would equal 1.  This is the same as saying that parallel transport along the constant path must equal the identity.  When we use the pullback along the $\partial_{(-,-)}$ we simply wedge the adjacent forms ($A\wedge \omega$, $A\wedge A$, or $\omega \wedge A$) and no further identifications are used.  Next we rewrite the first term on the right side of  \eqref{disdplusboundary} 
\begin{align}
\sum\limits_{\underline{q}} (-1)^{\epsilon_{\underline{q}} + \epsilon_E} \int\limits_{E} d(ev^*(Ins^A_{\underline{q}}(\underline{\omega}))) &= \sum\limits_{\underline{q}} (-1)^{\epsilon_{\underline{q}} + \epsilon_E}  \int\limits_{E} ev^*(d(Ins^A_{\underline{q}}(\underline{\omega}))) \nonumber \\
&= \sum\limits_{\underline{q}} (-1)^{\epsilon_{\underline{q}} + \epsilon_E}  \sum\limits_{l=1}^{k} \sum\limits_{i=1}^{n} \int\limits_{E} ev^*(d_{(i,l)} (Ins^A_{\underline{q}}(\underline{\omega}))) \label{tdsandtauds}\\
&+ \sum\limits_{\underline{q}} (-1)^{\epsilon_{\underline{q}} + \epsilon_E}  \sum\limits_{\tau_a^i} \int\limits_{E} ev^*(d_{\tau_a^i} (Ins^A_{\underline{q}}(\underline{\omega}))) \nonumber
\end{align}
where $d_{(i,l)}$ (in a slight abuse of notation) applies the DeRham differential to a form $\omega^{\mathcal{L}}$, $\omega_{(i,l)}$, or $\omega_l^{\dagger_l}$ and $d_{\tau_a^i}$ applies $d$ to an inserted $A$.  We now show where our term $\nabla(\omega_{(i,k)})$ comes from.  From equation \eqref{boundarytsandtaus} we get our suggestively-labeled ``$A\wedge\omega$'' and  ``$\omega \wedge A$'' terms, without sign\footnote{For the rest of this proof, we will proceed without sign as it is mostly straightforward to check that our sign conventions work out, but the details would  unnecessarily obfuscate the ideas.}, 
\begin{align*}
&A\wedge \omega  &\int_{\partial_{(\tau_a^{q_a}, t_i)}E} ev^*_{\partial_{(\tau_a^{q_a}, t_i)}(n,k,\underline{q})}(b_{(\tau_a^{q_a}, t_i)} (Ins^A_{\underline{q}} (\underline{\omega}))) \\
&\omega \wedge A  &\int_{\partial_{(t_i, \tau_b^{1})}E} ev^*_{\partial_{(t_i , \tau_b^{1})}(n,k,\underline{q})}(b_{(t_i, \tau_b^{1})} (Ins^A_{\underline{q}} (\underline{\omega}))) \\
\intertext{and from equation \eqref{tdsandtauds} we get another suggestively-labeled term} \\
 &d\omega   &\int\limits_{E} ev^*(d_{(i,k)} (Ins^A_{\underline{q}}(\underline{\omega}))) \end{align*}
Notice that for these three terms, all of the changes to the fiber $E$ involved removing a single $\tau$ or none at all.  Thus, when we sum over all choices of $\underline{q}$ we will recover any removed ${\tau}$ slots.  Now we can write 
\begin{align*}
&A \wedge \omega & \sum\limits_{\underline{q}} \pm \int_{\partial_{(\tau_a^{q_a}, t_i)}E} ev^*_{\partial_{(\tau_a^{q_a}, t_i)}(n,k,\underline{q})}(b_{(\tau_a^{q_a}, t_i)} (Ins^A_{\underline{q}} (\underline{\omega}))) \\
&\omega \wedge A &+ \sum\limits_{\underline{q}}  \pm \int_{\partial_{(t_i, \tau_b^{1})}E} ev^*_{\partial_{(t_i , \tau_b^{1})}(n,k,\underline{q})}(b_{(t_i, \tau_b^{1})} (Ins^A_{\underline{q}} (\underline{\omega}))) \\
&d \omega&+ \sum\limits_{\underline{q}} \pm \int\limits_{E} ev^*(d_{(i,l)} (Ins^A_{\underline{q}}(\underline{\omega}))) \\
&d \omega + [A, \omega]& = It^A(\nabla_{(i,l)}(\underline{\omega})),
\end{align*}
where $\nabla_{(i,l)}$ applies $\nabla$ to the corresponding component of $\underline{\omega}$.  
Next we show where shuffling in the $R$'s comes from.  From equations \eqref{boundarytsandtaus} and \eqref{tdsandtauds} we obtain
\begin{align*}
&A \wedge A \quad  &  \int_{\partial_{(\tau_a^j,\tau_a^{j+1})} E} (ev_{\partial_{(\tau_a^j,\tau_a^{j+1})}(n,k, \underline{q})}^*(b_{(\tau_a^j,\tau_a^{j+1})}(Ins^A_{\underline{q}}(\underline{\omega}))))\\
& dA \quad & \int\limits_{E} ev^*(d_{\tau_a^j} (Ins^A_{\underline{q}}(\underline{\omega}))) 
\end{align*}
and so by a similar argument as above we can write
\begin{align*}
 &A \wedge A & \sum\limits_{\underline{q}}  \sum\limits_{j=1, \ldots, q_a} \pm \int_{\partial_{(\tau_a^j,\tau_a^{j+1})} E} (ev_{\partial_{(\tau_a^j,\tau_a^{j+1})}(n,k, \underline{q})}^*(b_{(\tau_a^j,\tau_a^{j+1})}(Ins^A_{\underline{q}}(\underline{\omega}))))\\
&d A &+  \sum\limits_{\underline{q}}  \sum\limits_{j=1, \ldots, q_a} \pm  \int\limits_{E} ev^*(d_{\tau_a^j} (Ins^A_{\underline{q}}(\underline{\omega})))\\
& R = dA + A \wedge A  &=It^A(c_{(i,l)}(\underline{\omega}))
\end{align*}
Note that the sum of these kinds of terms gives $It^A(c (\underline{\omega}))$ since these terms insert exactly one $R = dA + A\wedge A$ in all positions and $1$'s in the corresponding positions $\tau_a^i$ of the other zig zags.  Finally we note that the first term on the right hand side of \eqref{boundarytsandtaus} gives
\begin{align*}
&\omega \wedge \omega'  & \int_{\partial_{i} E} (ev_{\partial_{i}(n,k, \underline{q})}^*(b_{i}(Ins^A_{\partial_i(\underline{q})}(\underline{\omega}))) \\
\intertext{where $b_i$ collapses the $i$\textsuperscript{ th} and $i+1$\textsuperscript{ th} columns in $\underline{\omega}$.  This provides us with the $b_i$ part of our differential again after summing over all configurations $\underline{q}$ to give:}
&\omega \wedge \omega' & \sum\limits_{\underline{q}} \pm \int_{\partial_{i} E} (ev_{\partial_{i}(n,k, \underline{q})}^*(b_{i}(Ins^A_{\partial_i(\underline{q})}(\underline{\omega}))) = It^A(b_i(\underline{\omega}))
\end{align*}
Similar calculations can be made at the endpoints.  In particular, if we focus our attention to $\omega^{\mathcal{L}}$, when we apply $d$ we obtain $d \omega^{\mathcal{L}}$, considering the boundary $\partial_{(0,\tau^1_1)}$ yields a term of the form  $\omega^{\mathcal{L}} \wedge A$, and $\partial_{0}$ yields $\omega^{\mathcal{L}} \wedge \omega_{(1,1)}$.  Similar terms arise when we focus our attention to $\omega_k^{\mathcal{L}}$.  Focusing our attention to $\omega_{i}^{\mathcal{R}}$, when we apply $d$ we obtain $d \omega_{i}^{\mathcal{R}}$, considering the boundary $\partial_{(\tau_a^{q_a},1)}$ we obtain $A \wedge \omega_{i}^{\mathcal{R}}$, applying the boundary $\partial_{(0, \tau_b^{1})}$ yields $\omega_{i}^{\mathcal{R}} \wedge A$, and the boundary $\partial_{(t_n ,1)}$ yields $\omega \wedge \omega_{i}^{\mathcal{R}} \wedge \omega'$.  Similar terms arise when we focus our attention to $\omega_{i}^{\mathcal{L}}$ when $i$ is even and $1<i<k$.  If we sum over all choices of $i, k, \underline{q}, \tau_a^i$, etc, we obtain all terms in $D(\underline{\omega})$.  Thus we have shown that 
$$(d_{DR} \circ It^A)(\underline{\omega}) = (It^A \circ D)(\underline{\omega}).$$
\end{proof}\end{prop}
\subsection{The two-dimensional case}
For the two dimensional curved case, we can follow the transition from the non-curved 1-d case to the non-curved 2-d case with one small addition.  Although one can proceed by collapsing the vertical left and right boundaries of our squares, and work on bigons, we will keep the square un-identified and so we need to account for parallel transport along the vertical paths moving from one zigzag to the next.  See Figure \ref{fig:It^A(exp(B))} on page \pageref{fig:It^A(exp(B))} for the idea.
\begin{definition}The curved rectangular zigzag Hochschild complex, $CH_{Rec}^{ZZ}(\A)$, has the same underlying vector space as in Definition \ref{def rect zig zag} with differential $D = \star + b + c + \nabla$.  Here, $\star$ and $b$ are the same as in the non-curved $CH_{Rec}^{ZZ}(\A)$ and we use the two-dimensional analog of our $\nabla$ and $c$ defined in the curved $CH^{ZZ}(\A)$ where now $c$ may also add a curvature term, $R$, at a vertical path on the left of the square.  In particular
\begin{align*}
&D((\underline{x}^0_{k_0, n}) \otimes \ldots \otimes (\underline{x}^{m+1}_{k_{m+1}, n}))\\
&:= \sum\limits_{r=1}^{m-1}(-1)^{m+r}((\underline{x}_{k_1, n}) \otimes \ldots \otimes (\underline{x}_{k_r, n}) \star (\underline{x}_{k_{r+1}, n}) \otimes \ldots \otimes (\underline{x}_{k_m, n}))\\
&+\sum\limits_{j=1}^m \sum\limits_{p=1}^{n-1} (-1)^{n+m+p}\ldots \otimes b_p(\underline{x}_{k_i, n})  \otimes \ldots \otimes b_p(\underline{x}_{k_j, n}) \otimes \ldots \\
&+\sum\limits_{j=1}^m \sum\limits_{p=1}^{n-1}(-1)^{n + \beta_{j,p}} \ldots \otimes \nabla_p(\underline{x}_{k_j, n}) \otimes \ldots \\
&+c((\underline{x}_{k_1, n}) \otimes \ldots \otimes (\underline{x}_{k_m, n}))
\end{align*} 
where \begin{align*}
c((\underline{x}_{k_1, n}) \otimes \ldots \otimes (\underline{x}_{k_m, n}))&= 
\sum\limits_{j=1}^m (-1)^{m} \sum\limits_{\sigma \in S_{n,1}} \ldots \otimes c_{\sigma}^1(\underline{x}_{k_i, n}) \otimes \ldots \otimes c_{\sigma}^{R}(\underline{x}_{k_j, n}) \otimes \ldots\\
&+ \sum\limits_{j=1}^{m-1} (-1)^{m+j} \ldots \otimes (\underline{x}_{k_j, n}) \otimes R \otimes (\underline{x}_{k_{j+1}, n}) \otimes \ldots.
\end{align*}
Here we used $c_{\sigma}^R$ to represent the usual component of our differential, $c$, which inserts an $R$ into one zig or zag and inserts $1$'s everywhere else in that new column.  The $c_{\sigma}^1$ mimics $c^R$ except that it only inserts $1$'s along the entire column.  Note also that the sign on the first line is only an $m$ since $c^R$ will be a sum of terms inserting $R$ between different columns, each term having an additional sign of $(-1)^{n + \sigma(n+1)}$ just as in Definition \ref{def of curved CH^ZZ}.
\end{definition}

For the sake of having a complete figure without all of the $A$'s inserted for the moment, we recall the following figure to work through the definition of our Iterated Integral in this case.
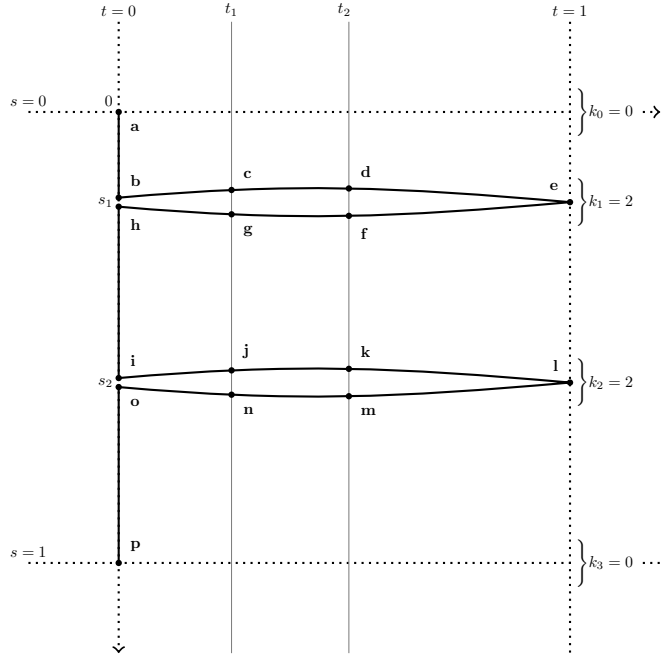
\begin{equation*}
\begin{tikzpicture}[thick,scale=0.6, every node/.style={transform shape}]
%help lines grid
%\draw[help lines] (0.1,0.1) grid (9.9,9.9);

% s axis
\draw[thick][dotted][->] (0,12) -- (0,-2);

% t axis
\draw[thick][dotted][->] (-2,10) -- (12,10);

%label axes
\node [above left] at (0,10) {0};
\node[above] at (0,12) {$t=0$};
\draw[thick][dotted] (10,12) -- (10,-2);
\node[above] at(10,12) {$t=1$};
\draw[thick][dotted] (-2,0) -- (12,0);
\node[above] at(-2,0) {$s=1$};
\node[above] at (-2,10) {$s=0$};

% s=0 zigzags
\draw[name path=zigzag] (0,10)
% s=s_1 zigzazags
-- (0,8.1) to [out=5, in=175] (10,8) to [out=185, in=-5](0,7.9)

%s=s_2 zigzags
-- (0,4.1) to [out=5, in=175] (10,4) to [out=185, in=-5](0,3.9)

%s=1 zigzags
-- (0,0);% to [out=5, in=175] (10,0) to [out=185, in=-5](0,-0.1);

%label si's
\node[left] at (0,8) {$s_1$};
\node[left] at (0,4) {$s_2$};

%label ki's
\node[right, fill=white] at (10,10) { $\Bigg\} k_0 = 0$};
\node[right, fill=white] at (10,8) { $\Bigg\} k_1 = 2$};
\node[right, fill=white] at (10,4) { $\Bigg\} k_2 = 2$};
\node[right, fill=white] at (10,0) { $\Bigg\} k_3 = 0$};

%t-lines

%t1
\draw[name path=t1, thin, gray] (2.5,-2) -- (2.5,12);
\node[above] at (2.5,12) {$t_1$};

%t2
\draw[name path=t2, thin, gray] (5.1,-2) -- (5.1,12);
\node[above] at (5,12) {$t_2$};

%t3
%\draw[name path=t3, thin, gray] (7.5,-2) -- (7.5,12);
%\node[above] at (7.5,12) {$t_3$};

\fill [name intersections={of=t1 and zigzag, by={a1, a2, a3, a4}}]
        \foreach \y in {1, 2, ...,4}
        {(a\y) circle (2pt)};
        \node[label=above right:$\bold{c}$] at (a1) {};
        \node[label=above right:$\bold{j}$] at (a3) {};
        \node[label=below right:$\bold{n}$] at (a4) {};
        \node[label=below right:$\bold{g}$] at (a2) {};

% \fill [name intersections={of=t3 and zigzag, by={a1, a2, a3, a4, a5, a6, a7, a8, a9, a10, a11, a12, a13, a14}}]
   %     \foreach \y in {1, 2, ..., 14}
      %  {(a\y) circle (2pt)};
        
\fill [name intersections={of=t2 and zigzag, by={a1, a2, a3, a4}}]
        \foreach \y in {1, 2, ...,4}
        {(a\y) circle (2pt)};
        \node[label=above right:$\bold{d}$] at (a1) {};
        \node[label=below right:$\bold{f}$] at (a2) {};
        \node[label=above right:$\bold{k}$] at (a3) {};
        \node[label=below right:$\bold{m}$] at (a4) {};

\fill (0,10) circle (2pt);
\node[label=below right:$\bold{a}$] at (0,10) {};
\fill (0,8.1) circle (2pt);
\node[label=above right:$\bold{b}$] at (0,8.1) {};
\fill (10,8) circle (2pt);
\node[label=above left:$\bold{e}$] at (10,8) {};
\fill (0,7.9) circle (2pt);
\node[label=below right:$\bold{h}$] at (0,7.9) {};
\fill (0,4.1) circle (2pt);
\node[label=above right:$\bold{i}$] at (0,4.1) {};
\fill (10,4) circle (2pt);
\node[label=above left:$\bold{l}$] at (10,4) {};
\fill (0,3.9) circle (2pt);
\node[label=below right:$\bold{o}$] at (0,3.9) {};
\fill (0,0) circle (2pt);
\node[label=above right:$\bold{p}$] at (0,0) {};

\end{tikzpicture}
\end{equation*}
First, notice that we have 15 sections in which to insert an arbitrary number of $A$'s.  So we consider those choices $\underline{q} = (q_1, \ldots, q_{15})$ and using our previous notation we have $(\phi = 16)$-many forms on $M$, namely $a, \ldots, p$ and $\phi^A = 16 + \sum q_i$.  Next, it becomes important to pull-back our form $A$ along a 1-path rather than a 2-path so that there is a single 1-path along which parallel transport is performed.  So for a 2-path $\Gamma(s,t)$ we mean define the 1-paths $\gamma_s(t) = \gamma^t(s) = \Gamma(t,s)$.  Note that our fiber $E$ will be quite cumbersome to write down in general.  There is certainly a formula, but the diagram gives that formula more easily than symbols.  For the particular diagram we are considering above, we have $t_1, t_2, s_1, s_2$ and 15 different choices of inserting $A$'s, labeled $q_1, \ldots q_{15}$.  We then have
\begin{align*}
&E \subset  (\prod \Delta^{q_i} ) \times \Delta^n \times \Delta^m  \quad \text{where}\\
&E= \{(\underline{\sigma}_1, \underline{\tau}_2, \underline{\tau}_3, \ldots, \underline{\tau}_7, \underline{\sigma}_8, \underline{\tau}_9, \underline{\tau}_{10}, \ldots, \underline{\tau}_{14}, \underline{\sigma}_{15}, t_1, t_2,  s_1,s_2),)|\\
&0 \le t_1 \le t_2 \le 1 \quad \text{and} \quad 0 \le s_2 \le s_2 \le 1  \\
& 0 \le \sigma_{(1,1)} \le \ldots \le \sigma_{(1,q_1)} \le s_1  \\
& 0 \le \tau_{(2,1)} \le \ldots \le  \tau_{(2,q_2)} \le t_1  \le \tau_{(3,1)} \le \ldots \le  \tau_{(3,q_3)} \le t_2  \le \tau_{(4,1)} \le \ldots \le  \tau_{(4,q_4)} \le 1  \\
& 0 \le \tau_{(7,q_7)} \le \ldots \le  \tau_{(7,1)} \le t_1  \le \tau_{(6,q_6)} \le \ldots \le  \tau_{(6,1)} \le t_2  \le \tau_{(5,q_5)} \le \ldots \le  \tau_{(5,1)} \le 1  \\
 & s_1 \le \sigma_{(8,1)} \le \ldots \le  \sigma_{(8,q_8)} \le s_2  \\
& 0 \le \tau_{(9,1)} \le \ldots \le  \tau_{(9,q_{9})} \le t_1  \le \tau_{(10,1)} \le \ldots \le  \tau_{(10,q_{10})} \le t_2  \le \tau_{(11,1)} \le \ldots \le  \tau_{(11,q_{11})} \le 1  \\
& 0 \le \tau_{(14,q_{14})} \le \ldots \le  \tau_{(14,1)} \le t_1  \le \tau_{(13,q_{13})} \le \ldots \le  \tau_{(13,1)} \le t_2  \le \tau_{(12,q_{12})} \le \ldots \le  \tau_{(12,1)} \le 1  \\
& s_2 \le \sigma_{(15,1)} \le \ldots \le  \sigma_{(15,q_{15})} \le 1 \}  
\end{align*}
using the convention $\underline{\sigma}_i = ( \sigma_{(1,1)}, \ldots, \sigma_{(1, q_i)})$ and $\underline{\tau}_i = ( \tau_{(1,1)}, \ldots, \tau_{(1, q_i)})$.  Observe that to each choice of $n, m$, $\underline{k}$, and $\underline{q}$, where $\underline{k}:=(k_0, \ldots, k_{m+1})$ and $\underline{q}:= (q_1, \ldots, q_{(n+1)(k_0 + \ldots + k_{m+1}) + (m+1)})$, we have a uniquely determined zigzag diagram with $A$'s inserted.  Our fiber $E$ is determined by these choices as well.  With this in mind, the evaluation map for the above element can be written:
\begin{align*}
&BM \times E \xrightarrow{ev} M^{\phi^A}\\
&( \gamma, \underline{\sigma}_1, \underline{\tau}_2, \underline{\tau}_3, \ldots,  \underline{\tau}_7, \underline{\sigma}_8,  \underline{\tau}_9,\ldots, \underline{\tau}_{14}, \underline{\sigma}_{15}, t_1,t_2, s_1,s_2) \\
\mapsto & ( \Gamma(0,0), ( \gamma^0(\sigma_{(1,i)}))_{i=1}^{q_1}, \Gamma(0, s_1), ( \gamma_{s_1}(\tau_{(2,i)}))_{i=1}^{q_2}, \Gamma(t_1, s_1), ( \gamma_{s_1}(\tau_{(3,i)}))_{i=1}^{q_3},\\
& \ldots, ( \gamma_{s_1}(\tau_{(5,i)}))_{i=1}^{q_5}, \ldots , ( \gamma_{s_2}(\tau_{(14,i)}))_{i=1}^{q_{14}}, \Gamma(0, s_2), (\gamma^0( \sigma_{(15,i)}))_{i=1}^{q_{15}}, \Gamma(0,1) )
\end{align*}
\begin{definition}
We define the curved iterated integral $It:CH^{ZZ}_{Rec}( \OmegaMat(M)^{\phi}) \to \OmegaMat(BM)$ by $$It: = \sum\limits_{(\underline{q})} (-1)^{\sum |\omega_j| Q_j} It_{(\underline{q})}^A$$ where the component of the iterated integral $It_{(\underline{q})}^A: \OmegaMat(M)^{\otimes  \phi} \to \OmegaMat(BM)$ is given by the diagram
\begin{equation*}
\xymatrix{\OmegaMat(BM \times E)    \ar[d]^{\int_{E}} &\ar[l]_{ev^*}\OmegaMat(M^{\phi^A}) &\ar[l]  \OmegaMat(M)^{\otimes {\phi^A}}  & \ar[llld]^{It^A_{(\underline{q})}}  \ar[l]_{ Ins^A_{(\underline{q})}}  \OmegaMat(M)^{\otimes {\phi}} \\
\OmegaMat(B(M)) &  & &}
\end{equation*}
\end{definition}
Since our differential is simply a mix of terms from the 2-d non-curved case and the 1-d curved case, we can combine all of those arguments to obtain the following theorem.
\begin{theorem}
The curved iterated integral $It: CH_{Rec}^{ZZ}(\OmegaMat(M)) \to \OmegaMat(B(M))$ is a chain map.
\end{theorem}
\section{Holonomy}
We first remark on elements in $CH^{ZZ}(\OmegaMat(M))$ which map to 1-dimensional holonomy.  Given a 1-form $A \in \OmegaMat^1(M)$ we denote by $P_{\gamma}(t)$ the parallel transport along a path $\gamma$ from $0$ to $t$.  By the construction of our zigzag Hochschild complex and its shuffle product we note that 
$$P_{\gamma}(t) = \sum\limits_{n \ge 0} It(\tilde{A})^{\odot n}$$
where $\tilde{A}:= (1 \otimes A \otimes 1 \otimes 1 \otimes 1)$ has $n=1$ and $k=2$ and each  $\tilde{A}^{\odot n} \in CH^{ZZ}(\OmegaMat(M))$.

Finally, we show that we have elements in the completion of $CH^{ZZ}_{Rec}(\OmegaMat(M))$ which map to 2-dimensional holonomy.  In the remainder of this chapter we restrict our space $BM$ to the subspace of $\{ \Sigma:[0,1]^2 \to M \}$ where for each $\Sigma(t,s): [0,1]^2 \to M$ there exist $x, y \in M$ so that $\Sigma(0,s) = x$ and $\Sigma(1,s) = y$ for all $ s \in [0,1]$ (i.e. ``bigons'').  Let $B \in \OmegaMat^2(M)$ be a matrix-valued 2-form.  We define an element $\exp(B) = \sum\limits_{m \ge 0} B^{\tilde{\odot}m}$, so that each $B^{\tilde{\odot}m} \in CH^{ZZ}_{Rec}(\OmegaMat(M))$, as in the completed zigzag Hochschild complex, and show that the curved iterated integral of this element is the well-known 2-holonomy as defined in \cite{BaSc}, \cite{PM1}, \cite{PM2}, and \cite{ScWa}.  Define the element $\exp(B):= \sum\limits_{m \ge 0} B^{\tilde{\odot}m}$ via Figure \ref{fig:expB} below.  Note that although $It^A(\exp(B)):= \sum\limits_{m\ge 0} It^A(B^{\tilde{\odot} m})$ is an infinite sum of forms in $\OmegaMat(BM)$, a simple boundedness condition on $B$ guarantees that this sum converges in $\OmegaMat(BM)$; compare the argument for the well-definedness of $It^A$ above Proposition \ref{It^A chain and algebra}.
Following Baez, Martins, Picken, Schreiber, and Waldorf (\cite{PM1}, \cite{PM2}, \cite{ScWa}, and \cite{BaSc}), we show that $It^A(\exp(B))$ solves the differential equation which governs 2-holonomy.  We do so using more familiar notation, suppressing the evaluation pullback notation used previously in this paper.   
\begin{prop}
Let $A \in \OmegaMat^1(M)$ and $B \in \OmegaMat^2(M)$, then  $$\frac{\partial}{\partial s} (It^A(\exp(B(t,s)))) = It^A(\exp(B(t,s))) \wedge \int\limits_0^1 P_{(t',s)}(B(t',s))P_{(t',s)}^{-1} dt.$$where $P_{(t,s)}$ is the 1-holonomy obtained from $A$ via the path $\restr{\gamma_s }{[0,t]} \circ \restr{\gamma^0}{[0,s]}$.
\begin{proof}
Let $tra(B)(a,b):=  P_{(a,b)}(B(a,b))P^{-1}_{(a,b)}$, then 
\begin{align*}
&\frac{\partial}{\partial s} (It^A (\exp(B(t,s)))) =  \sum\limits_{m \ge 0} \frac{\partial}{\partial s} It^A( B(t,s)^{\tilde{\odot}m})\\
= &\sum\limits_{m \ge 0} \sum\limits_{\sigma \in S_m} sgn(\sigma) \frac{\partial}{\partial s}  \int\limits_{\Delta^{m} \times \Delta^{m}} tra(B)(t_{\sigma^{-1}(1)},s_1)  \\ \quad \quad \quad \dots &  tra(B)(t_{\sigma^{-1}(m)},s_m) dt_1 \dots dt_m ds_1 \dots ds_m \\
= &\sum\limits_{m \ge 0} \sum\limits_{\sigma \in S_m}  sgn(\sigma)  \int\limits_{\Delta^{m} \times \Delta^{m-1}}   tra(B)(t_{\sigma^{-1}(1)},s_1)  \\ \quad \quad \quad \dots &tra(B)(t_{\sigma^{-1}(m-1)},s_{m-1}) tra(B)(t_{\sigma^{-1}(m)},s)dt_1 \dots dt_m ds_1 \dots ds_{m-1}\\
= &\sum\limits_{m \ge 0}  \sum\limits_{\sigma \in S_{m-1}} sgn(\sigma) \int\limits_{\Delta^{m} \times \Delta^{m-1}}  tra(B)(t_{\sigma^{-1}(1)},s_1)\\  \quad \quad \quad  \dots &tra(B)(t_{\sigma^{-1}(m-1)},s_{m-1}) dt_1 \dots dt_{m-1} ds_1 \dots ds_{m-1}\\  \wedge &\int\limits_0^t  tra(B)(t',s) dt'\\
 = &It^A(\exp(B)) \wedge \int\limits_0^t P_{(t',s)}(B(t',s))P^{-1}_{(t',s)} dt' 
\end{align*}  \end{proof}
\end{prop}

\begin{figure}[H]
\resizebox{9cm}{!}{\begin{tabular}{ccc}
&$1$  +
&\resizebox{5cm}{!}{\begin{tikzpicture}[scale= 0.8, baseline={([yshift=-.5ex]current bounding box.center)},vertex/.style={anchor=base,circle,fill=black!25,minimum size=18pt,inner sep=2pt}]
%help lines grid
%\draw[help lines] (0.1,0.1) grid (9.9,9.9);
% s axis
\draw[thick][dotted][->] (0,12) -- (0,-2);
\node [left, lightgray] at (0,-2) {$s$};
% t axis
\draw[thick, dotted, lightgray][->] (-2,10) -- (12,10);
\node [above, lightgray] at (12,10) {$t$};
%label axes
\node [above left] at (0,10) {0};
\node[above, lightgray] at (0,12) {$t=0$};
\draw[thick][dotted] (10,12) -- (10,-2);
\node[above, lightgray] at(10,12) {$t=1$};
\draw[thick][dotted] (-2,0) -- (12,0);
\node[above,lightgray] at(-2,0) {$s=1$};
\node[above, lightgray] at (-2,10) {$s=0$};
% s=s_1 zigzags
\draw[name path=zigzag] (0,10) -- (0,5.1)  to [out=5, in=175]  (10,5) to [out=185, in=-5]    (0,4.9)  
%s=1 zigzags
--   (0,0);% to [out=5, in=175] (10,0) to [out=185, in=-5](0,-0.1);
%label si's
\node[left, lightgray] at (0,8) {$s_1$};
%label ki's
%\node[right, fill=white] at (10,10) { $\Bigg\} k_0 = 4$};
%\node[right, fill=white] at (10,8) { $\Bigg\} k_{s_1} = 2$};
%\node[right, fill=white] at (10,4) { $\Bigg\} k_{s_2} = 6$};
%\node[right, fill=white] at (10,0) { $\Bigg\} k_1 = 2$};
%t-lines
%t1
\draw[thin, lightgray, name path=t1] (6,-2) -- (6,12);
\node[above, lightgray] at (6,12) {$t_1$};
\fill [name intersections={of=t1 and zigzag, by={a,b}}]
        (a) circle (2pt)
        (b) circle (2pt);
\node[above left] at (a) {\huge$B$};
\node[below left = -0.02cm] at (b) {$1$};
\end{tikzpicture}}\\
$+$
&\resizebox{5cm}{!}{\begin{tikzpicture}[scale= 0.8, baseline={([yshift=-.5ex]current bounding box.center)},vertex/.style={anchor=base,circle,fill=black!25,minimum size=18pt,inner sep=2pt}]
%help lines grid
%\draw[help lines] (0.1,0.1) grid (9.9,9.9);
% s axis
\draw[thick][dotted][->] (0,12) -- (0,-2);
\node [left, lightgray] at (0,-2) {$s$};
% t axis
\draw[thick, dotted, lightgray][->] (-2,10) -- (12,10);
\node [above, lightgray] at (12,10) {$t$};
%label axes
\node [above left] at (0,10) {0};
\node[above, lightgray] at (0,12) {$t=0$};
\draw[thick][dotted] (10,12) -- (10,-2);
\node[above, lightgray] at(10,12) {$t=1$};
\draw[thick][dotted] (-2,0) -- (12,0);
\node[above,lightgray] at(-2,0) {$s=1$};
\node[above, lightgray] at (-2,10) {$s=0$};
% s=0 zigzags
\draw[name path=zigzag] (0,10) -- (0,8.1)  to [out=5, in=175]  (10,8) to [out=185, in=-5]    (0,7.9)  
%s=s_2 zigzags
--   (0,6.1) to [out=5, in=175]  (10,6) to [out=185, in=-5]  (0,5.9)
%--(0, 4.2) to  [out=10, in=170] (10,4.2) to [out=175, in=7](0,4.1)
%to [out=3, in=177] (10,4) to [out=183, in=-3] (0,3.9)
%to [out=-7, in=187] (10,3.8) to [out=190, in=-10](0,3.8)
%s=1 zigzags
--   (0,0);% to [out=5, in=175] (10,0) to [out=185, in=-5](0,-0.1);
%label si's
\node[left, lightgray] at (0,8) {$s_1$};
\node[left, lightgray] at (0,6) {$s_2$};
%label ki's
%\node[right, fill=white] at (10,10) { $\Bigg\} k_0 = 4$};
%\node[right, fill=white] at (10,8) { $\Bigg\} k_{s_1} = 2$};
%\node[right, fill=white] at (10,4) { $\Bigg\} k_{s_2} = 6$};
%\node[right, fill=white] at (10,0) { $\Bigg\} k_1 = 2$};
%t-lines
%t1
\draw[thin, lightgray, name path=t1] (3,-2) -- (3,12);
\node[above, lightgray] at (3,12) {$t_1$};
%t2
\draw[thin, lightgray, name path=t2] (7,-2) -- (7,12);
\node[above, lightgray] at (7,12) {$t_2$};
\fill [name intersections={of=t1 and zigzag, by={a,b, c, d}}]
        (a) circle (2pt)
        (b) circle (2pt)
        (c) circle (2pt)
        (d) circle (2pt);
\node[above left] at (a) {\huge$B$};
\node[below left = -0.02cm] at (b) {$1$};
\node[above left = -0.02cm] at (c) {$1$};
\node[below left = -0.02cm] at (d) {$1$};
\fill [name intersections={of=t2 and zigzag, by={a,b, c, d}}]
        (a) circle (2pt)
        (b) circle (2pt)
        (c) circle (2pt)
        (d) circle (2pt);       
\node[above left = -0.02cm] at (a) {$1$};
\node[below left = -0.02cm] at (b) {$1$};
\node[above left] at (c) {\huge$B$};
\node[below left = -0.02cm] at (d) {$1$};
\end{tikzpicture}}

$-$
&\resizebox{5cm}{!}{\begin{tikzpicture}[scale= 0.8, baseline={([yshift=-.5ex]current bounding box.center)},vertex/.style={anchor=base,circle,fill=black!25,minimum size=18pt,inner sep=2pt}]
%help lines grid
%\draw[help lines] (0.1,0.1) grid (9.9,9.9);

% s axis
\draw[thick][dotted][->] (0,12) -- (0,-2);
\node [left, lightgray] at (0,-2) {$s$};
% t axis
\draw[thick, dotted, lightgray][->] (-2,10) -- (12,10);
\node [above, lightgray] at (12,10) {$t$};
%label axes
\node [above left] at (0,10) {0};
\node[above, lightgray] at (0,12) {$t=0$};
\draw[thick][dotted] (10,12) -- (10,-2);
\node[above, lightgray] at(10,12) {$t=1$};
\draw[thick][dotted] (-2,0) -- (12,0);
\node[above,lightgray] at(-2,0) {$s=1$};
\node[above, lightgray] at (-2,10) {$s=0$};

% s=0 zigzags
\draw[name path=zigzag] (0,10) -- (0,8.1)  to [out=5, in=175]  (10,8) to [out=185, in=-5]    (0,7.9)

%s=s_2 zigzags
--   (0,6.1) to [out=5, in=175]  (10,6) to [out=185, in=-5]  (0,5.9)
%--(0, 4.2) to  [out=10, in=170] (10,4.2) to [out=175, in=7](0,4.1)
%to [out=3, in=177] (10,4) to [out=183, in=-3] (0,3.9)
%to [out=-7, in=187] (10,3.8) to [out=190, in=-10](0,3.8)

%s=1 zigzags
--   (0,0);% to [out=5, in=175] (10,0) to [out=185, in=-5](0,-0.1);

%label si's
\node[left, lightgray] at (0,8) {$s_1$};
\node[left, lightgray] at (0,6) {$s_2$};

%label ki's
%\node[right, fill=white] at (10,10) { $\Bigg\} k_0 = 4$};
%\node[right, fill=white] at (10,8) { $\Bigg\} k_{s_1} = 2$};
%\node[right, fill=white] at (10,4) { $\Bigg\} k_{s_2} = 6$};
%\node[right, fill=white] at (10,0) { $\Bigg\} k_1 = 2$};

%t-lines

%t1
\draw[thin, lightgray, name path=t1] (3,-2) -- (3,12);
\node[above, lightgray] at (3,12) {$t_1$};

%t2
\draw[thin, lightgray, name path=t2] (7,-2) -- (7,12);
\node[above, lightgray] at (7,12) {$t_2$};

\fill [name intersections={of=t1 and zigzag, by={a,b, c, d}}]
        (a) circle (2pt)
        (b) circle (2pt)
        (c) circle (2pt)
        (d) circle (2pt);

         \node[above left = -0.02cm] at (a) {$1$};
\node[below left = -0.02cm] at (b) {$1$};
\node[above left] at (c) {\huge$B$};
\node[below left = -0.02cm] at (d) {$1$};

\fill [name intersections={of=t2 and zigzag, by={a,b, c, d}}]
        (a) circle (2pt)
        (b) circle (2pt)
        (c) circle (2pt)
        (d) circle (2pt);

\node[above left] at (a) {\huge$B$};
\node[below left = -0.02cm] at (b) {$1$};
\node[above left = -0.02cm] at (c) {$1$};
\node[below left = -0.02cm] at (d) {$1$};\

\end{tikzpicture}}\\
+
&\resizebox{5cm}{!}{\begin{tikzpicture}[scale= 0.8, baseline={([yshift=-.5ex]current bounding box.center)},vertex/.style={anchor=base,circle,fill=black!25,minimum size=18pt,inner sep=2pt}]
%help lines grid
%\draw[help lines] (0.1,0.1) grid (9.9,9.9);

% s axis
\draw[thick][dotted][->] (0,12) -- (0,-2);
\node [left, lightgray] at (0,-2) {$s$};
% t axis
\draw[thick, dotted, lightgray][->] (-2,10) -- (12,10);
\node [above, lightgray] at (12,10) {$t$};
%label axes
\node [above left] at (0,10) {0};
\node[above, lightgray] at (0,12) {$t=0$};
\draw[thick][dotted] (10,12) -- (10,-2);
\node[above, lightgray] at(10,12) {$t=1$};
\draw[thick][dotted] (-2,0) -- (12,0);
\node[above,lightgray] at(-2,0) {$s=1$};
\node[above, lightgray] at (-2,10) {$s=0$};

% s=0 zigzags
\draw[name path=zigzag] (0,10) -- (0,8.1)  to [out=5, in=175]  (10,8) to [out=185, in=-5]    (0,7.9)

%s=s_2 zigzags
--   (0,6.1) to [out=5, in=175]  (10,6) to [out=185, in=-5]  (0,5.9)
%--(0, 4.2) to  [out=10, in=170] (10,4.2) to [out=175, in=7](0,4.1)
%to [out=3, in=177] (10,4) to [out=183, in=-3] (0,3.9)
%to [out=-7, in=187] (10,3.8) to [out=190, in=-10](0,3.8)

%s=s_3 zigzags
--   (0,2.1) to [out=5, in=175]  (10,2) to [out=185, in=-5]  (0,1.9)

%s=1 zigzags
--   (0,0);% to [out=5, in=175] (10,0) to [out=185, in=-5](0,-0.1);

%label si's
\node[left, lightgray] at (0,8) {$s_1$};
\node[left, lightgray] at (0,6) {$s_2$};
\node[left, lightgray] at (0,2) {$s_3$};

%label ki's
%\node[right, fill=white] at (10,10) { $\Bigg\} k_0 = 4$};
%\node[right, fill=white] at (10,8) { $\Bigg\} k_{s_1} = 2$};
%\node[right, fill=white] at (10,4) { $\Bigg\} k_{s_2} = 6$};
%\node[right, fill=white] at (10,0) { $\Bigg\} k_1 = 2$};

%t-lines

%t1
\draw[thin, lightgray, name path=t1] (2.5,-2) -- (2.5,12);
\node[above, lightgray] at (2.5,12) {$t_1$};

%t2
\draw[thin, lightgray, name path=t2] (5,-2) -- (5,12);
\node[above, lightgray] at (5,12) {$t_2$};

%t3
\draw[thin, lightgray, name path=t3] (7.5,-2) -- (7.5,12);
\node[above, lightgray] at (7.5,12) {$t_3$};

\fill [name intersections={of=t1 and zigzag, by={a,b, c, d, e, f}}]
        (a) circle (2pt)
        (b) circle (2pt)
        (c) circle (2pt)
        (d) circle (2pt)
        (e) circle (2pt)
         (f) circle (2pt);
\node[above left] at (a) {\huge$B$};
\node[below left = -0.02cm] at (b) {$1$};
\node[above left = -0.02cm] at (c) {$1$};
\node[below left = -0.02cm] at (d) {$1$};
\node[above left = -0.02cm] at (e) {$1$};
\node[below left = -0.02cm] at (f) {$1$};

\fill [name intersections={of=t2 and zigzag, by={a,b, c, d, e, f}}]
        (a) circle (2pt)
        (b) circle (2pt)
        (c) circle (2pt)
        (d) circle (2pt)
        (e) circle (2pt)
         (f) circle (2pt);
         
         \node[above left = -0.02cm] at (a) {$1$};
\node[below left = -0.02cm] at (b) {$1$};
\node[above left] at (c) {\huge$B$};
\node[below left = -0.02cm] at (d) {$1$};
\node[above left = -0.02cm] at (e) {$1$};
\node[below left = -0.02cm] at (f) {$1$};

\fill [name intersections={of=t3 and zigzag, by={a,b, c, d, e, f}}]
        (a) circle (2pt)
        (b) circle (2pt)
        (c) circle (2pt)
        (d) circle (2pt)
        (e) circle (2pt)
         (f) circle (2pt);
         
\node[above left = -0.02cm] at (a) {$1$};
\node[below left = -0.02cm] at (b) {$1$};
\node[above left = -0.02cm] at (c) {$1$};
\node[below left = -0.02cm] at (d) {$1$};
\node[above left] at (e) {\huge$B$};
\node[below left = -0.02cm] at (f) {$1$};

\end{tikzpicture}}

$-$
&\resizebox{5cm}{!}{\begin{tikzpicture}[scale= 0.8, baseline={([yshift=-.5ex]current bounding box.center)},vertex/.style={anchor=base,circle,fill=black!25,minimum size=18pt,inner sep=2pt}]
%help lines grid
%\draw[help lines] (0.1,0.1) grid (9.9,9.9);

% s axis
\draw[thick][dotted][->] (0,12) -- (0,-2);
\node [left, lightgray] at (0,-2) {$s$};
% t axis
\draw[thick, dotted, lightgray][->] (-2,10) -- (12,10);
\node [above, lightgray] at (12,10) {$t$};
%label axes
\node [above left] at (0,10) {0};
\node[above, lightgray] at (0,12) {$t=0$};
\draw[thick][dotted] (10,12) -- (10,-2);
\node[above, lightgray] at(10,12) {$t=1$};
\draw[thick][dotted] (-2,0) -- (12,0);
\node[above,lightgray] at(-2,0) {$s=1$};
\node[above, lightgray] at (-2,10) {$s=0$};

% s=0 zigzags
\draw[name path=zigzag] (0,10) -- (0,8.1)  to [out=5, in=175]  (10,8) to [out=185, in=-5]    (0,7.9)

%s=s_2 zigzags
--   (0,6.1) to [out=5, in=175]  (10,6) to [out=185, in=-5]  (0,5.9)
%--(0, 4.2) to  [out=10, in=170] (10,4.2) to [out=175, in=7](0,4.1)
%to [out=3, in=177] (10,4) to [out=183, in=-3] (0,3.9)
%to [out=-7, in=187] (10,3.8) to [out=190, in=-10](0,3.8)

%s=s_3 zigzags
--   (0,2.1) to [out=5, in=175]  (10,2) to [out=185, in=-5]  (0,1.9)

%s=1 zigzags
--   (0,0);% to [out=5, in=175] (10,0) to [out=185, in=-5](0,-0.1);

%label si's
\node[left, lightgray] at (0,8) {$s_1$};
\node[left, lightgray] at (0,6) {$s_2$};
\node[left, lightgray] at (0,2) {$s_3$};

%label ki's
%\node[right, fill=white] at (10,10) { $\Bigg\} k_0 = 4$};
%\node[right, fill=white] at (10,8) { $\Bigg\} k_{s_1} = 2$};
%\node[right, fill=white] at (10,4) { $\Bigg\} k_{s_2} = 6$};
%\node[right, fill=white] at (10,0) { $\Bigg\} k_1 = 2$};

%t-lines

%t1
\draw[thin, lightgray, name path=t1] (2.5,-2) -- (2.5,12);
\node[above, lightgray] at (2.5,12) {$t_1$};

%t2
\draw[thin, lightgray, name path=t2] (5,-2) -- (5,12);
\node[above, lightgray] at (5,12) {$t_2$};

%t3
\draw[thin, lightgray, name path=t3] (7.5,-2) -- (7.5,12);
\node[above, lightgray] at (7.5,12) {$t_3$};

\fill [name intersections={of=t1 and zigzag, by={a,b, c, d, e, f}}]
        (a) circle (2pt)
        (b) circle (2pt)
        (c) circle (2pt)
        (d) circle (2pt)
        (e) circle (2pt)
         (f) circle (2pt);

         \node[above left = -0.02cm] at (a) {$1$};
\node[below left = -0.02cm] at (b) {$1$};
\node[above left] at (c) {\huge$B$};
\node[below left = -0.02cm] at (d) {$1$};
\node[above left = -0.02cm] at (e) {$1$};
\node[below left = -0.02cm] at (f) {$1$};

\fill [name intersections={of=t2 and zigzag, by={a,b, c, d, e, f}}]
        (a) circle (2pt)
        (b) circle (2pt)
        (c) circle (2pt)
        (d) circle (2pt)
        (e) circle (2pt)
         (f) circle (2pt);
      \node[above left] at (a) {\huge$B$};
\node[below left = -0.02cm] at (b) {$1$};
\node[above left = -0.02cm] at (c) {$1$};
\node[below left = -0.02cm] at (d) {$1$};
\node[above left = -0.02cm] at (e) {$1$};
\node[below left = -0.02cm] at (f) {$1$};

\fill [name intersections={of=t3 and zigzag, by={a,b, c, d, e, f}}]
        (a) circle (2pt)
        (b) circle (2pt)
        (c) circle (2pt)
        (d) circle (2pt)
        (e) circle (2pt)
         (f) circle (2pt);
         
\node[above left = -0.02cm] at (a) {$1$};
\node[below left = -0.02cm] at (b) {$1$};
\node[above left = -0.02cm] at (c) {$1$};
\node[below left = -0.02cm] at (d) {$1$};
\node[above left] at (e) {\huge$B$};
\node[below left = -0.02cm] at (f) {$1$};

\end{tikzpicture}}
$+ \ldots$
\end{tabular}}
\caption{Our element $\exp(B)$.  Here the signs come from the permutation governing where the $B$'s appear along the zigzags.}
\label{fig:expB}
\end{figure}
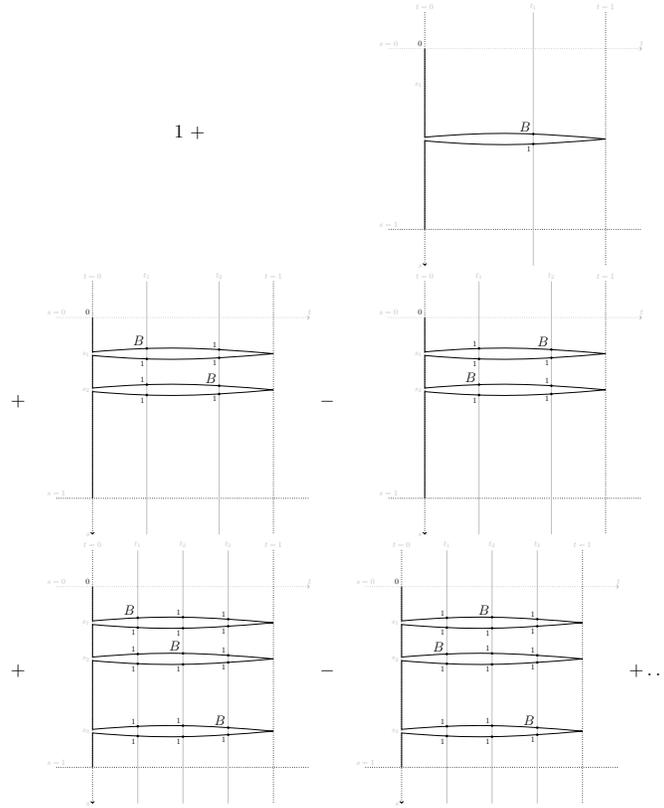

\begin{cor}
For a gerbe $\mathcal{G}$ with local data given by a 1-form $A \in \OmegaMat^1(U)$ and $B \in \OmegaMat^2(U)$ on an open set $U \subset M$ of a manifold $M$, the form $It^A(\exp(B)) \in \OmegaMat^0(BU)$ represents the local 2-holonomy on $U$.
\end{cor}

\begin{remark}
If we did not restrict our $BM$ to true bigons and instead continued to use squares, then our $\exp(B)$ would be off by $P_{(0,s)}^{-1}$.  One could simply then use the expression $It^A(\exp(B)) \cdot P_{(0,s)}^{-1}$ to obtain 2-dimensional holonomy.
\end{remark}

\end{document}